\documentclass[11pt, leqno]{amsart}

\usepackage{amssymb,amsmath,amsthm}
\usepackage{mathrsfs}
\usepackage{enumitem}
\usepackage[all]{xy}
\usepackage{placeins}
\usepackage{tikz}
\usepackage[normalem]{ulem}

\usepackage[colorlinks=true, pdfstartview={FitV null},
 linkcolor=blue,citecolor=blue,urlcolor=blue]{hyperref}

\allowdisplaybreaks[2]
\renewcommand{\le}{\leqslant}
\renewcommand{\ge}{\geqslant}
\numberwithin{equation}{section}

\theoremstyle{plain}
\newtheorem{lemma}{Lemma}[section]
\newtheorem{proposition}[lemma]{Proposition}
\newtheorem{theorem}[lemma]{Theorem}
\newtheorem{corollary}[lemma]{Corollary}

\theoremstyle{definition}

\newtheorem{example}[lemma]{Example}

\newcommand{\Q}{\mathbb Q}
\newcommand{\Z}{\mathbb Z}
\newcommand{\g}{\mathfrak g}
\newcommand{\n}{\mathfrak n}
\newcommand{\Hom}{\operatorname{Hom}}
\newcommand{\one}{\mathbf 1}
\newcommand{\ep}{\varepsilon}
\newcommand{\te}{\widetilde e}
\newcommand{\tf}{\widetilde f}
\newcommand{\wt}{\operatorname{wt}}
\newcommand{\supp}{\operatorname{supp}}

\newcommand{\weyl}{\mathsf W}
\newcommand{\rR}{\mathrm R}
\newcommand{\bfa}{\mathbf a}

\newcommand{\bfi}{\mathbf i}
\newcommand{\bfj}{\mathbf j}
\newcommand{\bfk}{\mathbf k}
\newcommand{\bfv}{\mathbf v}

\newcommand{\qA}{\mathrm A}
\newcommand{\Gup}{\mathrm G^{\mathrm{up}}}
\newcommand{\tGup}{\widetilde{\mathrm G}^{\mathrm{up}}}

\newcommand{\PBW}{\mathrm{PBW}}
\newcommand{\STR}{\mathrm{STR}}
\newcommand{\tcP}{\widetilde{\mathsf P}}
\newcommand{\tcS}{\widetilde{\mathsf S}}
\newcommand{\pbwP}{\mathfrak P}
\newcommand{\strS}{\mathfrak S}
\newcommand{\gL}{\mathrm g^{\mathrm L}}
\newcommand{\gR}{\mathrm g^{\mathrm R}}

\newcommand{\gM}{\mathcal M}
\newcommand{\gN}{\mathcal N}
\newcommand{\LB}{\widetilde B}
\newcommand{\qmD}{\mathsf D}
\newcommand{\qmC}{\mathrm C}
\newcommand{\fz}{\mathsf c}
\newcommand{\eqfr}{{\,\equiv_{\mathrm{fr}}\,}}
\newcommand{\cdual}{\mathcal D}
\newcommand{\PD}{\varrho}

\title[Similarity maps and quantum twists]{Similarity maps and quantum twists on localized crystals}
\author[W.-S. Jung]{Woo-Seok Jung}
\thanks{The research of W.-S. Jung was supported by the National Research Foundation of Korea (NRF) grants funded by the Korea government (Ministry of Education and MSIT) (RS-2024-00451844 and RS-2023-00273425).}
\address[W.-S. Jung]{Department of Mathematics, University of Seoul, Seoul 02504, Korea}
\email{jungws@uos.ac.kr}
\subjclass[2020]{16T20, 17B37, 05E10}
\keywords{crystals, quantum groups, quantum twist automorphism, similarity maps, localized crystals}
\date{\today}

\begin{document}
\begin{abstract}

Kashiwara's similarity maps give embeddings between crystals associated with
different Cartan data. We extend the similarity maps to localized crystals
and establish compatibility with the crystal structures of Kashiwara and
Nakashima.
We describe the extended similarity maps in the Poincar\'e--Birkhoff--Witt
(PBW), string, and right and left $g$-vector parametrizations by the same
linear formula.
Using the four parametrizations, we prove that the extended similarity maps
intertwine the permutations induced by the quantum twist automorphisms of
Kimura and Oya.
As an application, we determine the periods of quantum twists on localized
crystals, modulo multiplication by frozen variables, for a family of Weyl group
elements in types $B_n$ and $C_n$.
\end{abstract}

\maketitle
\tableofcontents

\section*{Introduction}

Let $\g$ be a symmetrizable Kac--Moody algebra with negative part
$\n$, simple roots indexed by $I$, and Weyl group $\weyl$, and let
$q$ be an indeterminate.
The quantum unipotent coordinate ring $\qA_q(\n)$ contains, for each
$w\in\weyl$, a subalgebra $\qA_q(\n(w))$ quantizing the coordinate
ring of the unipotent subgroup associated with $w$ \cite{Kimura12}.
For the fundamental weights $\Lambda_i$, write
$\qmC_i^w=\qmD(w\Lambda_i,\Lambda_i)$ for the corresponding quantum
unipotent minors, indexed by $i\in I$.
The minors $\qmC_i^w$, $i\in I$, commute with one another up to powers of $q$
\cite[Corollary~6.18]{Kimura12}.
Adjoining $(\qmC_i^w)^{-1}$, for $i\in I$, gives the algebra
\[
\widetilde{\qA}_q(\n(w))
=\qA_q(\n(w))[(\qmC_i^w)^{-1}\mid i\in I].
\]
The algebra $\widetilde{\qA}_q(\n(w))$ is isomorphic to the quantum coordinate ring of the
unipotent cell associated with $w$ \cite[Theorem~4.13]{KiOy21}.

On $\widetilde{\qA}_q(\n(w))$, Kimura and Oya constructed the quantum twist
automorphism
\[
\eta_w\colon\widetilde{\qA}_q(\n(w))
\xrightarrow{\sim}\widetilde{\qA}_q(\n(w))
\]
\cite[Theorem~6.1]{KiOy21}.
At $q=1$, $\eta_w$ specializes to the automorphism of the coordinate
ring induced by the classical twist of the unipotent cell
\cite[Corollary~6.2]{KiOy21}.
In finite type, the classical twist introduced in \cite{BFZ96,BZ97}
is used in the Chamber Ansatz to recover factorization coordinates of
elements of the unipotent cell associated with $w$.
Kimura and Oya \cite[Theorem~6.1]{KiOy21} proved that $\eta_w$ permutes
the upper global basis $\tGup(w)$ of
$\widetilde{\qA}_q(\n(w))$
(see \cite{HK02,LusztigBook} for global bases).

The action of $\eta_w^{-1}$ on $\tGup(w)$ can be expressed in terms of
crystal elements.
Let $U_q^-(\g)$ be the negative part of the quantum enveloping
algebra of $\g$, and let $B(\infty)$ be the crystal of $U_q^-(\g)$.
The elements $b\in B(\infty)$ index the upper global basis elements
$\Gup(b)$ of $\qA_q(\n)$, and we set
\[
B(w)=\{b\in B(\infty)\mid
\Gup(b)\in\qA_q(\n(w))\}.
\]
Assume that every index in $I$ occurs in a reduced expression of $w$.
For $i\in I$, let $\fz_i^w\in B(w)$ satisfy
$\Gup(\fz_i^w)=\qmC_i^w$.
Following \cite[Definition~2.3]{JP25}, let $\LB(w)$ be the
localization of $B(w)$ obtained by inverting the elements
$\fz_i^w$, $i\in I$.
We write an element of $\LB(w)$ as a formal product
\[
x=b\cdot\prod_{i\in I}(\fz_i^w)^{a_i}
\qquad(b\in B(w),\ a_i\in\Z).
\]
For the equivalence relation defining $\LB(w)$, see
\eqref{eq:localized-equiv}--\eqref{eq:localized-equiv-product}.
The upper global basis $\tGup(w)$ is parametrized by $\LB(w)$, and we
denote the basis element corresponding to $x\in\LB(w)$ by $\tGup(x)$.

The localization of $\qA_q(\n(w))$ has a categorical realization
that relates $\eta_w$ to duality.
Kashiwara, Kim, Oh, and Park \cite{KKOP21C} constructed a localization
of a category of finite-dimensional graded modules over quiver Hecke
algebras, in which the simple objects corresponding to the minors
$\qmC_i^w$ become invertible.
In the localization obtained by inverting the simple objects
corresponding to $\qmC_i^w$, every object has a left dual and a
right dual \cite[Theorem~3.9]{KKOP23}, and the simple objects, up to grading
shifts, are parametrized by $\LB(w)$.
Under the bijection $\tGup$ and the parametrization of simple objects
by $\LB(w)$,
the permutation induced by $\eta_w^{-1}$ corresponds to taking
right duals \cite[Theorem~B.22]{IO21_arXiv2}.\footnote{Appendix~B of
\cite{IO21_arXiv2} contains Theorem~B.22 and is not included in the published
version \cite{IO21}.}

When $\g$ is simple and of finite type and $w$ is the longest
Weyl group element, Nakashima \cite{Nakashima22} constructed a
crystal structure on the simple objects, up to grading shifts,
parametrized by $\LB(w)$.
Kashiwara and Nakashima \cite{KN25} constructed a crystal
structure on $\LB(w)$ for arbitrary symmetrizable Kac--Moody algebras
and arbitrary Weyl group elements.

We denote by $\cdual_w$ the permutation of $\LB(w)$ induced by
$\eta_w^{-1}$:
\[
\cdual_w\colon\LB(w)\xrightarrow{\sim}\LB(w),
\qquad
\tGup(\cdual_w(x))=\eta_w^{-1}(\tGup(x)).
\]
The author and Park \cite[Theorem~3.3]{JP25} described $\cdual_w$
using Poincar\'e--Birkhoff--Witt (PBW) data, string data, and left and
right $g$-vectors.
The PBW, string, left $g$-vector, and right $g$-vector parametrizations
assign integer coordinates to elements of
$\LB(w)$ relative to a chosen reduced expression of $w$.
The left and right $g$-vector parametrizations use the quantum
cluster algebra structure of $\widetilde{\qA}_q(\n(w))$
\cite{GLS13,GY21} (see Section~\ref{subsec:g-vector-parametrizations}).
Lee and Park \cite{LP26} used the crystal
description of quantum twists to construct one-parameter families of
imaginary elements in the upper global basis of $\qA_q(\n)$.

Kashiwara's similarity maps \cite{K96} give embeddings
between crystals associated with different Cartan data.
Foldings of Dynkin diagrams give examples of similarity maps:
crystals for the Cartan data obtained by folding can be realized as
subsets of $B(\infty)$ fixed by diagram automorphisms.
Related work includes Naito and Sagaki's study of fixed points of
diagram automorphisms on crystal bases \cite{NS04} and Park's study
of the compatibility of similarity maps with braid group actions
on extended crystals in finite type \cite{Park23}.

To describe Kashiwara's similarity map, let $J$ be a finite set, let
$\vartheta\colon I\twoheadrightarrow J$ be a surjection, and set
$F_j=\vartheta^{-1}(j)$ for $j\in J$.
Choose positive integers $m_i$, for $i\in I$, and impose
the conditions \eqref{eq:similarity-m-condition} on
$\vartheta$ and $(m_i)_{i\in I}$.
In particular, distinct indices in each $F_j$ are orthogonal.
Let $\g_\vartheta$, $\weyl_\vartheta$, and $B_\vartheta(\infty)$
be the Kac--Moody algebra, Weyl group, and crystal determined by
$\vartheta$ and $(m_i)_{i\in I}$ as in
Section~\ref{subsec:similarity-data}.
Write $\Lambda_j^\vartheta$, $j\in J$, for the fundamental weights
of $\g_\vartheta$.
Kashiwara's similarity map \cite{K96} is an injection
\[
\Upsilon_\vartheta\colon
B_\vartheta(\infty)\hookrightarrow B(\infty)
\]
that preserves the highest weight element and satisfies the following
relation between the Kashiwara operators $\tf_j$ on
$B_\vartheta(\infty)$ and $\tf_i$ on $B(\infty)$:
\[
\Upsilon_\vartheta(\tf_j b)
=\Bigl(\prod_{i\in F_j}\tf_i^{\,m_i}\Bigr)\Upsilon_\vartheta(b)
\qquad(j\in J,\ b\in B_\vartheta(\infty))
\]
(see Proposition~\ref{prop:similarity-map}).

In this paper, we extend similarity maps to localized crystals
and establish compatibility with quantum twists.
Let $s_i$ and $s_j^\vartheta$
be the simple reflections in $\weyl$ and $\weyl_\vartheta$,
respectively, and write
$F_j=\{j^{(1)},\ldots,j^{(d_j)}\}$, where $d_j=|F_j|$, in a fixed order.
The assignments
\[
\Theta\colon\weyl_\vartheta\hookrightarrow\weyl,\qquad
\Theta(s_j^\vartheta)=\prod_{i\in F_j}s_i
\qquad(j\in J)
\]
define an injective homomorphism
(Section~\ref{subsec:similarity-data}).
Fix $w\in\weyl_\vartheta$ such that every index in $J$ occurs in
a reduced expression of $w$.
Let $B_\vartheta(w)$ consist of the elements $b\in B_\vartheta(\infty)$
such that the upper global basis element $\Gup(b)$ corresponding to $b$
belongs to the quantum unipotent coordinate ring associated with
$\g_\vartheta$ and $w$.
For $j\in J$, let $\fz_j^w\in B_\vartheta(w)$ satisfy
$\Gup(\fz_j^w)=\qmD(w\Lambda_j^\vartheta,\Lambda_j^\vartheta)$,
and let $\LB_\vartheta(w)$ be the localization of
$B_\vartheta(w)$ at $\fz_j^w$, $j\in J$.
The inverse quantum twists induce permutations $\cdual_w$ on
$\LB_\vartheta(w)$ and $\cdual_{\Theta(w)}$ on
$\LB(\Theta(w))$.

We first prove that $\Upsilon_\vartheta$ restricts to an injection
$B_\vartheta(w)\hookrightarrow B(\Theta(w))$
(Theorem~\ref{thm:similarity-Bw}).
We then extend $\Upsilon_\vartheta|_{B_\vartheta(w)}$ to an injection
\[
\widetilde{\Upsilon}_\vartheta\colon
\LB_\vartheta(w)\hookrightarrow\LB(\Theta(w))
\]
(Proposition~\ref{prop:S-extends-to-localized}).
The sets $\LB_\vartheta(w)$ and $\LB(\Theta(w))$ carry the crystal
structures constructed by Kashiwara and Nakashima \cite{KN25}.
We prove that $\widetilde\Upsilon_\vartheta$ is compatible with the
crystal operators
(Theorem~\ref{thm:localized-crystal-similarity}).
Our main result, Theorem~\ref{thm:main-goal}, states that
\begin{equation}\label{eq:intro-twist-similarity}
\widetilde{\Upsilon}_\vartheta(\cdual_w(x))
=\cdual_{\Theta(w)}(\widetilde{\Upsilon}_\vartheta(x))
\qquad(x\in\LB_\vartheta(w)).
\end{equation}

We also describe $\widetilde{\Upsilon}_\vartheta$ in the PBW,
string, left $g$-vector, and right $g$-vector parametrizations.
Let $\ell$ and $L$ be the lengths of $w$ and $\Theta(w)$,
respectively.
Choose a reduced expression $\bfj=(j_1,\ldots,j_\ell)$ of $w$,
and replace each $j_k$ by
$j_k^{(1)},\ldots,j_k^{(d_{j_k})}$ to obtain
\[
\bfi=(i_1,\ldots,i_L)
=\bigl(j_1^{(1)},\ldots,j_1^{(d_{j_1})},\ldots,
       j_\ell^{(1)},\ldots,j_\ell^{(d_{j_\ell})}\bigr).
\]
The word $\bfi$ is a reduced expression of $\Theta(w)$
(Section~\ref{subsec:similarity-data}).
Define the linear map
\[
\begin{aligned}
\Gamma_\bfj\colon\Z^\ell&\longrightarrow\Z^L,\\
(a_1,\ldots,a_\ell)&\longmapsto
\bigl((m_i a_1)_{i\in F_{j_1}},\ldots,
      (m_i a_\ell)_{i\in F_{j_\ell}}\bigr).
\end{aligned}
\]
The coordinate lists $(m_i a_k)_{i\in F_{j_k}}$ are concatenated for
$k=1,\ldots,\ell$, using the fixed order on each $F_{j_k}$.
We write $\PBW_\bfj,\STR_\bfj,\gL_\bfj,\gR_\bfj$ for the PBW,
string, left $g$-vector, and right $g$-vector parametrizations of
$\LB_\vartheta(w)$, respectively.
For the parametrizations of $\LB(\Theta(w))$, we replace the subscript
$\bfj$ by $\bfi$.
Theorem~\ref{thm:four-param-localized} states that, for every
$x\in\LB_\vartheta(w)$,
\[
\begin{aligned}
\PBW_\bfi(\widetilde{\Upsilon}_\vartheta(x))
&=\Gamma_\bfj(\PBW_\bfj(x)),\\
\STR_\bfi(\widetilde{\Upsilon}_\vartheta(x))
&=\Gamma_\bfj(\STR_\bfj(x)),\\
\gL_\bfi(\widetilde{\Upsilon}_\vartheta(x))
&=\Gamma_\bfj(\gL_\bfj(x)),\\
\gR_\bfi(\widetilde{\Upsilon}_\vartheta(x))
&=\Gamma_\bfj(\gR_\bfj(x)).
\end{aligned}
\]

We prove \eqref{eq:intro-twist-similarity} using the PBW parametrization
$\PBW_\bfi$.
The description of $\PBW_\bfj(\cdual_w(x))$ in
\cite[Theorem~3.3(iv)]{JP25} uses the changes between the PBW,
string, left $g$-vector, and right $g$-vector parametrizations.
Theorem~\ref{thm:four-param-localized} and the linearity of $\Gamma_\bfj$
give equal values of $\PBW_\bfi$ on
$\widetilde{\Upsilon}_\vartheta(\cdual_w(x))$ and
$\cdual_{\Theta(w)}(\widetilde{\Upsilon}_\vartheta(x))$.
The injectivity of $\PBW_\bfi$ then gives \eqref{eq:intro-twist-similarity}.

The identity \eqref{eq:intro-twist-similarity} provides a way to
deduce results about quantum twists in other types from results in
simply laced types through similarity maps.
Let $\PD(w)$ and $\PD(\Theta(w))$ denote the periods of
$\cdual_w$ and $\cdual_{\Theta(w)}$ up to products of integer
powers of $\fz_j^w$, $j\in J$, and
$\fz_i^{\Theta(w)}$, $i\in I$, respectively
(Section~\ref{subsec:equivalence-periods}).
If $\PD(\Theta(w))$ is finite,
Corollary~\ref{prop:periodicity-descends-under-similarity} shows that
$\PD(w)$ is finite and
\[
\PD(w)\mid\PD(\Theta(w)).
\]

We use similarity maps from types $B_n$ and $C_n$ to type
$A_{2n-1}$ to extend the type $A$ period result of
\cite[Theorem~5.18]{JP25} to types $B_n$ and $C_n$.
For $n\ge2$, let $s_1,\ldots,s_n$ be the simple reflections in
types $B_n$ and $C_n$,
numbered along the Dynkin diagram with the multiple edge joining
$n-1$ and $n$.
Let $x_n^{B_n}$ and $x_n^{C_n}$ be the longest elements in the
respective sets of coset representatives of minimal length for
$\weyl_n\backslash\weyl(B_n)$ and
$\weyl_n\backslash\weyl(C_n)$, where
$\weyl_n=\langle s_1,\ldots,s_{n-1}\rangle$ is the maximal standard
parabolic subgroup obtained by deleting vertex $n$.
Theorem~\ref{thm:type-C-xn-period} proves the period formula
conjectured in \cite[Section~5.2.1, Conjecture~1\textup{(i)(b)}]{JP25}:
\[
\PD(x_n^{B_n})=\PD(x_n^{C_n})=
\begin{cases}
2 & \text{if }n=2,\\
4 & \text{if }n\ge3.
\end{cases}
\]

Section~\ref{sec:preliminaries} recalls crystal bases, PBW and string
data, and the Cartan data for similarity.
Section~\ref{sec:similarity-unlocalized} proves
$\Upsilon_\vartheta(B_\vartheta(w))\subset B(\Theta(w))$ and computes
the PBW and string data of $\Upsilon_\vartheta(b)$
for $b\in B_\vartheta(w)$.
Section~\ref{sec:similarity-localized} constructs $\widetilde{\Upsilon}_\vartheta$
and proves Theorems~\ref{thm:four-param-localized} and
\ref{thm:localized-crystal-similarity}.
Section~\ref{sec:similarity-twist} proves Theorem~\ref{thm:main-goal},
Proposition~\ref{prop:similarity-frozen-equivalence}, and
Corollary~\ref{prop:periodicity-descends-under-similarity}.
Section~\ref{sec:application-BC-periods} determines the periods in
types $B_n$ and $C_n$.
\section*{Acknowledgements}
The author thanks Euiyong Park for helpful discussions and for suggesting
the topic of Section~\ref{subsec:localized-crystal-structures}.

The author used SageMath~\cite{sage} to check computations and Codex
(OpenAI) for assistance with SageMath code, literature searches,
and drafting the exposition from the author's notes.

\section{Crystal bases and Cartan data}\label{sec:preliminaries}

We fix the notation for crystal bases and for PBW and string data,
and then specify the Cartan data used for the similarity map.
For integers $a,b$, we write $[a,b]:=\{k\in\Z\mid a\le k\le b\}$.

\subsection{Crystal bases and quantum unipotent coordinate rings}

Let $(A,P,\Pi,P^{\vee},\Pi^{\vee})$ be a symmetrizable Cartan datum with index set $I$.
Thus $A=(a_{i,j})_{i,j\in I}$ is a symmetrizable generalized Cartan matrix, $P$ is the
weight lattice, $P^{\vee}=\Hom_{\Z}(P,\Z)$ is the coweight lattice,
$\Pi=\{\alpha_i\mid i\in I\}\subset P$ is the set of simple roots, and
$\Pi^{\vee}=\{h_i\mid i\in I\}\subset P^{\vee}$ is the set of simple coroots. We assume
that there exists a symmetric bilinear form $(\ ,\ )$ on $P\otimes_{\Z}\Q$ such that
$\langle h_i,\alpha_j\rangle=a_{i,j}$ and
$\langle h_i,\lambda\rangle=2(\alpha_i,\lambda)/(\alpha_i,\alpha_i)$ for
$i,j\in I$ and $\lambda\in P$, and that for each $i\in I$ there exists a fundamental
weight $\Lambda_i\in P$ satisfying $\langle h_j,\Lambda_i\rangle=\delta_{j,i}$ for all
$j\in I$. We normalize the bilinear form $(\ ,\ )$ so that
$(\alpha_i,\alpha_i)/2\in\Z_{>0}$ for every $i\in I$.
We write
\[
P^+:=\{\lambda\in P\mid \langle h_i,\lambda\rangle\ge 0\text{ for all }i\in I\},
\qquad
Q_+:=\bigoplus_{i\in I}\Z_{\ge 0}\alpha_i,
\qquad
Q_-:=-Q_+.
\]
The Weyl group is denoted by
\[
\weyl:=\langle s_i\mid i\in I\rangle,
\]
where $s_i(\lambda)=\lambda-\langle h_i,\lambda\rangle\alpha_i$ for $\lambda\in P$. For
$w\in \weyl$, we denote the length of $w$ by $\ell(w)$ and the set of all reduced
expressions of $w$ by $\rR(w)$. We write $\supp(w)\subset I$ for the set of indices that appear in a
reduced expression of $w$.
We use $\le$ for the Bruhat order on $\weyl$.

Let $\g$ be the Kac--Moody algebra associated with
$(A,P,\Pi,P^{\vee},\Pi^{\vee})$, and let $\n$ denote the negative part of $\g$.
Let $q$ be an indeterminate and set $q_i:=q^{(\alpha_i,\alpha_i)/2}$. The quantized
enveloping algebra $U_q(\g)$ is the $\Q(q)$-algebra generated by $e_i$, $f_i$ $(i\in I)$
and $q^h$ $(h\in P^{\vee})$ subject to the standard Drinfeld--Jimbo relations
(see \cite{HK02,LusztigBook}).
We write $U_q^-(\g)$ for the negative half of $U_q(\g)$.
The quantum unipotent coordinate ring $\qA_q(\n)$ of $U_q(\g)$ is defined by
\begin{equation}\label{eq:qA-n}
\qA_q(\n):=\bigoplus_{\beta\in Q_-}\qA_q(\n)_\beta,
\quad\text{where}\quad
\qA_q(\n)_\beta:=\Hom_{\Q(q)}\bigl(U_q^-(\g)_\beta,\Q(q)\bigr)
\end{equation}
where the multiplication of $\qA_q(\n)$ is defined by using the
comultiplication of $U_q(\g)$ (see \cite[Section~8]{KKKO18} for details).
For a nonzero homogeneous element $x\in\qA_q(\n)_\beta$, write $\wt(x)=\beta$.

A crystal is a set $B$ together with maps
\[
\wt\colon B\to P,\qquad
\ep_i,\varphi_i\colon B\to \Z\sqcup\{-\infty\},\qquad
\te_i,\tf_i\colon B\to B\sqcup\{0\}\qquad (i\in I)
\]
satisfying the crystal axioms \cite[Definition~2.30]{Kimura12}.
For $i\in I$ and $b,b'\in B$, one has
\[
\varphi_i(b)=\ep_i(b)+\langle h_i,\wt(b)\rangle,
\qquad
b'=\te_i b\ \Longleftrightarrow\ b=\tf_i b',
\]
and $\te_i b=\tf_i b=0$ if $\varphi_i(b)=-\infty$.
If $\te_i b\ne0$, then
\[
\wt(\te_i b)=\wt(b)+\alpha_i,\qquad
\ep_i(\te_i b)=\ep_i(b)-1,\qquad
\varphi_i(\te_i b)=\varphi_i(b)+1.
\]
If $\tf_i b\ne0$, then
\[
\wt(\tf_i b)=\wt(b)-\alpha_i,\qquad
\ep_i(\tf_i b)=\ep_i(b)+1,\qquad
\varphi_i(\tf_i b)=\varphi_i(b)-1.
\]

For crystals $B_1$ and $B_2$, the tensor product
\[
B_1\otimes B_2:=\{b_1\otimes b_2\mid b_1\in B_1,\ b_2\in B_2\}
\]
has the crystal structure \cite[Section~1.3]{K93D}.
For $b_1\in B_1$, $b_2\in B_2$, and $i\in I$, 
\begin{equation}\label{eq:crystal-tensor-functions}
\begin{aligned}
\wt(b_1\otimes b_2)&=\wt(b_1)+\wt(b_2),\\
\ep_i(b_1\otimes b_2)&=
\max\{\ep_i(b_1),\ep_i(b_2)-\langle h_i,\wt(b_1)\rangle\},\\
\varphi_i(b_1\otimes b_2)&=
\max\{\varphi_i(b_2),\varphi_i(b_1)+\langle h_i,\wt(b_2)\rangle\}.
\end{aligned}
\end{equation}
We set $0\otimes b_2=b_1\otimes0=0$.
The crystal operators are
\begin{equation}\label{eq:crystal-tensor-operators}
\begin{aligned}
\te_i(b_1\otimes b_2)&=
\begin{cases}
(\te_i b_1)\otimes b_2&\text{if }\varphi_i(b_1)\ge\ep_i(b_2),\\
b_1\otimes(\te_i b_2)&\text{if }\varphi_i(b_1)<\ep_i(b_2),
\end{cases}\\
\tf_i(b_1\otimes b_2)&=
\begin{cases}
(\tf_i b_1)\otimes b_2&\text{if }\varphi_i(b_1)>\ep_i(b_2),\\
b_1\otimes(\tf_i b_2)&\text{if }\varphi_i(b_1)\le\ep_i(b_2).
\end{cases}
\end{aligned}
\end{equation}

Kashiwara constructed a crystal $B(\infty)$ from $U_q^-(\g)$
\cite[Theorem~4]{K91}.
The upper global basis (or dual canonical basis) of $\qA_q(\n)$
is indexed by $B(\infty)$ \cite{K91,K93}.
We use the notation of \cite[Section~3.1]{Kimura12}.
We write
\[
\Gup(\infty):=\{\Gup(b)\mid b\in B(\infty)\},
\qquad \Gup(b)\in\qA_q(\n)_{\wt(b)}.
\]
We denote by $\one$ the highest element of $B(\infty)$.
For $b\in B(\infty)$, one has
\[
\ep_i(b)=\max\{n\ge0\mid\te_i^n b\ne0\}.
\]

Kashiwara's involution $*$ on $B(\infty)$ is induced by the
$\Q(q)$-linear algebra anti-involution of $U_q^-(\g)$ fixing every $f_i$; see \cite[Theorem~2.1.1]{K93D} and the formulation in
\cite[(2.1) and Theorem~2.15]{Kimura12}.

The crystal operators and Kashiwara's involution are extended by zero.
We write
\[
\begin{aligned}
\te_i^*&:=*\,\te_i\,*, & \tf_i^*&:=*\,\tf_i\,*,\\
\ep_i^*(b)&:=\ep_i(b^*), & \varphi_i^*(b)&:=\varphi_i(b^*),\\
\te_i^{\max}(b)&:=\te_i^{\ep_i(b)}b, &
(\te_i^*)^{\max}(b)&:=(\te_i^*)^{\ep_i^*(b)}b.
\end{aligned}
\]

For $\Lambda\in P^+$, let $V(\Lambda)$ be the integrable highest weight
$U_q(\g)$-module with highest weight vector $u_\Lambda$.
We use the symmetric bilinear form $(\ ,\ )_\Lambda$ on $V(\Lambda)$, normalized by
$(u_\Lambda,u_\Lambda)_\Lambda=1$, for which $e_i$ and $f_i$ are adjoint
and $q^h$ is self-adjoint.
For $w\in\weyl$, let $u_{w\Lambda}$ be the vector of weight $w\Lambda$
in the lower global basis of $V(\Lambda)$; see
\cite[Definition~3.13 and Proposition~3.33]{KiOy21}.
Following \cite[Definition~3.40]{KiOy21}, for $v,w\in\weyl$ define the
quantum unipotent minor $\qmD(w\Lambda,v\Lambda)$ by
\[
\qmD(w\Lambda,v\Lambda)(x):=(u_{w\Lambda},xu_{v\Lambda})_\Lambda
\qquad(x\in U_q^-(\g)).
\]
If $v\le w$, then
\[
\qmD(w\Lambda,v\Lambda)\in\Gup(\infty),\qquad
\wt(\qmD(w\Lambda,v\Lambda))=w\Lambda-v\Lambda;
\]
see \cite[Proposition~3.46]{KiOy21}.

\subsection{PBW and string data}

Fix $w\in \weyl$ and a reduced expression
$\bfi=(i_1,\dots,i_\ell)\in \rR(w)$. We set
$w_{\le 0}:=1$, $w_{\le k}:=s_{i_1}\cdots s_{i_k}$, and $w_{<k}:=w_{\le k-1}$ for
$1\le k\le \ell$, and define the associated positive roots by
\begin{equation}\label{eq:beta-i-k}
\beta_{\bfi,k}:=w_{<k}(\alpha_{i_k})
\qquad (1\le k\le \ell).
\end{equation}
For $1\le k\le\ell$, the dual PBW root vector
$\mathsf{F}^{\mathrm{up}}_{\bfi}(\beta_{\bfi,k})$ is the quantum minor
\[
\mathsf{F}^{\mathrm{up}}_{\bfi}(\beta_{\bfi,k})
=\qmD(w_{\le k}\Lambda_{i_k},w_{<k}\Lambda_{i_k});
\]
see \cite[Sections~1.1 and~2.2]{JP25}.
The quantum unipotent coordinate ring attached to $w$ is the subalgebra
\begin{equation}\label{eq:qA-nw}
\qA_q(\n(w))\subset \qA_q(\n)
\end{equation}
generated by the dual PBW root vectors
$\mathsf{F}^{\mathrm{up}}_{\bfi}(\beta_{\bfi,k})$ $(1\le k\le \ell)$. The algebra $\qA_q(\n(w))$ is independent
of the choice of $\bfi$ \cite[Theorem~4.9 and Proposition~4.11]{Kimura12}.
We set
\begin{equation}\label{eq:Gup-w}
\Gup(w):=\qA_q(\n(w))\cap \Gup(\infty).
\end{equation}
By \cite[Theorem~4.25]{Kimura12},
$\Gup(w)$ is a basis of $\qA_q(\n(w))$.

For $1\le k\le \ell$ and $a\in\Z_{\ge0}$, we use the convention
\[
\mathsf{F}^{\mathrm{up}}_{\bfi}(\beta_{\bfi,k})^{\{a\}}
:=
q_{i_k}^{a(a-1)/2}\mathsf{F}^{\mathrm{up}}_{\bfi}(\beta_{\bfi,k})^a,
\]
as in \cite[Section~1.1]{JP25}.
For $\bfa=(a_1,\dots,a_\ell)\in \Z_{\ge 0}^{[1,\ell]}$, define the dual PBW monomial
\begin{equation}\label{eq:dual-pbw-monomial}
\mathsf{F}^{\mathrm{up}}_{\bfi}(\bfa):=
\mathsf{F}^{\mathrm{up}}_{\bfi}(\beta_{\bfi,\ell})^{\{a_\ell\}}
\cdots
\mathsf{F}^{\mathrm{up}}_{\bfi}(\beta_{\bfi,1})^{\{a_1\}}.
\end{equation}
The monomials $\mathsf{F}^{\mathrm{up}}_{\bfi}(\bfa)$, with
$\bfa\in\Z_{\ge0}^{[1,\ell]}$, form the dual PBW basis of
$\qA_q(\n(w))$ (see \cite[Proposition~40.2.1]{LusztigBook} and \cite[Section~4.7.2]{Kimura12}).
By \cite[Theorem~4.29]{Kimura12}, for each $x\in\Gup(w)$ there is a unique vector
$\bfa\in\Z_{\ge0}^{[1,\ell]}$ such that
\begin{equation}\label{eq:pbw-datum-leading-term}
\mathsf{F}^{\mathrm{up}}_{\bfi}(\bfa)\equiv x
\pmod{q\,\mathrm{span}_{\Z[q]}\Gup(\infty)}.
\end{equation}
The vector $\bfa$ is called the PBW datum of $x$ with respect to $\bfi$ and is denoted by
$\PBW_{\bfi}(x)$.

The crystal basis of $\qA_q(\n(w))$ is
\begin{equation}\label{eq:Bw-definition}
B(w):=\{b\in B(\infty)\mid \Gup(b)\in \qA_q(\n(w))\}.
\end{equation}
For $b\in B(w)$, the PBW datum of $b$ with respect to $\bfi$ is
\begin{equation}\label{eq:PBW-on-Bw}
\PBW_{\bfi}(b):=\PBW_{\bfi}(\Gup(b)).
\end{equation}

The string datum of $b\in B(w)$ with respect to
$\bfi$ is
\begin{equation}\label{eq:STR-on-Bw}
\STR_{\bfi}(b)=(t_1,\dots,t_\ell),
\qquad
t_k:=\ep_{i_k}\bigl(\te_{i_{k-1}}^{t_{k-1}}\cdots \te_{i_1}^{t_1}(b)\bigr)
\qquad (1\le k\le \ell).
\end{equation}

\subsection{The Cartan datum associated with \texorpdfstring{$\vartheta$}{vartheta}}
\label{subsec:similarity-data}

\begin{figure}[!htbp]
\centering
\begin{minipage}{0.47\textwidth}
\centering
\begin{tikzpicture}[x=1.25cm,y=0.95cm,baseline=(current bounding box.center),
 every node/.style={font=\small},
 vertex/.style={circle,draw,fill=white,inner sep=0pt,minimum size=4.5pt}]
\node[vertex,label={above:$1$}] (a) at (0,0.55) {};
\node[vertex,label={above:$2$}] (b) at (1,0.55) {};
\node[vertex,label={above:$3$}] (c) at (2,0) {};
\node[vertex,label={below:$4$}] (d) at (1,-0.55) {};
\node[vertex,label={below:$5$}] (e) at (0,-0.55) {};
\draw (a) -- (b);
\draw (b) -- (c);
\draw (c) -- (d);
\draw (d) -- (e);
\node[anchor=east] at (-0.85,0) {$A_5$};
\draw[->,densely dashed,gray] (0,-1.05) -- (0,-1.65);
\node[vertex,label={below:$1$}] (f1) at (0,-1.95) {};
\draw[->,densely dashed,gray] (1,-1.05) -- (1,-1.65);
\node[vertex,label={below:$2$}] (f2) at (1,-1.95) {};
\draw[->,densely dashed,gray] (2,-1.05) -- (2,-1.65);
\node[vertex,label={below:$3$}] (f3) at (2,-1.95) {};
\node[anchor=east] at (-0.45,-1.95) {$B_3$};
\draw (f1) -- (f2);
\draw (1.065,-1.9949999999999999) -- (1.935,-1.9949999999999999);
\draw (1.065,-1.905) -- (1.935,-1.905);
\draw (1.43,-1.83) -- (1.59,-1.95) -- (1.43,-2.07);
\end{tikzpicture}
\par\smallskip
(a)
\end{minipage}
\hfill
\begin{minipage}{0.47\textwidth}
\centering
\begin{tikzpicture}[x=1.25cm,y=0.95cm,baseline=(current bounding box.center),
 every node/.style={font=\small},
 vertex/.style={circle,draw,fill=white,inner sep=0pt,minimum size=4.5pt}]
\node[vertex,label={above:$1$}] (a) at (0,0.55) {};
\node[vertex,label={above:$2$}] (b) at (1,0.55) {};
\node[vertex,label={above:$3\ (2)$}] (c) at (2,0) {};
\node[vertex,label={below:$4$}] (d) at (1,-0.55) {};
\node[vertex,label={below:$5$}] (e) at (0,-0.55) {};
\draw (a) -- (b);
\draw (b) -- (c);
\draw (c) -- (d);
\draw (d) -- (e);
\node[anchor=east] at (-0.85,0) {$A_5$};
\draw[->,densely dashed,gray] (0,-1.05) -- (0,-1.65);
\node[vertex,label={below:$1$}] (f1) at (0,-1.95) {};
\draw[->,densely dashed,gray] (1,-1.05) -- (1,-1.65);
\node[vertex,label={below:$2$}] (f2) at (1,-1.95) {};
\draw[->,densely dashed,gray] (2,-1.05) -- (2,-1.65);
\node[vertex,label={below:$3$}] (f3) at (2,-1.95) {};
\node[anchor=east] at (-0.45,-1.95) {$C_3$};
\draw (f1) -- (f2);
\draw (1.065,-1.9949999999999999) -- (1.935,-1.9949999999999999);
\draw (1.065,-1.905) -- (1.935,-1.905);
\draw (1.57,-1.83) -- (1.41,-1.95) -- (1.57,-2.07);
\end{tikzpicture}
\par\smallskip
(b)
\end{minipage}
\par\smallskip
\begin{minipage}{0.47\textwidth}
\centering
\begin{tikzpicture}[x=1.25cm,y=0.95cm,baseline=(current bounding box.center),
 every node/.style={font=\small},
 vertex/.style={circle,draw,fill=white,inner sep=0pt,minimum size=4.5pt}]
\node[vertex,label={left:$1$}] (a) at (0,0.8) {};
\node[vertex,label={right:$2$}] (b) at (1,0) {};
\node[vertex,label={left:$3$}] (c) at (0,0) {};
\node[vertex,label={left:$4$}] (d) at (0,-0.8) {};
\draw (a) -- (b);
\draw (c) -- (b);
\draw (d) -- (b);
\node[anchor=east] at (-0.45,0) {$D_4$};
\draw[->,densely dashed,gray] (0,-1.05) -- (0,-1.65);
\node[vertex,label={below:$1$}] (f1) at (0,-1.95) {};
\draw[->,densely dashed,gray] (1,-1.05) -- (1,-1.65);
\node[vertex,label={below:$2$}] (f2) at (1,-1.95) {};
\node[anchor=east] at (-0.45,-1.95) {$G_2$};
\draw (0.065,-2.02) -- (0.935,-2.02);
\draw (0.065,-1.95) -- (0.935,-1.95);
\draw (0.065,-1.88) -- (0.935,-1.88);
\draw (0.43,-1.83) -- (0.59,-1.95) -- (0.43,-2.07);
\end{tikzpicture}
\par\smallskip
(c)
\end{minipage}
\hfill
\begin{minipage}{0.47\textwidth}
\centering
\begin{tikzpicture}[x=1.25cm,y=0.95cm,baseline=(current bounding box.center),
 every node/.style={font=\small},
 vertex/.style={circle,draw,fill=white,inner sep=0pt,minimum size=4.5pt}]
\node[vertex,label={left:$1$}] (a) at (0,0.8) {};
\node[vertex,label={right:$2$}] (b) at (1,0.45) {};
\node[vertex,label={left:$3\ (2)$}] (c) at (0,0) {};
\node[vertex,label={right:$4$}] (d) at (1,-0.45) {};
\node[vertex,label={left:$5$}] (e) at (0,-0.8) {};
\draw (a) -- (b);
\draw (b) -- (c);
\draw (c) -- (d);
\draw (d) -- (e);
\node[anchor=east] at (-0.85,0) {$A_5$};
\draw[->,densely dashed,gray] (0,-1.05) -- (0,-1.65);
\node[vertex,label={below:$1$}] (f1) at (0,-1.95) {};
\draw[->,densely dashed,gray] (1,-1.05) -- (1,-1.65);
\node[vertex,label={below:$2$}] (f2) at (1,-1.95) {};
\node[anchor=east] at (-0.45,-1.95) {$G_2$};
\draw (0.065,-2.02) -- (0.935,-2.02);
\draw (0.065,-1.95) -- (0.935,-1.95);
\draw (0.065,-1.88) -- (0.935,-1.88);
\draw (0.43,-1.83) -- (0.59,-1.95) -- (0.43,-2.07);
\end{tikzpicture}
\par\smallskip
(d)
\end{minipage}
\par\smallskip
\begin{minipage}{0.47\textwidth}
\centering
\begin{tikzpicture}[x=1.25cm,y=0.95cm,baseline=(current bounding box.center),
 every node/.style={font=\small},
 vertex/.style={circle,draw,fill=white,inner sep=0pt,minimum size=4.5pt}]
\node[vertex,label={above:$1$}] (a) at (0,0.55) {};
\node[vertex,label={above:$3$}] (b) at (1,0.55) {};
\node[vertex,label={above:$4$}] (c) at (2,0) {};
\node[vertex,label={below:$5$}] (d) at (1,-0.55) {};
\node[vertex,label={below:$6$}] (e) at (0,-0.55) {};
\node[vertex,label={above:$2$}] (f) at (3,0) {};
\draw (a) -- (b);
\draw (b) -- (c);
\draw (c) -- (d);
\draw (d) -- (e);
\draw (c) -- (f);
\node[anchor=east] at (-0.45,0) {$E_6$};
\draw[->,densely dashed,gray] (0,-1.05) -- (0,-1.65);
\node[vertex,label={below:$1$}] (f1) at (0,-1.95) {};
\draw[->,densely dashed,gray] (1,-1.05) -- (1,-1.65);
\node[vertex,label={below:$2$}] (f2) at (1,-1.95) {};
\draw[->,densely dashed,gray] (2,-1.05) -- (2,-1.65);
\node[vertex,label={below:$3$}] (f3) at (2,-1.95) {};
\draw[->,densely dashed,gray] (3,-1.05) -- (3,-1.65);
\node[vertex,label={below:$4$}] (f4) at (3,-1.95) {};
\node[anchor=east] at (-0.45,-1.95) {$F_4$};
\draw (f1) -- (f2);
\draw (1.065,-1.9949999999999999) -- (1.935,-1.9949999999999999);
\draw (1.065,-1.905) -- (1.935,-1.905);
\draw (1.43,-1.83) -- (1.59,-1.95) -- (1.43,-2.07);
\draw (f3) -- (f4);
\end{tikzpicture}
\par\smallskip
(e)
\end{minipage}
\par\smallskip
\begin{minipage}{0.47\textwidth}
\centering
\begin{tikzpicture}[x=1.25cm,y=0.95cm,baseline=(current bounding box.center),
 every node/.style={font=\small},
 vertex/.style={circle,draw,fill=white,inner sep=0pt,minimum size=4.5pt}]
\node[vertex,label={above:$0$}] (a) at (0,0) {};
\node[vertex,label={above:$1$}] (b) at (1,0.55) {};
\node[vertex,label={above:$2$}] (c) at (2,0.55) {};
\node[vertex,label={above:$3$}] (d) at (3,0) {};
\node[vertex,label={below:$4$}] (e) at (2,-0.55) {};
\node[vertex,label={below:$5$}] (f) at (1,-0.55) {};
\draw (a) -- (b);
\draw (b) -- (c);
\draw (c) -- (d);
\draw (d) -- (e);
\draw (e) -- (f);
\draw (f) -- (a);
\node[anchor=east] at (-0.45,0) {$A_5^{(1)}$};
\draw[->,densely dashed,gray] (0,-1.05) -- (0,-1.65);
\node[vertex,label={below:$0$}] (f0) at (0,-1.95) {};
\draw[->,densely dashed,gray] (1,-1.05) -- (1,-1.65);
\node[vertex,label={below:$1$}] (f1) at (1,-1.95) {};
\draw[->,densely dashed,gray] (2,-1.05) -- (2,-1.65);
\node[vertex,label={below:$2$}] (f2) at (2,-1.95) {};
\draw[->,densely dashed,gray] (3,-1.05) -- (3,-1.65);
\node[vertex,label={below:$3$}] (f3) at (3,-1.95) {};
\node[anchor=east] at (-0.45,-1.95) {$D_4^{(2)}$};
\draw (0.065,-1.9949999999999999) -- (0.935,-1.9949999999999999);
\draw (0.065,-1.905) -- (0.935,-1.905);
\draw (0.5700000000000001,-1.83) -- (0.41000000000000003,-1.95) -- (0.5700000000000001,-2.07);
\draw (f1) -- (f2);
\draw (2.065,-1.9949999999999999) -- (2.935,-1.9949999999999999);
\draw (2.065,-1.905) -- (2.935,-1.905);
\draw (2.43,-1.83) -- (2.59,-1.95) -- (2.43,-2.07);
\end{tikzpicture}
\par\smallskip
(f)
\end{minipage}
\hfill
\begin{minipage}{0.47\textwidth}
\centering
\begin{tikzpicture}[x=1.25cm,y=0.95cm,baseline=(current bounding box.center),
 every node/.style={font=\small},
 vertex/.style={circle,draw,fill=white,inner sep=0pt,minimum size=4.5pt}]
\node[vertex,label={above:$0\ (2)$}] (a) at (0,0) {};
\node[vertex,label={above:$1$}] (b) at (1,0.55) {};
\node[vertex,label={above:$2$}] (c) at (2,0.55) {};
\node[vertex,label={above:$3\ (2)$}] (d) at (3,0) {};
\node[vertex,label={below:$4$}] (e) at (2,-0.55) {};
\node[vertex,label={below:$5$}] (f) at (1,-0.55) {};
\draw (a) -- (b);
\draw (b) -- (c);
\draw (c) -- (d);
\draw (d) -- (e);
\draw (e) -- (f);
\draw (f) -- (a);
\node[anchor=east] at (-0.45,0) {$A_5^{(1)}$};
\draw[->,densely dashed,gray] (0,-1.05) -- (0,-1.65);
\node[vertex,label={below:$0$}] (f0) at (0,-1.95) {};
\draw[->,densely dashed,gray] (1,-1.05) -- (1,-1.65);
\node[vertex,label={below:$1$}] (f1) at (1,-1.95) {};
\draw[->,densely dashed,gray] (2,-1.05) -- (2,-1.65);
\node[vertex,label={below:$2$}] (f2) at (2,-1.95) {};
\draw[->,densely dashed,gray] (3,-1.05) -- (3,-1.65);
\node[vertex,label={below:$3$}] (f3) at (3,-1.95) {};
\node[anchor=east] at (-0.45,-1.95) {$C_3^{(1)}$};
\draw (0.065,-1.9949999999999999) -- (0.935,-1.9949999999999999);
\draw (0.065,-1.905) -- (0.935,-1.905);
\draw (0.43,-1.83) -- (0.59,-1.95) -- (0.43,-2.07);
\draw (f1) -- (f2);
\draw (2.065,-1.9949999999999999) -- (2.935,-1.9949999999999999);
\draw (2.065,-1.905) -- (2.935,-1.905);
\draw (2.57,-1.83) -- (2.41,-1.95) -- (2.57,-2.07);
\end{tikzpicture}
\par\smallskip
(g)
\end{minipage}
\caption{Cartan data obtained from $\vartheta$ and $(m_i)_{i\in I}$.}
\label{fig:similarity-data}
\end{figure}
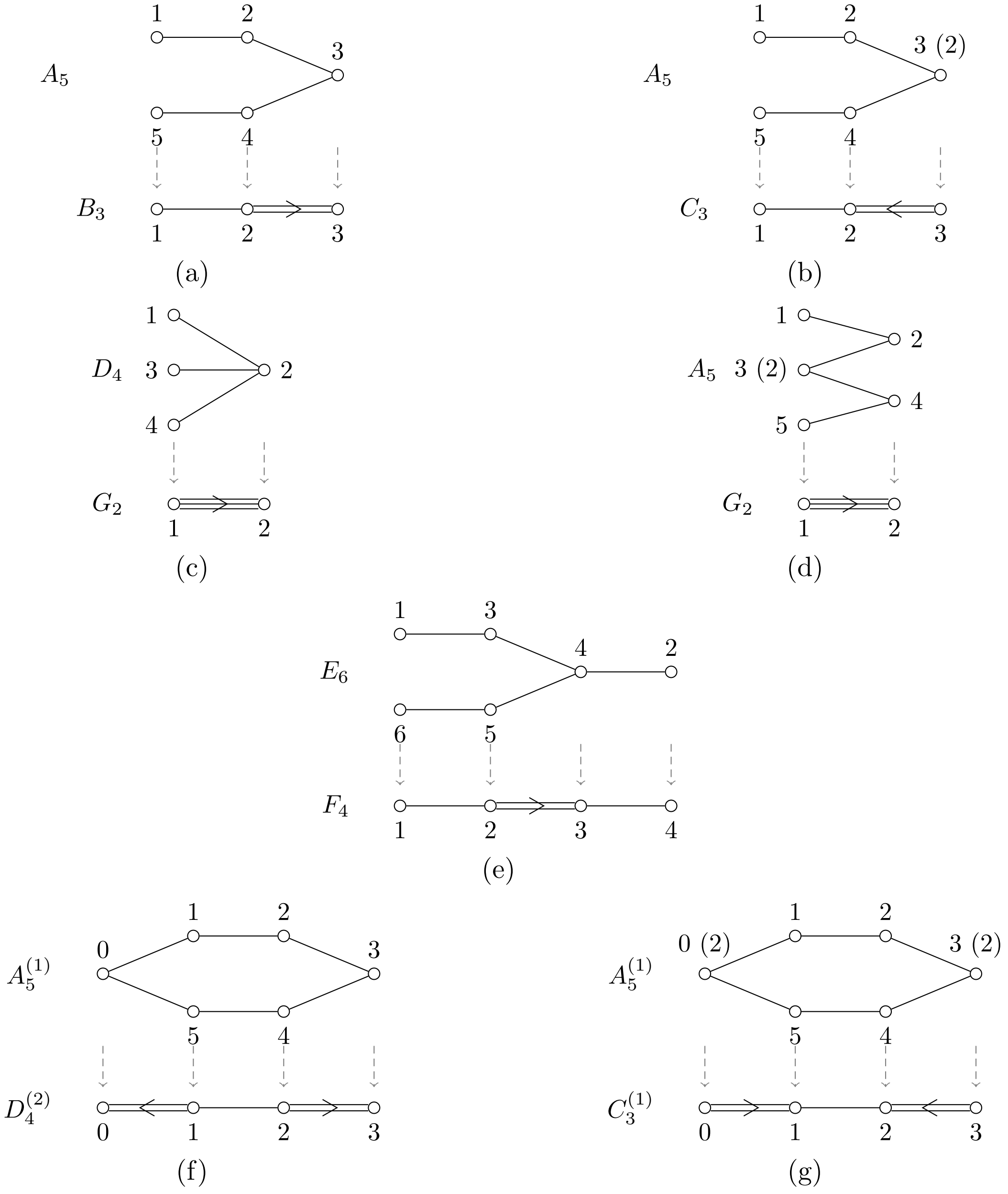

Following \cite[Section~5]{K96}, let
\[
\vartheta\colon I\twoheadrightarrow J
\]
be a surjective map onto a finite set $J$, and choose a positive integer
$m_i$ for each $i\in I$.
For $j\in J$, write
\[
F_j:=\vartheta^{-1}(j)=\{j^{(1)},\dots,j^{(d_j)}\},\qquad d_j:=|F_j|,
\]
with a fixed order on each fiber, and set
\begin{equation}\label{eq:similarity-alpha-hat}
\alpha_j^\vartheta:=\sum_{i\in F_j}m_i\alpha_i.
\end{equation}
Let $\widehat P$ be the subset of $P$ consisting of the weights $\lambda$
such that, for every $j\in J$, the number
$\langle h_i,\lambda\rangle/m_i$ is an integer and does not depend on
the choice of $i\in F_j$.
For $j\in J$, define $h_j^\vartheta\in\Hom_\Z(\widehat P,\Z)$ by
\begin{equation}\label{eq:similarity-h-hat}
\langle h_j^\vartheta,\lambda\rangle
=\frac1{m_i}\langle h_i,\lambda\rangle
\qquad(i\in F_j,\ \lambda\in\widehat P).
\end{equation}
We assume the two conditions in \cite[(5.1) and~(5.2)]{K96}:
\begin{equation}\label{eq:similarity-m-condition}
a_{i,i'}=0\quad(i\ne i',\ i,i'\in F_j),\qquad
\alpha_j^\vartheta\in\widehat P
\quad(j\in J).
\end{equation}
By \cite[Section~5]{K96}, under \eqref{eq:similarity-m-condition},
the lattice $\widehat P$ with simple roots $(\alpha_j^\vartheta)_{j\in J}$
and simple coroots $(h_j^\vartheta)_{j\in J}$ defines a symmetrizable
Cartan datum on $J$, with Cartan matrix
\begin{equation}\label{eq:A-vartheta}
A_\vartheta=(c_{j,k})_{j,k\in J},\qquad
c_{j,k}:=\langle h_j^\vartheta,\alpha_k^\vartheta\rangle.
\end{equation}
Let $\g_\vartheta$ and $\weyl_\vartheta$ be the Kac--Moody algebra and
Weyl group associated with $A_\vartheta$, respectively. We use the same notation $\ell$, $\rR$, and
$\supp$ for length, reduced expressions, and support in $\weyl_\vartheta$.
We identify $P(\g_\vartheta)$ with $\widehat P$ and write $P(\g)=P$.
The fundamental weights of $\g_\vartheta$ are
\[
\Lambda_j^\vartheta:=\sum_{i\in F_j}m_i\Lambda_i
\qquad(j\in J).
\]
 Figure~\ref{fig:similarity-data} gives examples of the Cartan data defined by
$\vartheta$ and $(m_i)_{i\in I}$ \cite[Section~5]{K96}.
The upper and lower Dynkin diagrams correspond to $\g$ and
$\g_\vartheta$, respectively.
The map $\vartheta$ sends each vertex of the upper diagram to the
vertically aligned vertex of the lower diagram.
We have $m_i=1$ unless indicated in parentheses.

Let $\iota\colon P(\g_\vartheta)=\widehat P\lhook\joinrel\longrightarrow P(\g)$
be the inclusion map. In particular,
\[
\iota(\Lambda_j^\vartheta)=\sum_{i\in F_j}m_i\Lambda_i
\qquad (j\in J).
\]
By \eqref{eq:similarity-h-hat},
\begin{equation}\label{eq:similarity-coroot-pairing}
\langle h_i,\iota(\lambda)\rangle
=m_i\langle h_j^\vartheta,\lambda\rangle
\qquad (j\in J,\ i\in F_j,\ \lambda\in P(\g_\vartheta)).
\end{equation}
For $j\in J$, let $s_j^\vartheta$ be the simple reflection of
$\weyl_\vartheta$.
For $j\in J$ and $\lambda\in P(\g_\vartheta)$, orthogonality
within $F_j$ and \eqref{eq:similarity-coroot-pairing} give
\[
\begin{aligned}
\Bigl(\prod_{i\in F_j}s_i\Bigr)\iota(\lambda)
&=\iota(\lambda)-\sum_{i\in F_j}
  \langle h_i,\iota(\lambda)\rangle\alpha_i\\
&=\iota(\lambda)-\langle h_j^\vartheta,\lambda\rangle
  \sum_{i\in F_j}m_i\alpha_i
=\iota(s_j^\vartheta\lambda).
\end{aligned}
\]
The products $\prod_{i\in F_j}s_i$, $j\in J$, therefore preserve
$\widehat P$ and act on $\widehat P$ as the simple reflections
$s_j^\vartheta$ of $\weyl_\vartheta$.
Set $\rho=\sum_{j\in J}\Lambda_j^\vartheta$.
Since $\langle h_i,\iota(\rho)\rangle=m_i>0$ for every $i\in I$,
the weight $\iota(\rho)$ is strictly dominant for $\g$.
A Weyl group element acting trivially on $\widehat P$ fixes
$\iota(\rho)$ and hence is the identity.
Consequently, restriction to $\widehat P$ identifies
$\langle\prod_{i\in F_j}s_i\mid j\in J\rangle$ with
$\weyl_\vartheta$.
The inverse of the restriction to $\widehat P$ gives an injective
homomorphism
\[
\Theta\colon\weyl_\vartheta\hookrightarrow\weyl,
\qquad
\Theta(s_j^\vartheta)=\prod_{i\in F_j}s_i
\quad(j\in J).
\]
For $u\in\weyl_\vartheta$ and $\lambda\in P(\g_\vartheta)$, we have
\begin{equation}\label{eq:similarity-weyl-equivariance}
\Theta(u)\iota(\lambda)=\iota(u\lambda)
\qquad(u\in\weyl_\vartheta,\ \lambda\in P(\g_\vartheta)).
\end{equation}

Let $w\in\weyl_\vartheta$, choose
$\bfj=(j_1,\ldots,j_\ell)\in\rR(w)$, and define
\begin{equation}\label{eq:expanded-reduced-expression}
\bfi:=\bigl(j_1^{(1)},\ldots,j_1^{(d_{j_1})},\ldots,
            j_\ell^{(1)},\ldots,j_\ell^{(d_{j_\ell})}\bigr).
\end{equation}
\begin{lemma}\label{lem:expanded-reduced-expression}
The word $\bfi$ defined in \eqref{eq:expanded-reduced-expression}
belongs to $\rR(\Theta(w))$. If $\supp(w)=J$, then $\supp(\Theta(w))=I$.
\end{lemma}

\begin{proof}
Let $u\in\weyl_\vartheta$ and $j\in J$ satisfy
$\ell(us_j^\vartheta)>\ell(u)$.
Then $u h_j^\vartheta$ is a positive coroot.
By \eqref{eq:similarity-coroot-pairing} and
\eqref{eq:similarity-weyl-equivariance},
\[
\langle\Theta(u)h_i,\iota(\rho)\rangle
=m_i\langle u h_j^\vartheta,\rho\rangle>0
\qquad(i\in F_j).
\]
Thus $\Theta(u)h_i$ is positive for every $i\in F_j$.
Since $s_i h_{i'}=h_{i'}$ for distinct $i,i'\in F_j$,
the length increases at each step when $\Theta(u)$ is multiplied
successively on the right by $s_{j^{(1)}},\ldots,s_{j^{(d_j)}}$.
Hence $\bfi\in\rR(\Theta(w))$.
If $\supp(w)=J$, then every index in $J$ occurs in $\bfj$.
Hence every index in $I$ occurs in $\bfi$, and
$\supp(\Theta(w))=I$.
\end{proof}

For instance, let $\g$ have type $A_3$ with $I=\{1,2,3\}$
numbered along the Dynkin diagram. Take $J=\{1,2\}$,
$\vartheta(1)=\vartheta(3)=1$, and $\vartheta(2)=2$. Then
\[
\Theta(s_1^\vartheta)=s_1s_3,\qquad
\Theta(s_2^\vartheta)=s_2.
\]

Write $B_\vartheta(\infty)$ for the crystal of
$U_q^-(\g_\vartheta)$, with highest element also denoted by $\one$.
On $B_\vartheta(\infty)$, we use the same symbols $*$, $\te_j^*$,
$\tf_j^*$, $\ep_j^*$, and $\varphi_j^*$ for Kashiwara's
involution and the corresponding starred operators and functions.

For $w\in\weyl_\vartheta$, $B_\vartheta(w)$ denotes the subset of
$B_\vartheta(\infty)$ obtained by applying the definition
\eqref{eq:Bw-definition} to $\g_\vartheta$ and the Weyl group element $w$.

\FloatBarrier

\subsection{Kashiwara's similarity map}

We recall the action of Kashiwara's similarity map $\Upsilon_\vartheta$ on weights
and crystal functions.

\begin{proposition}[{\cite[Section~5]{K96}; cf.\ \cite[Theorem~4.1]{Park23}}]
\label{prop:similarity-map}
There is a unique injective map
\begin{equation}\label{eq:crystal-similarity}
\Upsilon_\vartheta\colon B_\vartheta(\infty)\longrightarrow B(\infty),
\end{equation}
extended by $\Upsilon_\vartheta(0)=0$, satisfying
$\Upsilon_\vartheta(\one)=\one$ and
\begin{equation}\label{eq:similarity-e-f}
\begin{aligned}
\Upsilon_\vartheta(\te_j b)
&=\Bigl(\prod_{i\in F_j}\te_i^{m_i}\Bigr)\Upsilon_\vartheta(b),\\
\Upsilon_\vartheta(\tf_j b)
&=\Bigl(\prod_{i\in F_j}\tf_i^{m_i}\Bigr)\Upsilon_\vartheta(b)
\end{aligned}
\qquad (j\in J,\ b\in B_\vartheta(\infty)).
\end{equation}
The map $\Upsilon_\vartheta$ commutes with Kashiwara's involution:
\[
\Upsilon_\vartheta(b^*)=\Upsilon_\vartheta(b)^*
\qquad (b\in B_\vartheta(\infty)).
\]
The map $\Upsilon_\vartheta$ also satisfies
\begin{equation}\label{eq:similarity-e-f-star}
\begin{aligned}
\Upsilon_\vartheta(\te_j^* b)
&=\Bigl(\prod_{i\in F_j}(\te_i^*)^{m_i}\Bigr)\Upsilon_\vartheta(b),\\
\Upsilon_\vartheta(\tf_j^* b)
&=\Bigl(\prod_{i\in F_j}(\tf_i^*)^{m_i}\Bigr)\Upsilon_\vartheta(b)
\end{aligned}
\qquad (j\in J,\ b\in B_\vartheta(\infty)).
\end{equation}
For every $b\in B_\vartheta(\infty)$,
\begin{equation}\label{eq:similarity-weight}
\wt(\Upsilon_\vartheta(b))=\iota(\wt(b)).
\end{equation}
For $j\in J$, $i\in F_j$, and $b\in B_\vartheta(\infty)$, one has
\begin{align}
\ep_i(\Upsilon_\vartheta(b))&=m_i\ep_j(b),
\qquad
\ep_i^*(\Upsilon_\vartheta(b))=m_i\ep_j^*(b),
\label{eq:similarity-epsilon-fiber}\\
\varphi_i(\Upsilon_\vartheta(b))&=m_i\varphi_j(b),
\qquad
\varphi_i^*(\Upsilon_\vartheta(b))=m_i\varphi_j^*(b).
\label{eq:similarity-phi-fiber}
\end{align}
\end{proposition}

The products in \eqref{eq:similarity-e-f} and \eqref{eq:similarity-e-f-star}
are independent of the order of multiplication.
Indeed, for distinct $i,i'\in F_j$, the operators $\te_i,\te_{i'}$ commute,
as do $\tf_i,\tf_{i'}$, by \cite[Section~5]{K96}.
Conjugating $\te_i\te_{i'}=\te_{i'}\te_i$ and $\tf_i\tf_{i'}=\tf_{i'}\tf_i$
by Kashiwara's involution gives $\te_i^*\te_{i'}^*=\te_{i'}^*\te_i^*$
and $\tf_i^*\tf_{i'}^*=\tf_{i'}^*\tf_i^*$.

\section{The similarity map on \texorpdfstring{$B(w)$}{B(w)}}
\label{sec:similarity-unlocalized}
For $w\in\weyl_\vartheta$ with $\supp(w)=J$, we prove that
$\Upsilon_\vartheta(B_\vartheta(w))\subset B(\Theta(w))$ and compute
the PBW and string data of $\Upsilon_\vartheta(b)$ for $b\in B_\vartheta(w)$
(Theorem~\ref{thm:similarity-Bw} and
Proposition~\ref{prop:string-Bw}).
We use the Cartan data and the map $\Theta$ fixed in
Section~\ref{subsec:similarity-data}.

\subsection{Saito reflections}

Following the notation of \cite[Section~2]{Park23}, for $i\in I$ set
\begin{equation}\label{eq:saito-domains}
{}_{i}B(\infty):=\{b\in B(\infty)\mid \ep_i(b)=0\},
\qquad
B_i(\infty):=\{b\in B(\infty)\mid \ep_i^*(b)=0\}.
\end{equation}
The \emph{Saito crystal reflection} \cite[Corollary~3.4.8]{Sai94} is the bijection
\begin{equation}\label{eq:saito-reflection}
\mathrm{T}_i^*\colon B_i(\infty)\xrightarrow{\sim} {}_{i}B(\infty),
\qquad
\mathrm{T}_i^*(b):=(\tf_i^*)^{\varphi_i(b)}(\te_i)^{\ep_i(b)}(b).
\end{equation}
For $b\in B(\infty)$, set
\begin{equation}\label{eq:saito-projection-extension}
\pi_i(b):=(\te_i^*)^{\max}(b),\qquad
\widetilde{\mathrm T}_i^*:=\mathrm{T}_i^*\circ\pi_i.
\end{equation}

On $B_\vartheta(\infty)$, the subsets $B_j(\infty)$ and ${}_jB(\infty)$
and the maps $\pi_j$, $\mathrm T_j^*$, and $\widetilde{\mathrm T}_j^*$
are defined by \eqref{eq:saito-domains}--\eqref{eq:saito-projection-extension},
with $\g_\vartheta$ and $j$ in place of $\g$ and $i$.

Let $w\in\weyl$ and $\bfi=(i_1,\dots,i_\ell)\in\rR(w)$.
For $b\in B(\infty)$, define
\begin{equation}\label{eq:pbw-saito-recursion}
b_0=b,\qquad
b_k=\widetilde{\mathrm T}_{i_k}^*(b_{k-1}),\qquad
p_k=\ep_{i_k}^*(b_{k-1})
\qquad (1\le k\le \ell).
\end{equation}
By \cite[Theorem~4.29]{Kimura12}, the PBW parametrization is a bijection
\[
\PBW_{\bfi}\colon B(w)\xrightarrow{\sim}\Z_{\ge0}^{[1,\ell]}.
\]
By \cite[Section~3B2]{K17}, one has
\begin{equation}\label{eq:pbw-saito-coordinates}
\PBW_{\bfi}(b)=(p_1,\dots,p_\ell)
\qquad (b\in B(w)).
\end{equation}

\begin{lemma}[{cf.\ \cite[Proposition~3.12]{K17}}]\label{prop:pbw-saito-formula}
For $b\in B(\infty)$, one has $b\in B(w)$ if and only if
$b_\ell=\one$.
\end{lemma}

\begin{proof}
If $\ell=0$, then $B(1)=\{\one\}$ and $b_\ell=b$.
Suppose that $\ell>0$ and $b\in B(w)$.
For $1\le k\le\ell$, the definition of $\widetilde{\mathrm T}_{i_k}^*$ gives
\[
\wt(b_k)=s_{i_k}\bigl(\wt(b_{k-1})+p_k\alpha_{i_k}\bigr),
\]
and hence
\[
\wt(b)=w\,\wt(b_\ell)-\sum_{k=1}^{\ell}p_k\beta_{\bfi,k}.
\]
On the other hand, \eqref{eq:pbw-saito-coordinates} gives
$\PBW_\bfi(b)=(p_1,\dots,p_\ell)$. The dual PBW monomial
$\mathsf{F}^{\mathrm{up}}_\bfi(p_1,\dots,p_\ell)$ and the upper global basis element
$\Gup(b)$ in \eqref{eq:pbw-datum-leading-term} both have weight
$-\sum_{k=1}^{\ell}p_k\beta_{\bfi,k}$.
Thus $w\,\wt(b_\ell)=0$. Since $\one$ is the only element of
$B(\infty)$ of weight zero, $b_\ell=\one$.

Conversely, suppose that $b_\ell=\one$.
Since $\PBW_\bfi$ is surjective, choose $b'\in B(w)$ such that
$\PBW_\bfi(b')=(p_1,\dots,p_\ell)$.
Define $b'_0=b'$ and
$b'_k=\widetilde{\mathrm T}_{i_k}^*(b'_{k-1})$ for $1\le k\le\ell$.
Since $b'\in B(w)$, the implication
$b\in B(w)\Longrightarrow b_\ell=\one$ and \eqref{eq:pbw-saito-coordinates} give
\[
\ep_{i_k}^*(b'_{k-1})=p_k\quad(1\le k\le\ell),
\qquad b'_\ell=\one=b_\ell.
\]
For $k=\ell,\dots,1$, the inverse of the bijection $\mathrm T_{i_k}^*$ gives
\[
b_{k-1}=(\tf_{i_k}^*)^{p_k}(\mathrm T_{i_k}^*)^{-1}(b_k),
\qquad
b'_{k-1}=(\tf_{i_k}^*)^{p_k}(\mathrm T_{i_k}^*)^{-1}(b'_k).
\]
Induction from $b_\ell=b'_\ell$ yields $b=b'\in B(w)$.
\end{proof}

Equation~\eqref{eq:pbw-saito-coordinates} and
Lemma~\ref{prop:pbw-saito-formula} also hold for the Cartan datum of
$\g_\vartheta$.

Lemma~\ref{lem:orthogonal-commute} records commutation relations for
crystal operators and Saito reflections at orthogonal indices.

\begin{lemma}[{cf.\ \cite[Section~5]{K96};
\cite[Corollary~2.2.2]{K93D}}]\label{lem:orthogonal-commute}
Let $j\in J$, and let $i,i'\in F_j$ be distinct.
\begin{enumerate}[label=\textup{(\roman*)},ref=\roman*,beginpenalty=10000]
\item\label{item:orthogonal-crystal-operators} For
\[
X\in\{\te_i,\tf_i,\te_i^*,\tf_i^*\},\qquad
Y\in\{\te_{i'},\tf_{i'},\te_{i'}^*,\tf_{i'}^*\},
\]
one has $XY=YX$ on $B(\infty)\cup\{0\}$.
If $b\in B(\infty)$ and $Xb\ne0$, then
\begin{equation}\label{eq:orthogonal-crystal-functions}
\ep_{i'}(Xb)=\ep_{i'}(b),\qquad
\ep_{i'}^*(Xb)=\ep_{i'}^*(b),\qquad
\varphi_{i'}(Xb)=\varphi_{i'}(b).
\end{equation}
The identities \eqref{eq:orthogonal-crystal-functions} hold with
$i$ and $i'$ interchanged.
\item\label{item:orthogonal-projections-reflections} The maps $\pi_i$ and $\widetilde{\mathrm T}_i^*$ preserve
$\ep_{i'},\ep_{i'}^*,\varphi_{i'}$, and the maps
$\pi_{i'}$ and $\widetilde{\mathrm T}_{i'}^*$ preserve
$\ep_i,\ep_i^*,\varphi_i$.
The maps $\pi_i$ and $\widetilde{\mathrm T}_i^*$ satisfy
\begin{equation}\label{eq:orthogonal-projections-reflections}
\begin{gathered}
\pi_i\,\pi_{i'}=\pi_{i'}\,\pi_i,
\qquad
\pi_{i'}\widetilde{\mathrm T}_i^*=\widetilde{\mathrm T}_i^*\pi_{i'},\\
\widetilde{\mathrm T}_i^*\widetilde{\mathrm T}_{i'}^*
=\widetilde{\mathrm T}_{i'}^*\widetilde{\mathrm T}_i^*.
\end{gathered}
\end{equation}
\end{enumerate}
In particular, $\pi_{i'}$ and $\widetilde{\mathrm T}_{i'}^*$ preserve
$B_i(\infty)$, and $\pi_i$ and $\widetilde{\mathrm T}_i^*$ preserve
$B_{i'}(\infty)$. One has
\begin{equation}\label{eq:orthogonal-saito-projection}
\begin{gathered}
\mathrm T_i^*\pi_{i'}(b)=\pi_{i'}\mathrm T_i^*(b)
\qquad(b\in B_i(\infty)),\\
\mathrm T_{i'}^*\pi_i(b)=\pi_i\mathrm T_{i'}^*(b)
\qquad(b\in B_{i'}(\infty)),
\end{gathered}
\end{equation}
and
\begin{equation}\label{eq:orthogonal-saito-reflections}
\mathrm T_i^*\mathrm T_{i'}^*(b)=\mathrm T_{i'}^*\mathrm T_i^*(b)
\qquad(b\in B_i(\infty)\cap B_{i'}(\infty)).
\end{equation}
\end{lemma}

\begin{proof}

\textup{(\ref{item:orthogonal-crystal-operators})} For orthogonal indices, the defining relations and crystal operator
formulas in \cite[Section~2.2.4]{Kimura12} give

\[
XY=YX\qquad
\bigl(X\in\{\te_i,\tf_i\},\ Y\in\{\te_{i'},\tf_{i'}\}\bigr)
\]

on $B(\infty)\cup\{0\}$.
Conjugating $XY=YX$ by $*$ gives $XY=YX$ on $B(\infty)\cup\{0\}$ for
$X\in\{\te_i^*,\tf_i^*\}$ and $Y\in\{\te_{i'}^*,\tf_{i'}^*\}$.

For the mixed commutations, let
$\mathcal B_{i'}=\{b_{i'}(t)\mid t\in\Z\}$ be the elementary crystal with
\[
\begin{gathered}
\wt(b_{i'}(t))=t\alpha_{i'},\qquad
\ep_{i'}(b_{i'}(t))=-t,\qquad \varphi_{i'}(b_{i'}(t))=t,\\
\te_{i'}b_{i'}(t)=b_{i'}(t+1),\qquad
\tf_{i'}b_{i'}(t)=b_{i'}(t-1).
\end{gathered}
\]
For $r\ne i'$ and $t\in\Z$, one has
$\ep_r(b_{i'}(t))=\varphi_r(b_{i'}(t))=-\infty$ and
$\te_r b_{i'}(t)=\tf_r b_{i'}(t)=0$.
Kashiwara's embedding theorem \cite[Theorem~2.2.1]{K93D} gives an
injective map
\[
\Psi_{i'}\colon B(\infty)\longrightarrow B(\infty)\otimes\mathcal B_{i'},
\qquad
\Psi_{i'}(b)=\pi_{i'}(b)\otimes b_{i'}(-\ep_{i'}^*(b)).
\]
 The map $\Psi_{i'}$ commutes with $\te_h$ and $\tf_h$
for every $h\in I$.
We set $\Psi_{i'}(0)=0$ and a tensor with a zero factor equal to $0$.
Write $b'=\pi_{i'}(b)$ and $n=\ep_{i'}^*(b)$.
For $b'\otimes b_{i'}(-n)$, the tensor product rule applies $\te_i$
to $b'$ when $\varphi_i(b')\ge\ep_i(b_{i'}(-n))$, and $\tf_i$ to $b'$
when $\varphi_i(b')>\ep_i(b_{i'}(-n))$ \cite[(2.33)]{Kimura12}.
Since $i\ne i'$ and $\ep_i(b_{i'}(-n))=-\infty$, one has
$\varphi_i(b')>\ep_i(b_{i'}(-n))$.
Thus, by \cite[Theorem~2.36\textup{(2)}]{Kimura12},
\begin{equation}\label{eq:embedding-orthogonal-actions}
\begin{aligned}
\Psi_{i'}(Xb)&=(Xb')\otimes b_{i'}(-n)
&&\bigl(X\in\{\te_i,\tf_i\}\bigr),\\
\Psi_{i'}(\te_{i'}^*b)&=
\begin{cases}b'\otimes b_{i'}(1-n)&n>0,\\0&n=0,\end{cases}\\
\Psi_{i'}(\tf_{i'}^*b)&=b'\otimes b_{i'}(-n-1).
\end{aligned}
\end{equation}
In \eqref{eq:embedding-orthogonal-actions}, $X\in\{\te_i,\tf_i\}$
changes only $b'$, whereas $\te_{i'}^*$ and $\tf_{i'}^*$ change only $n$.
The injectivity of $\Psi_{i'}$ gives $XY=YX$ for
$X\in\{\te_i,\tf_i\}$ and $Y\in\{\te_{i'}^*,\tf_{i'}^*\}$,
including zero values. Interchanging $i$ and $i'$ gives $XY=YX$ for
$X\in\{\te_i^*,\tf_i^*\}$ and $Y\in\{\te_{i'},\tf_{i'}\}$.

 Let $X\in\{\te_i,\tf_i,\te_i^*,\tf_i^*\}$ and $b\in B(\infty)$
satisfy $Xb\ne0$.
Writing $a=\ep_{i'}(b)$, we have
\[
\te_{i'}^{a+1}Xb
=X\te_{i'}^{a+1}b=0,
\]
and hence $\ep_{i'}(Xb)\le\ep_{i'}(b)$.
Applying the bound $\ep_{i'}(Xb)\le\ep_{i'}(b)$ to the partial
inverse of $X$ on $Xb$ gives $\ep_{i'}(b)\le\ep_{i'}(Xb)$.
Since $X\te_{i'}^*=\te_{i'}^*X$, one has
$(\te_{i'}^*)^{\ep_{i'}^*(b)+1}Xb
=X(\te_{i'}^*)^{\ep_{i'}^*(b)+1}b=0$,
so $\ep_{i'}^*(Xb)\le\ep_{i'}^*(b)$.
Applying $\ep_{i'}^*(Xb)\le\ep_{i'}^*(b)$ to the partial inverse of
$X$ on $Xb$ gives $\ep_{i'}^*(b)\le\ep_{i'}^*(Xb)$.
Hence $\ep_{i'}^*(Xb)=\ep_{i'}^*(b)$.

Since $\wt(Xb)-\wt(b)=\pm\alpha_i$ and $\langle h_{i'},\alpha_i\rangle=0$,
\[
\begin{aligned}
\varphi_{i'}(Xb)
&=\ep_{i'}(Xb)+\langle h_{i'},\wt(Xb)\rangle\\
&=\ep_{i'}(b)+\langle h_{i'},\wt(b)\rangle
=\varphi_{i'}(b).
\end{aligned}
\]
Thus Lemma~\ref{lem:orthogonal-commute}\textup{(\ref{item:orthogonal-crystal-operators})} holds.

\textup{(\ref{item:orthogonal-projections-reflections})} The maps $\pi_i$ and $\widetilde{\mathrm T}_i^*$ preserve
$\ep_{i'},\ep_{i'}^*$, and $\varphi_{i'}$ by Lemma~\ref{lem:orthogonal-commute}\textup{(\ref{item:orthogonal-crystal-operators})}.
Equations~\eqref{eq:saito-reflection} and \eqref{eq:saito-projection-extension}
express $\pi_i(b)$ and $\widetilde{\mathrm T}_i^*(b)$ as compositions of
$\te_i,\tf_i,\te_i^*,\tf_i^*$ applied to $b$.
Since $\pi_i(b)\ne0$ and $\widetilde{\mathrm T}_i^*(b)\ne0$,
Lemma~\ref{lem:orthogonal-commute}\textup{(\ref{item:orthogonal-crystal-operators})} applies to each operator
in \eqref{eq:saito-reflection} and \eqref{eq:saito-projection-extension}.
The maps $\pi_{i'}$ and $\widetilde{\mathrm T}_{i'}^*$ preserve
$\ep_i,\ep_i^*$, and $\varphi_i$.
In particular, $\pi_{i'}$ and $\widetilde{\mathrm T}_{i'}^*$ preserve
$B_i(\infty)$, and $\pi_i$ and $\widetilde{\mathrm T}_i^*$ preserve
$B_{i'}(\infty)$.

For $X\in\{\te_i,\tf_i,\te_i^*,\tf_i^*\}$ with $Xb\ne0$, we also have
\begin{equation}\label{eq:orthogonal-projection-operator}
\pi_{i'}(Xb)
=(\te_{i'}^*)^{\ep_{i'}^*(b)}Xb
=X(\te_{i'}^*)^{\ep_{i'}^*(b)}b
=X\pi_{i'}(b).
\end{equation}
Applying \eqref{eq:orthogonal-projection-operator} repeatedly and
using preservation of $\ep_i^*$ by $\pi_{i'}$ gives
$\pi_i\pi_{i'}=\pi_{i'}\pi_i$.
Now put $d=\pi_i(b)$.
Since $\pi_{i'}$ preserves $\ep_i$ and $\varphi_i$,
\eqref{eq:orthogonal-projection-operator} gives
\[
\begin{aligned}
\pi_{i'}\widetilde{\mathrm T}_i^*(b)
&=
(\tf_i^*)^{\varphi_i(d)}
\te_i^{\ep_i(d)}\pi_{i'}(d)\\
&=
\widetilde{\mathrm T}_i^*\pi_{i'}(b),
\end{aligned}
\]
where $\pi_{i'}(d)=\pi_i\pi_{i'}(b)$ by
$\pi_i\pi_{i'}=\pi_{i'}\pi_i$.

For the commutation of $\widetilde{\mathrm T}_i^*$ and
$\widetilde{\mathrm T}_{i'}^*$, set $c=\pi_i\pi_{i'}(b)$.
Then $c\in B_i(\infty)\cap B_{i'}(\infty)$.
The preservation of $B_i(\infty)$ by $\widetilde{\mathrm T}_{i'}^*$
and of $B_{i'}(\infty)$ by $\widetilde{\mathrm T}_i^*$ ensures that
$\mathrm T_i^*\mathrm T_{i'}^*(c)$ and
$\mathrm T_{i'}^*\mathrm T_i^*(c)$ are defined.
Moreover, $\mathrm T_{i'}^*$ preserves $\ep_i,\varphi_i$,
and $\mathrm T_i^*$ preserves $\ep_{i'},\varphi_{i'}$.

We obtain
\[
\begin{aligned}
\widetilde{\mathrm T}_i^*
\widetilde{\mathrm T}_{i'}^*(b)
&=\mathrm T_i^*\mathrm T_{i'}^*(c)\\
&=
(\tf_i^*)^{\varphi_i(c)}\te_i^{\ep_i(c)}
(\tf_{i'}^*)^{\varphi_{i'}(c)}
\te_{i'}^{\ep_{i'}(c)}c\\
&=
(\tf_{i'}^*)^{\varphi_{i'}(c)}
\te_{i'}^{\ep_{i'}(c)}
(\tf_i^*)^{\varphi_i(c)}\te_i^{\ep_i(c)}c\\
&=\mathrm T_{i'}^*\mathrm T_i^*(c)\\
&=\widetilde{\mathrm T}_{i'}^*
\widetilde{\mathrm T}_i^*(b).
\end{aligned}
\]
Thus Lemma~\ref{lem:orthogonal-commute}\textup{(\ref{item:orthogonal-projections-reflections})} holds.
Restricting \eqref{eq:orthogonal-projections-reflections} gives
\eqref{eq:orthogonal-saito-projection} and
\eqref{eq:orthogonal-saito-reflections}, since
$\widetilde{\mathrm T}_i^*=\mathrm T_i^*$ on $B_i(\infty)$ and
$\widetilde{\mathrm T}_{i'}^*=\mathrm T_{i'}^*$ on $B_{i'}(\infty)$.
\end{proof}

\begin{proposition}\label{prop:saito-similarity}
For every $j\in J$ and every $b\in B_\vartheta(\infty)$, one has
\[
\Upsilon_\vartheta\bigl(\widetilde{\mathrm T}_j^*(b)\bigr)
=\Bigl(\prod_{u=1}^{d_j}\widetilde{\mathrm T}_{j^{(u)}}^*\Bigr)
\Upsilon_\vartheta(b).
\]
\end{proposition}

\begin{proof}

Put $c=\pi_j(b)$. By \eqref{eq:similarity-epsilon-fiber},
\[
m_i\ep_j^*(b)=\ep_i^*(\Upsilon_\vartheta(b))\qquad(i\in F_j).
\]
Lemma~\ref{lem:orthogonal-commute}\textup{(\ref{item:orthogonal-projections-reflections})} gives
$\ep_i^*(\pi_{i'}x)=\ep_i^*(x)$ for distinct $i,i'\in F_j$ and
$x\in B(\infty)$.
By \eqref{eq:similarity-e-f-star}, \eqref{eq:similarity-epsilon-fiber},
and Lemma~\ref{lem:orthogonal-commute}\textup{(\ref{item:orthogonal-projections-reflections})},
\begin{equation}\label{eq:saito-similarity-projection-mi}
\begin{aligned}
\Upsilon_\vartheta(c)
&=\Bigl(\prod_{i\in F_j}(\te_i^*)^{m_i\ep_j^*(b)}\Bigr)
  \Upsilon_\vartheta(b)\\
&=\Bigl(\prod_{i\in F_j}(\te_i^*)^{\ep_i^*(\Upsilon_\vartheta(b))}\Bigr)
  \Upsilon_\vartheta(b)
=\Bigl(\prod_{i\in F_j}\pi_i\Bigr)\Upsilon_\vartheta(b).
\end{aligned}
\end{equation}
Since $\ep_i^*(\Upsilon_\vartheta(c))=m_i\ep_j^*(c)=0$ for
$i\in F_j$, by \eqref{eq:similarity-epsilon-fiber}, we have
\[
\Upsilon_\vartheta(c)\in\bigcap_{i\in F_j}B_i(\infty).
\]

For $i\in F_j$, \eqref{eq:similarity-epsilon-fiber} and
\eqref{eq:similarity-phi-fiber} give
\[
\ep_i(\Upsilon_\vartheta(c))=m_i\ep_j(c),\qquad
\varphi_i(\Upsilon_\vartheta(c))=m_i\varphi_j(c).
\]
For distinct $i,i'\in F_j$, the reflection $\mathrm T_{i'}^*$
preserves $B_i(\infty)$, $\ep_i$, and $\varphi_i$ by
Lemma~\ref{lem:orthogonal-commute}.
Hence $\prod_{i\in F_j}\mathrm T_i^*$ is defined on
$\Upsilon_\vartheta(c)$.
Using \eqref{eq:saito-reflection},
Lemma~\ref{lem:orthogonal-commute}\textup{(\ref{item:orthogonal-crystal-operators})}, and
\eqref{eq:similarity-e-f}--\eqref{eq:similarity-e-f-star}, we obtain
\[
\begin{aligned}
\Bigl(\prod_{i\in F_j}\mathrm T_i^*\Bigr)\Upsilon_\vartheta(c)
&=\Bigl(\prod_{i\in F_j}(\tf_i^*)^{m_i\varphi_j(c)}\Bigr)
  \Bigl(\prod_{i\in F_j}\te_i^{m_i\ep_j(c)}\Bigr)\Upsilon_\vartheta(c)\\
&=\Upsilon_\vartheta\bigl((\tf_j^*)^{\varphi_j(c)}
  \te_j^{\ep_j(c)}c\bigr)
=\Upsilon_\vartheta(\mathrm T_j^*c).
\end{aligned}
\]
For distinct $i,i'\in F_j$, Lemma~\ref{lem:orthogonal-commute} also gives
$\mathrm T_i^*\pi_{i'}=\pi_{i'}\mathrm T_i^*$ on $B_i(\infty)$.
Using \eqref{eq:orthogonal-saito-projection} and
\eqref{eq:saito-similarity-projection-mi}, we obtain
\[
\begin{aligned}
\Upsilon_\vartheta(\widetilde{\mathrm T}_j^*b)
&=\Bigl(\prod_{i\in F_j}\mathrm T_i^*\Bigr)
  \Bigl(\prod_{i\in F_j}\pi_i\Bigr)\Upsilon_\vartheta(b)\\
&=\Bigl(\prod_{i\in F_j}(\mathrm T_i^*\circ\pi_i)\Bigr)\Upsilon_\vartheta(b)
=\Bigl(\prod_{u=1}^{d_j}\widetilde{\mathrm T}_{j^{(u)}}^*\Bigr)
  \Upsilon_\vartheta(b).
\end{aligned}
\]
 \end{proof}

\subsection{The PBW and string formulas}\label{subsec:pbw-string-formulas}

Fix $w\in\weyl_\vartheta$ with $\supp(w)=J$ and a reduced expression

\begin{equation}\label{eq:bfi-vartheta}
\bfj=(j_1,\dots,j_\ell)\in\rR(w),\qquad j_k\in J.
\end{equation}
Concatenating
$(j_k^{(1)},\dots,j_k^{(d_{j_k})})$ for $1\le k\le\ell$ gives
$\bfi=(i_1,\dots,i_L)$ with $i_r\in I$, where
\begin{equation}\label{eq:ell-bar}
L=\sum_{k=1}^{\ell}d_{j_k}=\ell(\Theta(w)).
\end{equation}
Set $L_0=0$ and $L_k=\sum_{t=1}^{k}d_{j_t}$ for $1\le k\le\ell$.
For $1\le k\le\ell$, define
\begin{equation}\label{eq:block-interval}
D_k=\{L_{k-1}+1,\dots,L_k\}.
\end{equation}
Thus $D_k$ is the $k$-th block of coordinate positions in $\bfi$.
We define the $\Z$-linear map
\begin{equation}\label{eq:upsilon-def}
\Gamma_{\bfj}\colon\Z^{[1,\ell]}
\longrightarrow\Z^{[1,L]}
\end{equation}
by
\begin{equation}\label{eq:upsilon-formula}
\bigl(\Gamma_{\bfj}(a_1,\dots,a_\ell)\bigr)_r=m_{i_r}a_k
\qquad(1\le k\le\ell,\ r\in D_k).
\end{equation}

\begin{example}\label{ex:similarity-C2}
Let $\g$ have type $A_3$, with consecutive vertex labels $1,2,3$, and take
\[
F_1=\{1,3\},\qquad F_2=\{2\},\qquad
m_1=m_3=1,\quad m_2=2.
\]
The algebra $\g_\vartheta$ has type $C_2$.
Order $F_1$ so that $1$ precedes $3$.
The reduced expression $\bfj=(1,2)$ then expands to
$\bfi=(1,3,2)$, and
\[
\Gamma_{\bfj}(a_1,a_2)=(a_1,a_1,2a_2).
\]
\end{example}

\begin{theorem}\label{thm:similarity-Bw}
\label{prop:restriction-Bw}\label{thm:PBW-Bw}
For every $b\in B_\vartheta(w)$, one has
$\Upsilon_\vartheta(b)\in B(\Theta(w))$ and
\begin{equation}\label{eq:PBW-Bw}
\PBW_{\bfi}\bigl(\Upsilon_\vartheta(b)\bigr)
=\Gamma_{\bfj}\bigl(\PBW_\bfj(b)\bigr).
\end{equation}
\end{theorem}

\begin{proof}
Let $b\in B_\vartheta(w)$, and define
\[
\begin{aligned}
b_0&=b,&
b_k&=\widetilde{\mathrm T}_{j_k}^*(b_{k-1})
,\qquad p_k=\ep_{j_k}^*(b_{k-1})
&&(1\le k\le\ell),\\
\overline b_0&=\Upsilon_\vartheta(b),&
\overline b_r&=\widetilde{\mathrm T}_{i_r}^*(\overline b_{r-1})
,\qquad \overline p_r=\ep_{i_r}^*(\overline b_{r-1})
&&(1\le r\le L).
\end{aligned}
\]
By \eqref{eq:pbw-saito-coordinates} and Lemma~\ref{prop:pbw-saito-formula},
\[
\PBW_\bfj(b)=(p_1,\dots,p_\ell),\qquad b_\ell=\one.
\]

We prove by induction on $k$ that
\[
\overline b_{L_k}=\Upsilon_\vartheta(b_k)\quad(0\le k\le\ell),
\qquad
\overline p_r=m_{i_r}p_k\quad(1\le k\le\ell,\ r\in D_k).
\]
The initial equality is $\overline b_0=\Upsilon_\vartheta(b_0)$.
Suppose that $\overline b_{L_{k-1}}=\Upsilon_\vartheta(b_{k-1})$.
For $r\in D_k$, the recursion for $\overline b_r$ gives
\[
\overline b_{r-1}
=\widetilde{\mathrm T}_{i_{r-1}}^*\cdots
 \widetilde{\mathrm T}_{i_{L_{k-1}+1}}^*
 \bigl(\overline b_{L_{k-1}}\bigr).
\]
For $L_{k-1}+1\le s<r$, the indices $i_s$ and $i_r$ are
distinct elements of $F_{j_k}$, so
Lemma~\ref{lem:orthogonal-commute}\textup{(\ref{item:orthogonal-projections-reflections})} shows that
$\widetilde{\mathrm T}_{i_s}^*$ preserves $\ep_{i_r}^*$.
Using $\overline b_{L_{k-1}}=\Upsilon_\vartheta(b_{k-1})$ and
\eqref{eq:similarity-epsilon-fiber}, we obtain
\[
\begin{aligned}
\overline p_r
&=\ep_{i_r}^*(\overline b_{r-1})
=\ep_{i_r}^*(\overline b_{L_{k-1}})\\
&=m_{i_r}\ep_{j_k}^*(b_{k-1})=m_{i_r}p_k.
\end{aligned}
\]
Proposition~\ref{prop:saito-similarity} gives
\[
\begin{aligned}
\overline b_{L_k}
&=\Bigl(\prod_{u=1}^{d_{j_k}}\widetilde{\mathrm T}_{j_k^{(u)}}^*\Bigr)
  \Upsilon_\vartheta(b_{k-1})\\
&=\Upsilon_\vartheta(\widetilde{\mathrm T}_{j_k}^*b_{k-1})
=\Upsilon_\vartheta(b_k).
\end{aligned}
\]

Thus $\overline b_L=\Upsilon_\vartheta(b_\ell)=\one$, so
Lemma~\ref{prop:pbw-saito-formula} implies
$\Upsilon_\vartheta(b)\in B(\Theta(w))$.
Equations~\eqref{eq:pbw-saito-coordinates} and \eqref{eq:upsilon-formula}
now give
\[
\PBW_\bfi(\Upsilon_\vartheta(b))=(\overline p_1,\dots,\overline p_L)
=\Gamma_\bfj(p_1,\dots,p_\ell)
=\Gamma_\bfj(\PBW_\bfj(b)).
\]
\end{proof}

\begin{proposition}\label{prop:string-Bw}
For every $b\in B_\vartheta(w)$, one has
\[
\STR_{\bfi}\bigl(\Upsilon_\vartheta(b)\bigr)
=\Gamma_{\bfj}\bigl(\STR_\bfj(b)\bigr).
\]
\end{proposition}

\begin{proof}
Let $b\in B_\vartheta(w)$.
By Theorem~\ref{thm:similarity-Bw},
$\Upsilon_\vartheta(b)\in B(\Theta(w))$.
Define
\[
\begin{aligned}
b_0&=b,&
b_k&=\te_{j_k}^{\max}(b_{k-1}),\qquad
t_k=\ep_{j_k}(b_{k-1})
&&(1\le k\le\ell),\\
\overline b_0&=\Upsilon_\vartheta(b),&
\overline b_r&=\te_{i_r}^{\max}(\overline b_{r-1}),\qquad
\overline t_r=\ep_{i_r}(\overline b_{r-1})
&&(1\le r\le L).
\end{aligned}
\]
By \eqref{eq:STR-on-Bw},
\[
\STR_\bfj(b)=(t_1,\dots,t_\ell),\qquad
\STR_\bfi(\Upsilon_\vartheta(b))=(\overline t_1,\dots,\overline t_L).
\]
Each $\overline b_r$ is nonzero because
$\overline t_r=\ep_{i_r}(\overline b_{r-1})$.

We prove by induction on $k$ that
\[
\overline b_{L_k}=\Upsilon_\vartheta(b_k)\quad(0\le k\le\ell),
\qquad
\overline t_r=m_{i_r}t_k\quad(1\le k\le\ell,\ r\in D_k).
\]
The initial equality is $\overline b_0=\Upsilon_\vartheta(b_0)$.
Suppose that $\overline b_{L_{k-1}}=\Upsilon_\vartheta(b_{k-1})$.
For $r\in D_k$, the indices
$i_{L_{k-1}+1},\dots,i_{r-1},i_r$ are distinct elements of $F_{j_k}$.
Lemma~\ref{lem:orthogonal-commute}\textup{(\ref{item:orthogonal-crystal-operators})} and
\eqref{eq:similarity-epsilon-fiber} therefore give
\[
\begin{aligned}
\overline t_r
&=\ep_{i_r}(\overline b_{r-1})
=\ep_{i_r}(\overline b_{L_{k-1}})\\
&=m_{i_r}\ep_{j_k}(b_{k-1})=m_{i_r}t_k.
\end{aligned}
\]
Equation~\eqref{eq:similarity-e-f} gives
\[
\begin{aligned}
\overline b_{L_k}
&=\Bigl(\prod_{i\in F_{j_k}}\te_i^{m_i t_k}\Bigr)
  \Upsilon_\vartheta(b_{k-1})\\
&=\Upsilon_\vartheta(\te_{j_k}^{t_k}b_{k-1})
=\Upsilon_\vartheta(b_k).
\end{aligned}
\]

Equation~\eqref{eq:upsilon-formula} now gives
\[
\STR_\bfi(\Upsilon_\vartheta(b))=(\overline t_1,\dots,\overline t_L)
=\Gamma_\bfj(t_1,\dots,t_\ell)
=\Gamma_\bfj(\STR_\bfj(b)).
\]
\end{proof}

\section{The localized similarity map}\label{sec:similarity-localized}\label{sec:remaining-linearity}

We construct $\widetilde{\Upsilon}_\vartheta$ and compute the PBW data,
string data, and left and right $g$-vectors of
$\widetilde{\Upsilon}_\vartheta(x)$ for $x\in\LB_\vartheta(w)$.
We also prove analogues of Proposition~\ref{prop:similarity-map}
for the crystal structures on $\LB_\vartheta(w)$ and
$\LB(\Theta(w))$ (Theorem~\ref{thm:localized-crystal-similarity}).
For $w\in\weyl$ and $\bfi\in\rR(w)$, we first recall $\LB(w)$ and
the extensions of $\PBW_\bfi$ and $\STR_\bfi$ from $B(w)$ to $\LB(w)$.

\subsection{Localized crystals}\label{sec:prelim-crystals}

For the rest of the paper, we assume $\supp(w)=I$ for
localizations associated with $w\in\weyl$, and $\supp(w)=J$
for localizations associated with $w\in\weyl_\vartheta$.
Fix $w\in\weyl$.
For each $i\in I$, let
\begin{equation}\label{eq:frozen-quantum-minor}
\qmC_i^w:=\qmD(w\Lambda_i,\Lambda_i)\in \Gup(w)
\end{equation}
be the quantum unipotent minor.
By \cite[Corollary~6.18]{Kimura12}, each $\qmC_i^w$ commutes with every
homogeneous element of $\qA_q(\n(w))$ up to a power of $q$.
We define the localized quantum unipotent coordinate ring by
\begin{equation}\label{eq:localized-algebra}
\widetilde{\qA}_q(\n(w))
=\qA_q(\n(w))[(\qmC_i^w)^{-1}\mid i\in I].
\end{equation}
For $\lambda\in P^+$, the minor $\qmD(w\lambda,\lambda)$ is, up to a
power of $q$, a product of the $\qmC_i^w$, and hence is invertible
in $\widetilde{\qA}_q(\n(w))$ \cite[Proposition~4.2]{KiOy21}.

The upper global basis of $\widetilde{\qA}_q(\n(w))$ is
\begin{equation}\label{eq:localized-upper-global-basis}
\tGup(w)=
\left\{
q^{(\lambda,\wt(b)+\lambda-w\lambda)}
\qmD(w\lambda,\lambda)^{-1}\Gup(b)
\ \middle|\ \lambda\in P^+,\ b\in B(w)
\right\},
\end{equation}
as defined in \cite[Definition~4.6]{KiOy21}.
For $x,y\in\widetilde{\qA}_q(\n(w))$, we write $x\equiv_q y$
if $x=q^m y$ for some $m\in\Z$.

For each $i\in I$, let $\fz_i^w\in B(w)$ be the crystal element determined by
\begin{equation}\label{eq:frozen-crystal-element}
\Gup(\fz_i^w)=\qmC_i^w.
\end{equation}

For the fixed $w$, we abbreviate $\fz_i^w$ to $\fz_i$.

Following \cite[Section~2]{JP25}, we define an equivalence relation on
$B(w)\times \Z^{\oplus I}$ as follows. For $(b,\bfa),(b',\bfa')\in B(w)\times \Z^{\oplus I}$,
write $\bfa=(a_i)_{i\in I}$ and $\bfa'=(a_i')_{i\in I}$. We declare that
\begin{equation}\label{eq:localized-equiv}
(b,\bfa)\sim (b',\bfa')
\end{equation}
if there exists $N\in \Z_{>0}$ such that $N+\min_{i\in I}\{a_i,a_i'\}\ge 0$ and
\begin{equation}\label{eq:localized-equiv-product}
\Gup(b)\prod_{i\in I}\Gup(\fz_i)^{N+a_i}
\equiv_q
\Gup(b')\prod_{i\in I}\Gup(\fz_i)^{N+a_i'}.
\end{equation}

Following \cite[Definition~2.3]{JP25}, we call the quotient
\begin{equation}\label{eq:localized-crystal-vartheta}
\LB(w)=\bigl(B(w)\times\Z^{\oplus I}\bigr)/\sim
\end{equation}
the localized crystal of $B(w)$.
The set $\LB(w)$ admits a crystal structure by
\cite[Theorem~6.13]{KN25} and the identification in \cite[(2.6)]{JP25}.
For $b\in B(w)$ and $\bfa=(a_i)_{i\in I}\in \Z^{\oplus I}$, we write the class of $(b,\bfa)$ as
\begin{equation}\label{eq:localized-notation}
x=b\cdot \prod_{i\in I}(\fz_i^w)^{a_i}=b\cdot \fz(\bfa).
\end{equation}
We identify $b\in B(w)$ with $b\cdot\fz(\mathbf0)\in\LB(w)$
and abbreviate $\one\cdot\fz(\bfa)$ to $\fz(\bfa)$, or to
$\prod_{i\in I}(\fz_i^w)^{a_i}$.

For $\mathbf d\in\Z^{\oplus I}$, we set
\[
\bigl(b\cdot\fz(\bfa)\bigr)\cdot\fz(\mathbf d)
:=b\cdot\fz(\bfa+\mathbf d).
\]
If $(b,\bfa)\sim(b',\bfa')$, then
$(b,\bfa+\mathbf d)\sim(b',\bfa'+\mathbf d)$, so
$x\cdot\fz(\mathbf d)$ is independent of the representative of $x$.

The bijection $\tGup\colon\LB(w)\to\tGup(w)$ is characterized by
\[
\tGup\bigl(b\cdot\fz(\bfa)\bigr)
\equiv_q \Gup(b)\prod_{i\in I}(\qmC_i^w)^{a_i};
\]
see \cite[Section~2]{JP25}.

Let $\bfi=(i_1,\dots,i_\ell)\in\rR(w)$.
For $b\in B(w)$ and $\bfa=(a_i)_{i\in I}\in\Z^{\oplus I}$, the PBW and string data on the
localized crystal are given, as in \cite[Section~2]{JP25}, by
\begin{align}
\STR_{\bfi}\bigl(b\cdot \fz(\bfa)\bigr)
&=
\STR_{\bfi}(b)+\sum_{i\in I}a_i\,\STR_{\bfi}(\fz_i^w), \label{eq:localized-affine-string}\\
\PBW_{\bfi}\bigl(b\cdot \fz(\bfa)\bigr)
&=
\PBW_{\bfi}(b)+\sum_{i\in I}a_i\,\PBW_{\bfi}(\fz_i^w). \label{eq:localized-affine-pbw}
\end{align}
Set
\begin{equation}\label{eq:lattices-vartheta}
\tcP_{\bfi}(w):=
\PBW_{\bfi}(B(w))+\sum_{i\in I}\Z\,\PBW_{\bfi}(\fz_i^w),
\qquad
\tcS_{\bfi}(w):=
\STR_{\bfi}(B(w))+\sum_{i\in I}\Z\,\STR_{\bfi}(\fz_i^w).
\end{equation}
By \cite[Sections~2.1--2.2]{JP25}, the right-hand sides of
\eqref{eq:localized-affine-string} and \eqref{eq:localized-affine-pbw}
depend only on the class of $(b,\bfa)$ in $\LB(w)$.
The extensions of $\STR_\bfi$ and $\PBW_\bfi$
defined by \eqref{eq:localized-affine-string} and
\eqref{eq:localized-affine-pbw} give bijections onto
$\tcS_{\bfi}(w)$ and $\tcP_{\bfi}(w)$, respectively.
Moreover, \cite[Remark~3.2(iii)]{JP25} gives
\begin{equation}\label{eq:g-vector-lattices-vartheta}
\tcP_{\bfi}(w)=\tcS_{\bfi}(w)=\Z^{[1,\ell]}.
\end{equation}
We write
\begin{equation}\label{eq:four-param-vartheta}
\pbwP_{\bfi}\colon \LB(w)\xrightarrow{\sim}\tcP_{\bfi}(w),
\quad x\longmapsto \PBW_{\bfi}(x),
\qquad
\strS_{\bfi}\colon \LB(w)\xrightarrow{\sim}\tcS_{\bfi}(w),
\quad x\longmapsto \STR_{\bfi}(x).
\end{equation}

For $w\in\weyl_\vartheta$, we write
$\LB_\vartheta(w)$ for the localization of $B_\vartheta(w)$ defined
by \eqref{eq:localized-equiv}--\eqref{eq:localized-crystal-vartheta}, with
$I$ replaced by $J$.
The PBW and string parametrizations of $\LB_\vartheta(w)$ are
defined using the Cartan datum of $\g_\vartheta$
and a reduced expression of $w$ in $\weyl_\vartheta$.
For $j\in J$, the corresponding quantum minor and crystal element are
denoted by $\qmC_j^w$ and $\fz_j^w\in B_\vartheta(w)$, respectively; thus
\[
\qmC_j^w=\qmD(w\Lambda_j^\vartheta,\Lambda_j^\vartheta),
\qquad
\Gup(\fz_j^w)=\qmC_j^w.
\]

\subsection{The extension to \texorpdfstring{$\LB_\vartheta(w)$}{LB-vartheta(w)}}\label{subsec:localized-similarity}

For the rest of Section~\ref{sec:similarity-localized}, fix $w\in\weyl_\vartheta$
and $\bfj=(j_1,\dots,j_\ell)\in\rR(w)$.
Let $\bfi$ be the expanded word, $D_k$ the blocks of positions in
$\bfi$, and $\Gamma_{\bfj}$ the linear map defined in
Section~\ref{subsec:pbw-string-formulas}.

\begin{lemma}\label{lem:frozen-image}
For every $j\in J$, the equality
\[
\Upsilon_\vartheta(\fz_j^w)=\prod_{i\in F_j}
\bigl(\fz_i^{\Theta(w)}\bigr)^{m_i}
\]
holds in $\LB(\Theta(w))$.
\end{lemma}

\begin{proof}
Fix $j\in J$. The PBW data of $\fz_j^w$ and $\fz_i^{\Theta(w)}$ are
\cite[Corollary~6.4]{Kimura12}
\[
\begin{aligned}
\bigl(\PBW_\bfj(\fz_j^w)\bigr)_k&=\delta_{j_k,j}
&& (1\le k\le\ell),\\
\bigl(\PBW_{\bfi}(\fz_i^{\Theta(w)})\bigr)_r&=\delta_{i_r,i}
&& (i\in I,\ 1\le r\le L).
\end{aligned}
\]
For $r\in D_k$, one has $i_r\in F_{j_k}$, so
\[
\begin{aligned}
\bigl(\Gamma_{\bfj}(\PBW_\bfj(\fz_j^w))\bigr)_r
&=m_{i_r}\delta_{j_k,j}
=\sum_{i\in F_j}m_i\delta_{i_r,i}\\
&=\Bigl(\sum_{i\in F_j}m_i
  \PBW_{\bfi}(\fz_i^{\Theta(w)})\Bigr)_r.
\end{aligned}
\]
By \eqref{eq:PBW-Bw} and \eqref{eq:localized-affine-pbw},
\begin{align*}
\PBW_{\bfi}\bigl(\Upsilon_\vartheta(\fz_j^w)\bigr)
&=\Gamma_{\bfj}\bigl(\PBW_{\bfj}(\fz_j^w)\bigr)\\
&=\sum_{i\in F_j}m_i\PBW_{\bfi}(\fz_i^{\Theta(w)})\\
&=\PBW_{\bfi}\Bigl(\prod_{i\in F_j}\bigl(\fz_i^{\Theta(w)}\bigr)^{m_i}\Bigr).
\end{align*}
The injectivity of $\pbwP_{\bfi}$ gives
$\Upsilon_\vartheta(\fz_j^w)=
\prod_{i\in F_j}(\fz_i^{\Theta(w)})^{m_i}$.
\end{proof}

For $\bfa=(a_j)_{j\in J}\in \Z^{\oplus J}$, set
\begin{equation}\label{eq:y-a}
y(\bfa):=\prod_{j\in J}\prod_{i\in F_j}\bigl(\fz_i^{\Theta(w)}\bigr)^{m_i a_j}.
\end{equation}

\begin{proposition}\label{prop:S-extends-to-localized}
There is a unique map
\[
\widetilde{\Upsilon}_\vartheta\colon \LB_\vartheta(w)\longrightarrow \LB(\Theta(w))
\]
such that
\begin{equation}\label{eq:localized-S-formula}
\widetilde{\Upsilon}_\vartheta\bigl(b\cdot\fz(\bfa)\bigr)
=\Upsilon_\vartheta(b)\cdot y(\bfa)
\qquad (b\in B_\vartheta(w),\ \bfa\in\Z^{\oplus J}).
\end{equation}
The map $\widetilde{\Upsilon}_\vartheta$ is injective and satisfies
\begin{equation}\label{eq:PBW-localized}
\PBW_{\bfi}\bigl(\widetilde{\Upsilon}_\vartheta(x)\bigr)
=\Gamma_{\bfj}\bigl(\PBW_{\bfj}(x)\bigr)
\qquad(x\in\LB_\vartheta(w)).
\end{equation}
\end{proposition}

\begin{proof}
For $b\in B_\vartheta(w)$ and
$\bfa=(a_j)_{j\in J}\in\Z^{\oplus J}$, we compute the PBW datum of
$\Upsilon_\vartheta(b)\cdot y(\bfa)$.
By Theorem~\ref{prop:restriction-Bw},
$\Upsilon_\vartheta(b)\in B(\Theta(w))$.
The definitions
\eqref{eq:y-a} and \eqref{eq:localized-crystal-vartheta}, applied to
$\Theta(w)$, therefore show that
$\Upsilon_\vartheta(b)\cdot y(\bfa)\in\LB(\Theta(w))$.
By \eqref{eq:localized-affine-pbw}, Theorem~\ref{thm:PBW-Bw}, and
Lemma~\ref{lem:frozen-image}, we compute
\begin{equation}\label{eq:localized-image-pbw}
\begin{aligned}
\PBW_{\bfi}\bigl(\Upsilon_\vartheta(b)\cdot y(\bfa)\bigr)
&=\PBW_{\bfi}\bigl(\Upsilon_\vartheta(b)\bigr)
  +\sum_{j\in J}\sum_{i\in F_j}m_i a_j
  \PBW_{\bfi}\bigl(\fz_i^{\Theta(w)}\bigr)\\
&=\Gamma_{\bfj}\bigl(\PBW_{\bfj}(b)\bigr)
  +\sum_{j\in J}a_j\,
  \Gamma_{\bfj}\bigl(\PBW_{\bfj}(\fz_j^w)\bigr)\\
&=\Gamma_{\bfj}\left(
  \PBW_{\bfj}(b)+\sum_{j\in J}a_j\PBW_{\bfj}(\fz_j^w)
  \right)\\
&=\Gamma_{\bfj}\bigl(\PBW_{\bfj}(b\cdot\fz(\bfa))\bigr).
\end{aligned}
\end{equation}

To prove independence of the representative, let
$b'\in B_\vartheta(w)$ and $\bfa'\in\Z^{\oplus J}$ satisfy
\[
b\cdot\fz(\bfa)=b'\cdot\fz(\bfa')
\qquad\text{in }\LB_\vartheta(w).
\]
Applying \eqref{eq:localized-image-pbw} to $b\cdot\fz(\bfa)$ and
$b'\cdot\fz(\bfa')$ gives
\begin{align*}
\PBW_{\bfi}\bigl(\Upsilon_\vartheta(b)\cdot y(\bfa)\bigr)
&=\Gamma_{\bfj}\bigl(\PBW_\bfj(b\cdot\fz(\bfa))\bigr)\\
&=\PBW_{\bfi}\bigl(\Upsilon_\vartheta(b')\cdot y(\bfa')\bigr).
\end{align*}
Since $\pbwP_{\bfi}$ is injective, we obtain
\[
\Upsilon_\vartheta(b)\cdot y(\bfa)
=\Upsilon_\vartheta(b')\cdot y(\bfa').
\]
Thus \eqref{eq:localized-S-formula} is independent of the representative.
Every element of $\LB_\vartheta(w)$ can be written as $b\cdot\fz(\bfa)$,
so \eqref{eq:localized-S-formula} defines a unique map
$\widetilde{\Upsilon}_\vartheta$.
Equation~\eqref{eq:localized-image-pbw} gives \eqref{eq:PBW-localized}.

To prove injectivity, let $x,x'\in\LB_\vartheta(w)$ satisfy
$\widetilde{\Upsilon}_\vartheta(x)=\widetilde{\Upsilon}_\vartheta(x')$.
By \eqref{eq:PBW-localized},
\[
\Gamma_{\bfj}\bigl(\PBW_\bfj(x)\bigr)
=\Gamma_{\bfj}\bigl(\PBW_\bfj(x')\bigr).
\]
Since $F_j\ne\varnothing$ for every $j\in J$ and $m_i>0$ for every
$i\in I$, \eqref{eq:upsilon-formula} makes $\Gamma_{\bfj}$ injective.
Hence $\PBW_\bfj(x)=\PBW_\bfj(x')$.
The injectivity of $\pbwP_\bfj$ gives $x=x'$.
\end{proof}

\subsection{The \texorpdfstring{$g$}{g}-vector parametrizations}\label{subsec:g-vector-parametrizations}

Fix $v\in\weyl$ and
$\bfk=(k_1,\dots,k_t)\in\rR(v)$.
We extend scalars from $\Q(q)$ to $\Q(q^{1/2})$
for the quantum cluster algebra structure on $\qA_q(\n(v))$
\cite[Theorem~12.3]{GLS13},
\cite[Theorem~7.3]{GY21}.
We recall the Gei\ss--Leclerc--Schr\"oer (GLS) seed associated with
$\bfk$ in the conventions of \cite[Section~1.3]{JP25}.
For $p\in[1,t]$, set
\begin{equation}\label{eq:p-plus-vartheta}
p^+:=\min\bigl(\{t+1\}\cup
\{q\in[p+1,t]\mid k_q=k_p\}\bigr),
\end{equation}
and put $K_{\mathrm{uf}}:=\{p\in[1,t]\mid p^+\le t\}$.
Write $v_{\le p}=s_{k_1}\cdots s_{k_p}$ and set
\[
\qmD_p^{\bfk}:=\qmD(v_{\le p}\Lambda_{k_p},\Lambda_{k_p})
\qquad(p\in[1,t]).
\]
The minors $\qmD_p^{\bfk}$, $p\in[1,t]$, pairwise $q$-commute.
The exchange matrix $\widetilde B_{\bfk}$ is defined, for
$p\in[1,t]$ and $q\in K_{\mathrm{uf}}$, by
\[
(\widetilde B_{\bfk})_{p,q}=
\begin{cases}
1 & \text{if }p=q^+,\\
a_{k_p,k_q} & \text{if }q<p<q^+<p^+,\\
-1 & \text{if }p^+=q,\\
-a_{k_p,k_q} & \text{if }p<q<p^+<q^+,\\
0 & \text{otherwise}.
\end{cases}
\]
Then
\begin{equation}\label{eq:GLS-seed}
\bigl(\{\qmD_p^{\bfk}\}_{p\in[1,t]},
\widetilde B_{\bfk}\bigr)
\end{equation}
is the GLS seed of $\qA_q(\n(v))$ associated with $\bfk$.
The frozen variables of the GLS seed \eqref{eq:GLS-seed} are
$\qmD_p^{\bfk}$ with $p^+=t+1$.
For $p\in[1,t]$ with $p^+=t+1$, every reflection $s_{k_r}$ with $r>p$
fixes $\Lambda_{k_p}$, so $v\Lambda_{k_p}=v_{\le p}\Lambda_{k_p}$
and hence $\qmD_p^{\bfk}=\qmC_{k_p}^v$.
We define the left and right $g$-vector parametrizations of $\LB(v)$
with respect to the GLS seed \eqref{eq:GLS-seed}.

Following \cite[Definition~3.1.1]{Qin20}, we define the dominance order
on $\Z^{[1,t]}$ by
\begin{equation}\label{eq:dominance-order}
\mathbf a\preceq_{\bfk}\mathbf b
\quad\Longleftrightarrow\quad
\mathbf a=\mathbf b+\widetilde B_{\bfk}\mathbf n
\quad\text{for some }\mathbf n\in\Z_{\ge0}^{K_{\mathrm{uf}}}.
\end{equation}

The quantum Laurent polynomial ring for the GLS seed \eqref{eq:GLS-seed}
is obtained from the subalgebra generated by
$\qmD_1^{\bfk},\dots,\qmD_t^{\bfk}$ by adjoining
$(\qmD_1^{\bfk})^{-1},\dots,(\qmD_t^{\bfk})^{-1}$
\cite[Section~2.2]{Qin20}.
For $z\in\tGup(v)$, the Laurent expansion of $z$ in powers of
$\qmD_1^{\bfk},\dots,\qmD_t^{\bfk}$ is
\[
z=\sum_{\mathbf a\in\Z^{[1,t]}}
c_{\mathbf a}(q)(\qmD_1^{\bfk})^{a_1}\cdots
(\qmD_t^{\bfk})^{a_t},
\qquad c_{\mathbf a}(q)\in\Q(q).
\]
By \cite[Proposition~3.4.5 and Theorem~9.5.1]{Qin20}, the exponent
vectors with nonzero coefficients have a unique maximal element and
a unique minimal element for $\preceq_{\bfk}$.
The maximal and minimal elements of
$\{\mathbf a\mid c_{\mathbf a}(q)\ne0\}$ for $\preceq_{\bfk}$
are called the \emph{degree} and \emph{codegree} of $z$ with respect
to the GLS seed \eqref{eq:GLS-seed}, respectively
\cite[Definition~3.1.3]{Qin20}.
Following \cite[Section~1.3]{JP25}, we set
\begin{equation}\label{eq:g-vector-extrema}
\gL_{\bfk}(z):=
\max_{\preceq_{\bfk}}\{\mathbf a\mid c_{\mathbf a}(q)\ne0\},
\qquad
\gR_{\bfk}(z):=
\min_{\preceq_{\bfk}}\{\mathbf a\mid c_{\mathbf a}(q)\ne0\},
\end{equation}
and call $\gL_{\bfk}(z)$ and $\gR_{\bfk}(z)$ the left and right
$g$-vectors of $z$, respectively.
For $x\in\LB(v)$, we set, as in \cite[Section~3]{JP25},
\begin{equation}\label{eq:g-vector-definition}
\gL_{\bfk}(x):=\gL_{\bfk}(\tGup(x)),
\qquad
\gR_{\bfk}(x):=\gR_{\bfk}(\tGup(x)).
\end{equation}

\begin{samepage}
Following \cite[(3.3)]{JP25}, define the integer matrices
\begin{equation}\label{eq:MN-matrices}
N_{\bfk}=(n_{p,q})_{p,q\in [1,t]},
\qquad
M_{\bfk}=(m_{p,q})_{p,q\in [1,t]}
\end{equation}
by
\begin{equation}\label{eq:MN-entries}
n_{p,q}=
\begin{cases}
1 & \text{if } p=q,\\
-1 & \text{if } p^+=q,\\
0 & \text{otherwise},
\end{cases}
\qquad
m_{p,q}=
\begin{cases}
\langle h_{k_p},s_{k_{p+1}}\cdots s_{k_q}\Lambda_{k_q}\rangle & \text{if } p\le q,\\
0 & \text{if } p>q.
\end{cases}
\end{equation}
\end{samepage}
The matrices $N_{\bfk}$ and $M_{\bfk}$ are upper triangular with
diagonal entries equal to $1$.

For $x\in\LB(v)$, \cite[Proposition~3.1\textup{(ii)--(iii)}]{JP25} gives
\begin{equation}\label{eq:g-vectors-explicit}
\gR_{\bfk}(x)=N_{\bfk}\PBW_{\bfk}(x),
\qquad
\gL_{\bfk}(x)=M_{\bfk}^{-1}\STR_{\bfk}(x).
\end{equation}

Since $N_{\bfk}$ and $M_{\bfk}$ are invertible over $\Z$ and
$\PBW_{\bfk}$ and $\STR_{\bfk}$ are bijective,
\eqref{eq:g-vectors-explicit} gives bijections

\begin{equation}\label{eq:four-param-g-vectors}
\gR_{\bfk}\colon \LB(v)\xrightarrow{\sim}\Z^{[1,t]},
\qquad
\gL_{\bfk}\colon \LB(v)\xrightarrow{\sim}\Z^{[1,t]}.
\end{equation}

As in \cite[(3.1) and Proposition~3.1\textup{(ii)--(iii)}]{JP25},
we define the transition maps among $\PBW_{\bfk}$, $\STR_{\bfk}$,
$\gL_{\bfk}$, and $\gR_{\bfk}$ by
\begin{equation}\label{eq:psi-def}
\psi_{\bfk}=\strS_{\bfk}\circ \pbwP_{\bfk}^{-1},
\qquad
\gM_{\bfk}=\strS_{\bfk}\circ (\gL_{\bfk})^{-1},
\qquad
\gN_{\bfk}=\gR_{\bfk}\circ \pbwP_{\bfk}^{-1}.
\end{equation}
In general, $\psi_{\bfk}$ is not $\Z$-linear
\cite[Remark~3.4\textup{(ii)}]{JP25}.
The transition maps $\gM_{\bfk}$ and $\gN_{\bfk}$ satisfy
\begin{equation}\label{eq:MN-action}
\gM_{\bfk}(\mathbf u)=M_{\bfk}\cdot \mathbf u
\qquad (\mathbf u\in \Z^{[1,t]}),
\qquad
\gN_{\bfk}(\mathbf a)=N_{\bfk}\cdot \mathbf a
\qquad (\mathbf a\in \tcP_{\bfk}(v)).
\end{equation}
As in \cite[(3.5)]{JP25}, one has the commutative diagram
\begin{equation}\label{eq:isom-four-param}
\begin{aligned}
\xymatrix{
&&& \LB(v) \ar[ld]_{\pbwP_{\bfk}} \ar[rd]^{\strS_{\bfk}} \\
\Z^{[1,t]} && \ar[ll]_{\gN_{\bfk}} \tcP_{\bfk}(v) \ar[rr]^{\psi_{\bfk}} && \tcS_{\bfk}(v) && \ar[ll]_{\gM_{\bfk}} \Z^{[1,t]}.
}
\end{aligned}
\end{equation}

We apply the definitions
\eqref{eq:GLS-seed}, \eqref{eq:dominance-order},
\eqref{eq:g-vector-extrema}--\eqref{eq:g-vector-definition},
\eqref{eq:MN-matrices}--\eqref{eq:MN-entries}, and \eqref{eq:psi-def}
to $\LB_\vartheta(w)$ using the Cartan datum of $\g_\vartheta$
and the reduced expression $\bfj$, and to $\LB(\Theta(w))$
using the Cartan datum of $\g$ and the reduced expression $\bfi$.

\subsection{The four parametrizations}\label{subsec:four-parametrizations}

\begin{theorem}\label{thm:four-param-localized}\label{thm:PBW-localized}\label{thm:string-localized}\label{thm:gR-localized}\label{thm:gL-localized}
For every $x\in\LB_\vartheta(w)$, one has
\begin{enumerate}[label=\textup{(\roman*)},ref=\roman*]
\item\label{item:four-param-PBW}
$\PBW_{\bfi}\bigl(\widetilde{\Upsilon}_\vartheta(x)\bigr)
=\Gamma_{\bfj}\bigl(\PBW_{\bfj}(x)\bigr)$.
\item\label{item:four-param-STR}
$\STR_{\bfi}\bigl(\widetilde{\Upsilon}_\vartheta(x)\bigr)
=\Gamma_{\bfj}\bigl(\STR_{\bfj}(x)\bigr)$.
\item\label{item:four-param-gR}
$\gR_{\bfi}\bigl(\widetilde{\Upsilon}_\vartheta(x)\bigr)
=\Gamma_{\bfj}\bigl(\gR_{\bfj}(x)\bigr)$.
\item\label{item:four-param-gL}
$\gL_{\bfi}\bigl(\widetilde{\Upsilon}_\vartheta(x)\bigr)
=\Gamma_{\bfj}\bigl(\gL_{\bfj}(x)\bigr)$.
\end{enumerate}
\end{theorem}

The proofs of Theorem~\ref{thm:four-param-localized}\textup{(\ref{item:four-param-gR})--(\ref{item:four-param-gL})}
use Lemma~\ref{lem:N-intertwining}.
\begin{lemma}\label{lem:N-intertwining}\label{lem:M-intertwining}
One has
\[
N_{\bfi}\Gamma_{\bfj}=\Gamma_{\bfj} N_{\bfj},
\qquad
M_{\bfi}\Gamma_{\bfj}=\Gamma_{\bfj} M_{\bfj}.
\]
\end{lemma}

\begin{proof}

For $1\le k\le\ell$ and $1\le r\le d_{j_k}$, let $(k,r)$ denote
the $r$-th position in $D_k$.
For $\bfv=(v_1,\dots,v_\ell)\in\Z^{[1,\ell]}$,
set $v_{\ell+1}=0$.
The next occurrence of $j_k^{(r)}$ in $\bfi$ is $(k^+,r)$ when
$k^+\le\ell$, and there is no next occurrence when $k^+=\ell+1$.
Consequently,
\[
\bigl(N_{\bfi}\Gamma_{\bfj} \bfv\bigr)_{(k,r)}
=m_{j_k^{(r)}}(v_k-v_{k^+})
=\bigl(\Gamma_{\bfj} N_\bfj \bfv\bigr)_{(k,r)}.
\]

For $M_{\bfi}\Gamma_{\bfj}=\Gamma_{\bfj}M_{\bfj}$, fix
$1\le p,q\le\ell$ and $1\le r\le d_{j_p}$.
Equation~\eqref{eq:upsilon-formula} gives
\[
\bigl(M_{\bfi}\Gamma_{\bfj}\bigr)_{(p,r),q}
=\sum_{u=1}^{d_{j_q}}(M_{\bfi})_{(p,r),(q,u)}m_{j_q^{(u)}},
\qquad
\bigl(\Gamma_{\bfj} M_\bfj\bigr)_{(p,r),q}
=m_{j_p^{(r)}}(M_\bfj)_{p,q}.
\]
If $p>q$, \eqref{eq:MN-entries} gives
\[
(M_{\bfi})_{(p,r),(q,u)}=0\quad(1\le u\le d_{j_q}),
\qquad (M_{\bfj})_{p,q}=0.
\]
For $p=q$, orthogonality within $F_{j_p}$ gives
\[
(M_{\bfi})_{(p,r),(p,u)}
=\begin{cases}
\langle h_{j_p^{(r)}},\Lambda_{j_p^{(u)}}\rangle=\delta_{r,u},&r\le u,\\
0,&r>u.
\end{cases}
\]
Thus
\[
\bigl(M_{\bfi}\Gamma_{\bfj}\bigr)_{(p,r),p}
=\bigl(\Gamma_{\bfj}M_{\bfj}\bigr)_{(p,r),p}
=m_{j_p^{(r)}}.
\]

Suppose that $p<q$. By \eqref{eq:MN-entries},
\begin{equation}\label{eq:M-expanded-entry}
(M_{\bfi})_{(p,r),(q,u)}
=\Bigl\langle h_{j_p^{(r)}},
 s_{j_p^{(r+1)}}\cdots s_{j_p^{(d_{j_p})}}
 \Theta(s_{j_{p+1}}^\vartheta\cdots s_{j_{q-1}}^\vartheta)
 s_{j_q^{(1)}}\cdots s_{j_q^{(u)}}\Lambda_{j_q^{(u)}}\Bigr\rangle.
\end{equation}
For $r<t\le d_{j_p}$ and $\lambda\in P$, orthogonality gives
\begin{equation}\label{eq:M-orthogonal-coroot}
\langle h_{j_p^{(r)}},s_{j_p^{(t)}}\lambda\rangle
=\langle h_{j_p^{(r)}},\lambda\rangle.
\end{equation}
For $u<t\le d_{j_q}$, one has $s_{j_q^{(t)}}\Lambda_{j_q^{(u)}}=\Lambda_{j_q^{(u)}}$,
so
\begin{equation}\label{eq:M-fiber-weight}
s_{j_q^{(1)}}\cdots s_{j_q^{(u)}}\Lambda_{j_q^{(u)}}
=\Theta(s_{j_q}^\vartheta)\Lambda_{j_q^{(u)}}.
\end{equation}
Substituting \eqref{eq:M-orthogonal-coroot} and \eqref{eq:M-fiber-weight}
into \eqref{eq:M-expanded-entry} yields
\[
(M_{\bfi})_{(p,r),(q,u)}
=\left\langle h_{j_p^{(r)}},
 \Theta(s_{j_{p+1}}^\vartheta\cdots s_{j_q}^\vartheta)
 \Lambda_{j_q^{(u)}}\right\rangle.
\]
Since $\sum_{u=1}^{d_{j_q}}m_{j_q^{(u)}}\Lambda_{j_q^{(u)}}
=\iota(\Lambda_{j_q}^\vartheta)$, \eqref{eq:similarity-weyl-equivariance}
and \eqref{eq:similarity-coroot-pairing} give
\begin{align*}
\bigl(M_{\bfi}\Gamma_{\bfj}\bigr)_{(p,r),q}
&=\left\langle h_{j_p^{(r)}},
 \Theta(s_{j_{p+1}}^\vartheta\cdots s_{j_q}^\vartheta)
 \iota(\Lambda_{j_q}^\vartheta)\right\rangle\\
&=m_{j_p^{(r)}}\left\langle h_{j_p}^\vartheta,
 s_{j_{p+1}}^\vartheta\cdots s_{j_q}^\vartheta
 \Lambda_{j_q}^\vartheta\right\rangle\\
&=m_{j_p^{(r)}}(M_\bfj)_{p,q}
=\bigl(\Gamma_{\bfj} M_\bfj\bigr)_{(p,r),q}.
\end{align*}
 \end{proof}

\begin{proof}[Proof of Theorem~\ref{thm:four-param-localized}]
\textup{(\ref{item:four-param-PBW})} The PBW identity is \eqref{eq:PBW-localized}.

\textup{(\ref{item:four-param-STR})} Write $x=b\cdot\fz(\bfa)$, where
$b\in B_\vartheta(w)$ and $\bfa=(a_j)_{j\in J}\in\Z^{\oplus J}$.
By Proposition~\ref{prop:S-extends-to-localized},
\eqref{eq:y-a}, and \eqref{eq:localized-affine-string},
\begin{align}
\STR_{\bfi}\bigl(\widetilde{\Upsilon}_\vartheta(x)\bigr)
&=\STR_{\bfi}\bigl(\Upsilon_\vartheta(b)\bigr)
  +\sum_{j\in J}\sum_{i\in F_j}m_i a_j\,
    \STR_{\bfi}(\fz_i^{\Theta(w)})\notag\\
&=\STR_{\bfi}\bigl(\Upsilon_\vartheta(b)\bigr)
  +\sum_{j\in J}a_j\,
    \STR_{\bfi}\bigl(\Upsilon_\vartheta(\fz_j^w)\bigr)\label{eq:localized-string-frozen-sum}\\
&=\Gamma_{\bfj}\bigl(\STR_{\bfj}(b)\bigr)
  +\sum_{j\in J}a_j\,\Gamma_{\bfj}\bigl(\STR_{\bfj}(\fz_j^w)\bigr)\label{eq:localized-string-scaled-sum}\\
&=\Gamma_{\bfj}\bigl(\STR_{\bfj}(x)\bigr).\label{eq:localized-string-coordinate}
\end{align}
Equality~\eqref{eq:localized-string-frozen-sum} uses
Lemma~\ref{lem:frozen-image} and \eqref{eq:localized-affine-string}.
Equality~\eqref{eq:localized-string-scaled-sum} uses
Proposition~\ref{prop:string-Bw} for $b$ and each $\fz_j^w$.
Equality~\eqref{eq:localized-string-coordinate} uses the linearity of
$\Gamma_{\bfj}$ and \eqref{eq:localized-affine-string}.

\textup{(\ref{item:four-param-gR}), (\ref{item:four-param-gL})} Lemma~\ref{lem:N-intertwining} and the invertibility of
$M_{\bfi}$ and $M_{\bfj}$ give
\[
N_{\bfi}\Gamma_{\bfj}=\Gamma_{\bfj}N_{\bfj},\qquad
M_{\bfi}^{-1}\Gamma_{\bfj}=\Gamma_{\bfj}M_{\bfj}^{-1}.
\]
Using \textup{(\ref{item:four-param-PBW})}, \textup{(\ref{item:four-param-STR})}, and
\eqref{eq:g-vectors-explicit}, we compute
\begin{align*}
\gR_{\bfi}\bigl(\widetilde{\Upsilon}_\vartheta(x)\bigr)
&=N_{\bfi}\PBW_{\bfi}\bigl(\widetilde{\Upsilon}_\vartheta(x)\bigr)\\
&=N_{\bfi}\Gamma_{\bfj}\bigl(\PBW_{\bfj}(x)\bigr)\\
&=\Gamma_{\bfj}\bigl(N_{\bfj}\PBW_{\bfj}(x)\bigr)
=\Gamma_{\bfj}\bigl(\gR_{\bfj}(x)\bigr),\\
\gL_{\bfi}\bigl(\widetilde{\Upsilon}_\vartheta(x)\bigr)
&=M_{\bfi}^{-1}\STR_{\bfi}\bigl(\widetilde{\Upsilon}_\vartheta(x)\bigr)\\
&=M_{\bfi}^{-1}\Gamma_{\bfj}\bigl(\STR_{\bfj}(x)\bigr)\\
&=\Gamma_{\bfj}\bigl(M_{\bfj}^{-1}\STR_{\bfj}(x)\bigr)
=\Gamma_{\bfj}\bigl(\gL_{\bfj}(x)\bigr).
\end{align*}
\end{proof}

\begin{corollary}\label{cor:intertwining}
The linear map $\Gamma_{\bfj}$ satisfies
\begin{enumerate}[label=\textup{(\roman*)},ref=\roman*]
\item\label{item:intertwining-psi}
$\psi_{\bfi}\circ\Gamma_{\bfj}=\Gamma_{\bfj}\circ\psi_{\bfj}$.
\item\label{item:intertwining-N}
$\gN_{\bfi}\circ\Gamma_{\bfj}=\Gamma_{\bfj}\circ\gN_{\bfj}$.
\item\label{item:intertwining-M}
$\gM_{\bfi}\circ\Gamma_{\bfj}=\Gamma_{\bfj}\circ\gM_{\bfj}$.
\end{enumerate}
\end{corollary}

\begin{proof}
\textup{(\ref{item:intertwining-psi})} For $x\in\LB_\vartheta(w)$,
Theorem~\ref{thm:four-param-localized}\textup{(\ref{item:four-param-PBW})--(\ref{item:four-param-STR})} gives
\begin{align*}
\psi_{\bfi}\bigl(\Gamma_{\bfj}(\PBW_\bfj(x))\bigr)
&=\psi_{\bfi}\bigl(\PBW_{\bfi}
  (\widetilde{\Upsilon}_\vartheta(x))\bigr)\\
&=\STR_{\bfi}\bigl(\widetilde{\Upsilon}_\vartheta(x)\bigr)\\
&=\Gamma_{\bfj}\bigl(\STR_\bfj(x)\bigr)
=\Gamma_{\bfj}\bigl(\psi_\bfj(\PBW_\bfj(x))\bigr).
\end{align*}
Since $\pbwP_\bfj$ is bijective, \textup{(\ref{item:intertwining-psi})} holds.

\textup{(\ref{item:intertwining-N})} and \textup{(\ref{item:intertwining-M})} follow
from Lemma~\ref{lem:N-intertwining} and \eqref{eq:MN-action}.
\end{proof}

\subsection{Compatibility with crystal operators}\label{subsec:localized-crystal-structures}

We prove analogues of Proposition~\ref{prop:similarity-map} for the
crystal structures on localized crystals constructed by Kashiwara and
Nakashima \cite{KN25}.

For $i\in I$, let $\mathcal B_i=\{b_i(n)\mid n\in\Z\}$ be the
elementary crystal used in Lemma~\ref{lem:orthogonal-commute}.
In particular,
\[
\wt(b_i(n))=n\alpha_i,\qquad
\ep_i(b_i(n))=-n,\qquad \te_i b_i(n)=b_i(n+1),
\]
and $\ep_h(b_i(n))=-\infty$ for $h\ne i$.

The following lemma follows from the tensor product formulas
\eqref{eq:crystal-tensor-functions} and \eqref{eq:crystal-tensor-operators}.

\begin{lemma}\label{lem:elementary-tensor-orthogonal}
Let $t\in\Z_{>0}$, let $k_1,\ldots,k_t\in I$, and let
$i,i'\in I$ satisfy $i\ne i'$ and $a_{i,i'}=0$.
For $z\in\mathcal B_{k_1}\otimes\cdots\otimes\mathcal B_{k_t}$,
the following hold.
\begin{enumerate}[label=\textup{(\roman*)},ref=\roman*]
\item\label{item:elementary-tensor-functions}
If $X\in\{\te_{i'},\tf_{i'}\}$ and $Xz\ne0$, then
\[
\ep_i(Xz)=\ep_i(z),\qquad \varphi_i(Xz)=\varphi_i(z).
\]
\item\label{item:elementary-tensor-commute}
One has
\[
\te_i\te_{i'}z=\te_{i'}\te_i z,\qquad
\tf_i\tf_{i'}z=\tf_{i'}\tf_i z.
\]
\end{enumerate}
\end{lemma}

\begin{proof}
Write $z=b_{k_1}(n_1)\otimes\cdots\otimes b_{k_t}(n_t)$, with
$n_1,\ldots,n_t\in\Z$.

\textup{(\ref{item:elementary-tensor-functions})} Suppose that $Xz\ne0$.
By \eqref{eq:crystal-tensor-operators}, for some $s$ with $k_s=i'$,
$X$ changes $b_{k_s}(n_s)$ to $b_{k_s}(n_s+1)$ or
$b_{k_s}(n_s-1)$.
For every $n\in\Z$, orthogonality gives
\[
\ep_i(b_{i'}(n))=\varphi_i(b_{i'}(n))=-\infty,
\qquad \langle h_i,\wt(b_{i'}(n))\rangle=0.
\]
Thus changing $n_s$ preserves $\ep_i$, $\varphi_i$, and the
pairing with $h_i$ on $b_{k_s}(n_s)$.
Repeated application of \eqref{eq:crystal-tensor-functions} gives
$\ep_i(Xz)=\ep_i(z)$ and $\varphi_i(Xz)=\varphi_i(z)$.

\textup{(\ref{item:elementary-tensor-commute})} We induct on $t$.
For $t=1$, at least one of $i,i'$ differs from $k_1$, so
\textup{(\ref{item:elementary-tensor-commute})} holds.
Suppose that $t>1$. If either $i$ or $i'$
does not occur in $(k_1,\ldots,k_t)$, then
$\te_i\te_{i'}z=\te_{i'}\te_i z=0$ and
$\tf_i\tf_{i'}z=\tf_{i'}\tf_i z=0$.
Otherwise, all four operators $\te_i,\tf_i,\te_{i'},\tf_{i'}$
are nonzero on every tensor in
$\mathcal B_{k_1}\otimes\cdots\otimes\mathcal B_{k_t}$.
Set
\[
x=b_{k_1}(n_1)\otimes\cdots\otimes b_{k_{t-1}}(n_{t-1}),
\qquad y=b_{k_t}(n_t),
\]
so that $z=x\otimes y$.

By \textup{(\ref{item:elementary-tensor-functions})}, applying $\te_{i'}$
to $x$ or $y$, when nonzero, leaves
$\varphi_i(x)$ and $\ep_i(y)$ unchanged.
Hence the inequality $\varphi_i(x)\ge\ep_i(y)$, which determines
whether $\te_i$ acts on $x$ or $y$ in
\eqref{eq:crystal-tensor-operators}, is unchanged.
Since $a_{i',i}=0$, applying $\te_i$ likewise does not change
whether $\te_{i'}$ acts on $x$ or $y$.
If $\te_i$ and $\te_{i'}$ both act on $x$, the induction
hypothesis gives $(\te_i\te_{i'}x)\otimes y=(\te_{i'}\te_i x)\otimes y$.
Since $i\ne i'$, the operators $\te_i$ and $\te_{i'}$ cannot
both act nontrivially on $y=b_{k_t}(n_t)$.
If $\te_i$ acts on $x$ and $\te_{i'}$ on $y$, then
$\te_i\te_{i'}z=\te_{i'}\te_i z=(\te_i x)\otimes(\te_{i'}y)$.
If $\te_i$ acts on $y$ and $\te_{i'}$ on $x$, then
$\te_i\te_{i'}z=\te_{i'}\te_i z=(\te_{i'}x)\otimes(\te_i y)$.
This proves $\te_i\te_{i'}z=\te_{i'}\te_i z$.

For $\tf_i$ and $\tf_{i'}$,
\eqref{eq:crystal-tensor-operators} uses $>$ in place of $\ge$.
Part~\textup{(\ref{item:elementary-tensor-functions})} also preserves
$\varphi_i(x)>\ep_i(y)$ and $\varphi_{i'}(x)>\ep_{i'}(y)$,
so the induction on $t$ gives $\tf_i\tf_{i'}z=\tf_{i'}\tf_i z$.
\end{proof}

For $j\in J$, let $\mathcal B_j=\{b_j(n)\mid n\in\Z\}$ be the
elementary crystal of $\g_\vartheta$ with index $j$.

Using Kashiwara's map \cite[p.~185]{K96} on each factor
$\mathcal B_{j_k}$, we define $T$ by
\begin{equation}\label{eq:cellular-factor-maps}
\begin{aligned}
T\colon\mathcal B_{j_1}\otimes\cdots\otimes\mathcal B_{j_\ell}
&\longrightarrow\mathcal B_{i_1}\otimes\cdots\otimes\mathcal B_{i_L},\\
\bigotimes_{k=1}^{\ell} b_{j_k}(n_k)
&\longmapsto
\bigotimes_{k=1}^{\ell}\ \bigotimes_{r=L_{k-1}+1}^{L_k}
b_{i_r}(m_{i_r}n_k).
\end{aligned}
\end{equation}

\begin{lemma}
\label{lem:cellular-tensor-similarity}
For every $j\in J$, one has
\[
T\circ\te_j=\Bigl(\prod_{i\in F_j}\te_i^{m_i}\Bigr)\circ T,\qquad
T\circ\tf_j=\Bigl(\prod_{i\in F_j}\tf_i^{m_i}\Bigr)\circ T.
\]
For $z\in\mathcal B_{j_1}\otimes\cdots\otimes\mathcal B_{j_\ell}$
and $i\in F_j$, one also has
\begin{equation}\label{eq:cellular-tensor-functions-similarity}
\wt(T(z))=\iota(\wt(z)),\qquad
\ep_i(T(z))=m_i\ep_j(z),\qquad
\varphi_i(T(z))=m_i\varphi_j(z).
\end{equation}
\end{lemma}

\begin{proof}
Fix $j\in J$ and set $T(0)=0$.
For an arbitrary nonempty finite sequence $(j_1,\ldots,j_\ell)$ in
$J$, form $(i_1,\ldots,i_L)$ by replacing each $j_k$ with the
fixed ordered list of elements of $F_{j_k}$, and define $T$ by
\eqref{eq:cellular-factor-maps}.
For $z=\bigotimes_{k=1}^{\ell}b_{j_k}(n_k)$,
\eqref{eq:similarity-alpha-hat} gives
\[
\wt(T(z))
=\sum_{k=1}^{\ell}\sum_{i\in F_{j_k}}m_i n_k\alpha_i
=\iota\Bigl(\sum_{k=1}^{\ell}n_k\alpha_{j_k}^{\vartheta}\Bigr)
=\iota(\wt(z)).
\]
We prove $T\circ\te_j=(\prod_{i\in F_j}\te_i^{m_i})\circ T$
by induction on $\ell$.
The induction also proves $\ep_i\circ T=m_i\ep_j$ and
$\varphi_i\circ T=m_i\varphi_j$ for $i\in F_j$.
For $\ell=1$, let $n\in\Z$.
If $j_1=j$, orthogonality within $F_j$ gives
\[
T(\te_j b_j(n))
=\bigotimes_{i\in F_j}b_i(m_i(n+1))
=\Bigl(\prod_{i\in F_j}\te_i^{m_i}\Bigr)T(b_j(n)).
\]
Moreover, $\ep_i(T(b_j(n)))=-m_i n=m_i\ep_j(b_j(n))$ and
$\varphi_i(T(b_j(n)))=m_i n=m_i\varphi_j(b_j(n))$.
If $j_1\ne j$, then $\te_j b_{j_1}(n)=0$ and
$\te_iT(b_{j_1}(n))=0$ for $i\in F_j$.
The values of $\ep_j,\varphi_j$ on $b_{j_1}(n)$ and of
$\ep_i,\varphi_i$ on $T(b_{j_1}(n))$ are all $-\infty$,
with the convention $m_i(-\infty)=-\infty$.

Suppose $\ell>1$ and the induction hypothesis holds for sequences
of length $\ell-1$.
Let $T'$ be the map defined by \eqref{eq:cellular-factor-maps}
for $(j_1,\ldots,j_{\ell-1})$, with $T'(0)=0$.
For $x\in\mathcal B_{j_1}\otimes\cdots\otimes
\mathcal B_{j_{\ell-1}}$ and $n\in\Z$, set
\[
y=b_{j_\ell}(n),\qquad
\overline x=T'(x),\qquad
\overline y=\bigotimes_{r=L_{\ell-1}+1}^{L_\ell}b_{i_r}(m_{i_r}n).
\]
Then $T(x\otimes y)=\overline x\otimes\overline y$.
For $i\in F_j$, the induction hypothesis for $x$ and the
case $\ell=1$ for $y$ give
\[
\begin{aligned}
\ep_i(\overline x)&=m_i\ep_j(x),&
\varphi_i(\overline x)&=m_i\varphi_j(x),\\
\ep_i(\overline y)&=m_i\ep_j(y),&
\varphi_i(\overline y)&=m_i\varphi_j(y).
\end{aligned}
\]
By \eqref{eq:similarity-alpha-hat},
$\wt(\overline x)=\iota(\wt(x))$ and
$\wt(\overline y)=\iota(\wt(y))$.
Thus \eqref{eq:crystal-tensor-functions} and
\eqref{eq:similarity-coroot-pairing} give
\[
\begin{aligned}
\ep_i(\overline x\otimes\overline y)
&=\max\{\ep_i(\overline x),\ep_i(\overline y)
 -\langle h_i,\wt(\overline x)\rangle\}\\
&=m_i\max\{\ep_j(x),\ep_j(y)-\langle h_j^\vartheta,\wt(x)\rangle\}\\
&=m_i\ep_j(x\otimes y).
\end{aligned}
\]
The formula $\varphi_i=\ep_i+\langle h_i,\wt\rangle$ also gives
$\varphi_i(\overline x\otimes\overline y)
=m_i\varphi_j(x\otimes y)$.

We now apply \eqref{eq:crystal-tensor-operators} to $x\otimes y$
and $\overline x\otimes\overline y$.
If $j\notin\{j_1,\ldots,j_\ell\}$, then $\te_j(x\otimes y)=0$
and $\te_iT(x\otimes y)=0$ for every $i\in F_j$.
Assume that $j\in\{j_1,\ldots,j_\ell\}$, and fix $i\in F_j$.
Since $i\in\{i_1,\ldots,i_L\}$, we have
$\te_i(\overline x\otimes\overline y)\ne0$.
\begin{enumerate}[label=\textup{(\roman*)}]
\item If $\varphi_j(x)\ge\ep_j(y)$, then
$\te_j(x\otimes y)=(\te_jx)\otimes y$.
The inequality $\varphi_i(\overline x)\ge\ep_i(\overline y)$
implies that $\te_i$ acts on $\overline x$.
For $s\ge0$,
$\varphi_i(\te_i^s\overline x)=m_i\varphi_j(x)+s
\ge\ep_i(\overline y)$.
Hence all $m_i$ applications of $\te_i$ act on $\overline x$.
\item If $\varphi_j(x)<\ep_j(y)$, then
$\te_j(x\otimes y)=x\otimes\te_jy$.
If $\varphi_j(x)=-\infty$, all $m_i$ applications of
$\te_i$ act on $\overline y$.
Otherwise, $\ep_j(y)-\varphi_j(x)\ge1$, so
\[
\varphi_i(\overline x)=m_i\varphi_j(x)
<m_i\ep_j(y)-s=\ep_i(\te_i^s\overline y)
\qquad(0\le s<m_i).
\]
Hence all $m_i$ applications of $\te_i$ act on $\overline y$.
\end{enumerate}

By Lemma~\ref{lem:elementary-tensor-orthogonal}\textup{(\ref{item:elementary-tensor-functions})},
a nonzero application of $\te_{i'}$, with $i'\in F_j\setminus\{i\}$,
to $\overline x$ or $\overline y$ preserves
$\varphi_i(\overline x)$ and $\ep_i(\overline y)$.
Hence all operators $\te_i^{m_i}$, $i\in F_j$, act on
$\overline x$ in case \textup{(i)} and on $\overline y$ in
case \textup{(ii)}.
In case \textup{(i)}, the induction hypothesis gives
\[
\begin{aligned}
\Bigl(\prod_{i\in F_j}\te_i^{m_i}\Bigr)T(x\otimes y)
&=\left(\Bigl(\prod_{i\in F_j}\te_i^{m_i}\Bigr)T'(x)\right)
 \otimes\overline y\\
&=T'(\te_jx)\otimes\overline y
=T\bigl((\te_jx)\otimes y\bigr).
\end{aligned}
\]
In case \textup{(ii)}, $\ep_j(y)>-\infty$, so $j_\ell=j$.
The calculation for $\ell=1$ gives
\[
\begin{aligned}
\Bigl(\prod_{i\in F_j}\te_i^{m_i}\Bigr)T(x\otimes y)
&=\overline x\otimes
 \left(\Bigl(\prod_{i\in F_j}\te_i^{m_i}\Bigr)\overline y\right)\\
&=\overline x\otimes\bigotimes_{r=L_{\ell-1}+1}^{L_\ell}
 b_{i_r}(m_{i_r}(n+1))\\
&=T\bigl(x\otimes\te_jy\bigr).
\end{aligned}
\]
Thus $T(\te_j(x\otimes y))
=(\prod_{i\in F_j}\te_i^{m_i})T(x\otimes y)$, completing the induction.

Finally, for the fixed reduced expression $\bfj$ of $w$,
every index of $J$ occurs in $\bfj$, and every index of $I$
occurs in $\bfi$, since $\supp(w)=J$.
The definitions of $\mathcal B_j$ and $\mathcal B_i$, together with
\eqref{eq:crystal-tensor-operators}, imply that $\te_j,\tf_j$
are mutually inverse on the domain of $T$, as are
$\te_i,\tf_i$ on the codomain of $T$.
For $z\in\mathcal B_{j_1}\otimes\cdots\otimes\mathcal B_{j_\ell}$,
we therefore have
\[
T(z)=T(\te_j\tf_jz)
=\Bigl(\prod_{i\in F_j}\te_i^{m_i}\Bigr)T(\tf_jz).
\]
Hence $\bigl(\prod_{i\in F_j}\tf_i^{m_i}\bigr)T(z)=T(\tf_jz)$.
\end{proof}

For $v\in\weyl$ with $\supp(v)=I$ and
$\bfk=(k_1,\ldots,k_t)\in\rR(v)$, let
\[
\mathrm K_\bfk\colon \LB(v)\xrightarrow{\sim}
\mathcal B_{k_1}\otimes\cdots\otimes\mathcal B_{k_t}
\]
be the crystal isomorphism of \cite[Main Theorem~9.3]{KN25}.
Write $\widetilde E_h,\widetilde F_h$ and
$\ep_h^{\mathrm{loc}},\varphi_h^{\mathrm{loc}}$, $h\in I$,
for the crystal operators and functions on $\LB(v)$ defined in
\cite[Definition~6.11]{KN25}.\footnote{ In \cite{KN25}, Kashiwara and Nakashima write
$\ep_h,\varphi_h$. We add the superscript $\mathrm{loc}$ to
distinguish the functions on $\LB(v)$ from the functions on
$B(\infty)$.}

For $z\in\LB(v)$ and $h\in I$, the crystal isomorphism $\mathrm K_\bfk$ gives
\begin{equation}\label{eq:KN-K-functions}
\begin{aligned}
\wt(\mathrm K_\bfk(z))&=\wt(z),&
\ep_h(\mathrm K_\bfk(z))&=\ep_h^{\mathrm{loc}}(z),&
\varphi_h(\mathrm K_\bfk(z))&=\varphi_h^{\mathrm{loc}}(z).
\end{aligned}
\end{equation}
\begin{proposition}[{\cite[Theorem~3.7]{KKOP23};
\cite[Section~6.4]{KN25}; \cite[Proposition~1.43]{TW16}}]
\label{prop:KN-localization-K}
For $v\in\weyl$, there are maps
\[
\mathrm Q_v\colon B(\infty)\longrightarrow\LB(v)\cup\{0\}
\]
and bijections
$\psi_*\colon\LB(v)\xrightarrow{\sim}\LB(v^{-1})$
satisfying the following properties.
\begin{enumerate}[label=\textup{(\roman*)},ref=\roman*]
\item\label{item:KN-localization-involution}
For $b\in B(v)$, one has $\mathrm Q_v(b)=b$.
For $c\in B(\infty)$ with $\mathrm Q_v(c)\ne0$, one has
\begin{equation}\label{eq:KN-localization-involution}
\psi_*(\mathrm Q_v(c))=\mathrm Q_{v^{-1}}(c^*).
\end{equation}
\item\label{item:KN-localization-K}
Let $\bfk=(k_1,\ldots,k_t)\in\rR(v)$ and
$c\in B(\infty)$ with $\mathrm Q_v(c)\ne0$.
Set
\begin{equation}\label{eq:KN-K-extraction}
c_t=c,\qquad c_{p-1}=\pi_{k_p}(c_p),\qquad
d_p=\ep_{k_p}^*(c_p)\qquad(p=t,\ldots,1).
\end{equation}
Then
\begin{equation}\label{eq:KN-K-localization-formula}
\mathrm K_\bfk(\mathrm Q_v(c))
=b_{k_1}(-d_1)\otimes\cdots\otimes b_{k_t}(-d_t).
\end{equation}
\end{enumerate}
\end{proposition}

\begin{lemma}\label{lem:localized-K-similarity}
For every $x\in\LB_\vartheta(w)$, one has
\begin{equation}\label{eq:cellular-K-similarity}
\mathrm K_\bfi\bigl(\widetilde\Upsilon_\vartheta(x)\bigr)
=T\bigl(\mathrm K_\bfj(x)\bigr).
\end{equation}
\end{lemma}

\begin{proof}
We first prove \eqref{eq:cellular-K-similarity} for $b\in B_\vartheta(w)$.
Define
\[
\begin{aligned}
b_\ell&=b,&
b_{k-1}&=\pi_{j_k}(b_k),\qquad
c_k=\ep_{j_k}^*(b_k)
&&(k=\ell,\ldots,1),\\
\overline b_L&=\Upsilon_\vartheta(b),&
\overline b_{r-1}&=\pi_{i_r}(\overline b_r),\qquad
\overline c_r=\ep_{i_r}^*(\overline b_r)
&&(r=L,\ldots,1).
\end{aligned}
\]
By Theorem~\ref{thm:similarity-Bw},
$\Upsilon_\vartheta(b)\in B(\Theta(w))$.

Applying \eqref{eq:KN-K-localization-formula} to $b$ and
$\Upsilon_\vartheta(b)$ gives
\[
\mathrm K_\bfj(b)=\bigotimes_{k=1}^\ell b_{j_k}(-c_k),\qquad
\mathrm K_\bfi(\Upsilon_\vartheta(b))
=\bigotimes_{r=1}^L b_{i_r}(-\overline c_r).
\]

We prove $\overline b_{L_k}=\Upsilon_\vartheta(b_k)$ by descending
induction on $k$.
The initial values give $\overline b_{L_\ell}=\Upsilon_\vartheta(b_\ell)$.
Suppose that $\overline b_{L_k}=\Upsilon_\vartheta(b_k)$ for some
$1\le k\le\ell$.
For $r\in D_k$, the definition of $\overline b_r$ gives
\[
\overline b_r
=\pi_{i_{r+1}}\cdots\pi_{i_{L_k}}(\overline b_{L_k}).
\]
For $r<s\le L_k$, the indices $i_r,i_s$ are distinct elements of
$F_{j_k}$.
Lemma~\ref{lem:orthogonal-commute}\textup{(\ref{item:orthogonal-projections-reflections})} and
\eqref{eq:similarity-epsilon-fiber} imply
\[
\overline c_r
=\ep_{i_r}^*(\overline b_r)
=\ep_{i_r}^*(\overline b_{L_k})
=m_{i_r}\ep_{j_k}^*(b_k)
=m_{i_r}c_k.
\]
The recursion and \eqref{eq:saito-similarity-projection-mi} give
\[
\begin{aligned}
\overline b_{L_{k-1}}
&=\Bigl(\prod_{i\in F_{j_k}}\pi_i\Bigr)\overline b_{L_k}\\
&=\Upsilon_\vartheta(\pi_{j_k}(b_k))
=\Upsilon_\vartheta(b_{k-1}).
\end{aligned}
\]
Thus $\overline c_r=m_{i_r}c_k$ for $1\le k\le\ell$ and
$r\in D_k$, and
\eqref{eq:cellular-factor-maps} yields
\begin{equation}\label{eq:cellular-K-similarity-Bw}
\mathrm K_\bfi(\Upsilon_\vartheta(b))
=\bigotimes_{k=1}^\ell\ \bigotimes_{r\in D_k}b_{i_r}(-m_{i_r}c_k)
=T(\mathrm K_\bfj(b)).
\end{equation}

Now let $x=b\cdot\fz(\bfa)$, where $b\in B_\vartheta(w)$ and
$\bfa=(a_h)_{h\in J}\in\Z^{\oplus J}$. Write
\[
\mathrm K_\bfj(b)=\bigotimes_{k=1}^\ell b_{j_k}(n_k),\qquad
\mathrm K_\bfj(\fz_h^w)=\bigotimes_{k=1}^\ell b_{j_k}(n_{h,k})
\quad(h\in J).
\]
By \cite[Section~6.4]{KN25}, one has
\begin{equation}\label{eq:cellular-K-translation}
\mathrm K_\bfj(x)
=\bigotimes_{k=1}^\ell b_{j_k}\Bigl(n_k+\sum_{h\in J}a_h n_{h,k}\Bigr).
\end{equation}
For $h\in J$ and $r\in D_k$, \eqref{eq:cellular-K-similarity-Bw} identifies the
$r$-th factors of $\mathrm K_\bfi(\Upsilon_\vartheta(b))$ and
$\mathrm K_\bfi(\Upsilon_\vartheta(\fz_h^w))$ as
$b_{i_r}(m_{i_r}n_k)$ and $b_{i_r}(m_{i_r}n_{h,k})$, respectively.
By Lemma~\ref{lem:frozen-image},
$\Upsilon_\vartheta(\fz_h^w)=
\prod_{i\in F_h}(\fz_i^{\Theta(w)})^{m_i}$.
For $z\in\LB(\Theta(w))$, write
\[
\mathrm K_\bfi(z)=\bigotimes_{r=1}^L b_{i_r}(t_r),
\qquad t_1,\ldots,t_L\in\Z.
\]
Since $\mathrm K_\bfi(\one)=\bigotimes_{r=1}^L b_{i_r}(0)$,
\cite[Section~6.4]{KN25} gives, for $h\in J$,
\[
\mathrm K_\bfi\left(z\cdot
 \prod_{i\in F_h}(\fz_i^{\Theta(w)})^{m_i a_h}\right)
=\bigotimes_{k=1}^\ell\ \bigotimes_{r\in D_k}
 b_{i_r}(t_r+a_hm_{i_r}n_{h,k}).
\]
Equations~\eqref{eq:localized-S-formula} and
\eqref{eq:cellular-K-translation} give
\[
\begin{aligned}
\mathrm K_\bfi(\widetilde\Upsilon_\vartheta(x))
&=\bigotimes_{k=1}^\ell\ \bigotimes_{r\in D_k}
 b_{i_r}\Bigl(m_{i_r}n_k+\sum_{h\in J}a_h m_{i_r}n_{h,k}\Bigr)\\
&=T(\mathrm K_\bfj(x)).
\end{aligned}
\]
\end{proof}

For $v\in\weyl$, write $\widetilde E_h^*,\widetilde F_h^*$ and
$\ep_h^{\mathrm{loc},*},\varphi_h^{\mathrm{loc},*}$, $h\in I$,
for the operators and functions on $\LB(v)$ defined
in \cite[Definition~6.11]{KN25}.
We express $\widetilde E_h^*,\widetilde F_h^*$ and
$\ep_h^{\mathrm{loc},*},\varphi_h^{\mathrm{loc},*}$ using the
string parametrization.

\begin{lemma}[{cf.\ \cite[Main Theorem~9.3]{KN25}}]
\label{lem:KN-starred-string}
Let $v\in\weyl$ and $\bfk=(k_1,\ldots,k_t)\in\rR(v)$.
For $z\in\LB(v)$, write
$\STR_\bfk(z)=(a_1,\ldots,a_t)$.
The map
\begin{equation}\label{eq:cellular-starred-string-map}
\begin{aligned}
\mathrm K_\bfk^*\colon\LB(v)
&\longrightarrow\mathcal B_{k_t}\otimes\cdots\otimes\mathcal B_{k_1},\\
z&\longmapsto b_{k_t}(-a_t)\otimes\cdots\otimes b_{k_1}(-a_1)
\end{aligned}
\end{equation}
is a bijection satisfying
\begin{equation}\label{eq:cellular-starred-string-operators}
\mathrm K_\bfk^*(\widetilde E_h^*z)=\te_h\mathrm K_\bfk^*(z),\qquad
\mathrm K_\bfk^*(\widetilde F_h^*z)=\tf_h\mathrm K_\bfk^*(z)
\qquad(h\in I).
\end{equation}
Moreover, for $z\in\LB(v)$ and $h\in I$, one has
\begin{equation}\label{eq:cellular-starred-string-functions}
\begin{aligned}
\wt(\mathrm K_\bfk^*(z))&=\wt(z),\\
\ep_h(\mathrm K_\bfk^*(z))&=\ep_h^{\mathrm{loc},*}(z),&
\varphi_h(\mathrm K_\bfk^*(z))&=\varphi_h^{\mathrm{loc},*}(z).
\end{aligned}
\end{equation}
\end{lemma}

\begin{proof}
Set $\bfk^{\mathrm{rev}}=(k_t,\ldots,k_1)\in\rR(v^{-1})$.
We prove
$\mathrm K_\bfk^*=\mathrm K_{\bfk^{\mathrm{rev}}}\circ\psi_*$.
The bijection $\psi_*$ in Proposition~\ref{prop:KN-localization-K}
satisfies the following identities by
\cite[Theorem~5.9 and Definition~6.11]{KN25}:
\begin{equation}\label{eq:KN-starred-exchange}
\psi_*(\widetilde E_h^*z)=\widetilde E_h\psi_*(z),\qquad
\psi_*(\widetilde F_h^*z)=\widetilde F_h\psi_*(z).
\end{equation}
\begin{equation}\label{eq:KN-starred-exchange-functions}
\begin{aligned}
\wt(\psi_*(z))&=\wt(z),\\
\ep_h^{\mathrm{loc}}(\psi_*(z))&=\ep_h^{\mathrm{loc},*}(z),&
\varphi_h^{\mathrm{loc}}(\psi_*(z))&=\varphi_h^{\mathrm{loc},*}(z).
\end{aligned}
\end{equation}

Let $b\in B(v)$.
Since $\mathrm Q_v(b)=b$, \eqref{eq:KN-localization-involution} gives
\begin{equation}\label{eq:KN-starred-localized-input}
\mathrm Q_{v^{-1}}(b^*)=\psi_*(b)\ne0.
\end{equation}
Write $\STR_\bfk(b)=(a_1,\ldots,a_t)$ and define
\[
b_0=b,\qquad b_r=\te_{k_r}^{a_r}b_{r-1},\qquad
a_r=\ep_{k_r}(b_{r-1})\qquad(1\le r\le t).
\]
The identities $\te_h^*=*\te_h*$ and $\ep_h^*=\ep_h\circ *$
on $B(\infty)$ give
\[
b_0^*=b^*,\qquad
b_r^*=(\te_{k_r}^*)^{a_r}b_{r-1}^*,\qquad
a_r=\ep_{k_r}^*(b_{r-1}^*)\qquad(1\le r\le t).
\]
Apply \eqref{eq:KN-K-extraction} to $c=b^*$ and
$\bfk^{\mathrm{rev}}$, and denote the resulting elements and integers by
$c_p,d_p$.
The initial value is $c_t=b_0^*$.
For $p=t+1-r$, the equality
$\pi_{k_r}(b_{r-1}^*)=b_r^*$ gives, by induction on $r$,
\[
c_{t+1-r}=b_{r-1}^*,\qquad d_{t+1-r}=a_r
\qquad(1\le r\le t).
\]
Equations~\eqref{eq:KN-K-localization-formula} and
\eqref{eq:KN-starred-localized-input} therefore give
\begin{equation}\label{eq:KN-reversed-string-Bv}
\mathrm K_{\bfk^{\mathrm{rev}}}(\psi_*(b))
=b_{k_t}(-a_t)\otimes\cdots\otimes b_{k_1}(-a_1)
=\mathrm K_\bfk^*(b).
\end{equation}

For $\mathbf d=(d_h)_{h\in I}\in\Z^{\oplus I}$, set
$z=b\cdot\fz(\mathbf d)$, and write
$\STR_\bfk(\fz_h^v)=(a_{h,1},\ldots,a_{h,t})$.
By \cite[Theorem~5.9, Lemma~5.8 and \textup{(5.7)}]{KN25},
each $\psi_*(\fz_h^v)$ is a product of integer powers of
$\fz_s^{v^{-1}}$, $s\in I$, and
\[
\psi_*(z)=\psi_*(b)\cdot\prod_{h\in I}\psi_*(\fz_h^v)^{d_h}.
\]
Equation~\eqref{eq:KN-reversed-string-Bv}, applied to $\fz_h^v$, gives
\[
\mathrm K_{\bfk^{\mathrm{rev}}}(\psi_*(\fz_h^v))
=b_{k_t}(-a_{h,t})\otimes\cdots\otimes b_{k_1}(-a_{h,1}).
\]
For $u\in\LB(v^{-1})$, write
$\mathrm K_{\bfk^{\mathrm{rev}}}(u)
=b_{k_t}(n_t)\otimes\cdots\otimes b_{k_1}(n_1)$,
with $n_1,\ldots,n_t\in\Z$.
By \cite[Section~6.4]{KN25}, for $h\in I$ and $d\in\Z$, one has
\[
\mathrm K_{\bfk^{\mathrm{rev}}}
 \bigl(u\cdot\psi_*(\fz_h^v)^d\bigr)
=b_{k_t}(n_t-d a_{h,t})\otimes\cdots\otimes
 b_{k_1}(n_1-d a_{h,1}).
\]
Taking $u=\psi_*(b)$ and multiplying successively by
$\psi_*(\fz_h^v)^{d_h}$ for $h\in I$, we obtain
\[
\mathrm K_{\bfk^{\mathrm{rev}}}(\psi_*(z))
=b_{k_t}\Bigl(-a_t-\sum_{h\in I}d_h a_{h,t}\Bigr)
 \otimes\cdots\otimes
 b_{k_1}\Bigl(-a_1-\sum_{h\in I}d_h a_{h,1}\Bigr).
\]
By \eqref{eq:localized-affine-string},
$\STR_\bfk(z)=(a_r+\sum_{h\in I}d_h a_{h,r})_{r=1}^t$.
Therefore $\mathrm K_\bfk^*=\mathrm K_{\bfk^{\mathrm{rev}}}\circ\psi_*$.
The bijectivity of $\mathrm K_{\bfk^{\mathrm{rev}}}$ and
$\psi_*$ implies that $\mathrm K_\bfk^*$ is bijective.
Since $\mathrm K_{\bfk^{\mathrm{rev}}}$ is a crystal isomorphism by
\cite[Main Theorem~9.3]{KN25}, equations
\eqref{eq:KN-starred-exchange} and
\eqref{eq:KN-starred-exchange-functions} give
\eqref{eq:cellular-starred-string-operators} and
\eqref{eq:cellular-starred-string-functions}.
\end{proof}

The element $\one\in\LB_\vartheta(w)$ is the class of
$(\one,\mathbf0)$ and satisfies $\widetilde E_j\one\ne0$
for every $j\in J$.
The map $\widetilde\Upsilon_\vartheta$ defined in
Proposition~\ref{prop:S-extends-to-localized} satisfies the following
analogues of the identities in Proposition~\ref{prop:similarity-map}.
\begin{theorem}\label{thm:localized-crystal-similarity}
The map
\[
\widetilde\Upsilon_\vartheta\colon
\LB_\vartheta(w)\longrightarrow\LB(\Theta(w))
\]
is injective, satisfies $\widetilde\Upsilon_\vartheta(\one)=\one$,
and has the following properties.
For $x\in\LB_\vartheta(w)$, $j\in J$, and $i\in F_j$, one has
\begin{enumerate}[label=\textup{(\roman*)},ref=\roman*]
\item\label{item:localized-similarity-ordinary}
\begin{equation}\label{eq:localized-crystal-E-F}
\begin{aligned}
\widetilde\Upsilon_\vartheta(\widetilde E_jx)
&=\Bigl(\prod_{i\in F_j}\widetilde E_i^{m_i}\Bigr)
  \widetilde\Upsilon_\vartheta(x),\\
\widetilde\Upsilon_\vartheta(\widetilde F_jx)
&=\Bigl(\prod_{i\in F_j}\widetilde F_i^{m_i}\Bigr)
  \widetilde\Upsilon_\vartheta(x),
\end{aligned}
\end{equation}
\begin{equation}\label{eq:localized-crystal-functions}
\ep_i^{\mathrm{loc}}(\widetilde\Upsilon_\vartheta(x))
=m_i\ep_j^{\mathrm{loc}}(x),\qquad
\varphi_i^{\mathrm{loc}}(\widetilde\Upsilon_\vartheta(x))
=m_i\varphi_j^{\mathrm{loc}}(x).
\end{equation}
\item\label{item:localized-similarity-starred}
\begin{equation}\label{eq:localized-crystal-E-F-star}
\begin{aligned}
\widetilde\Upsilon_\vartheta(\widetilde E_j^*x)
&=\Bigl(\prod_{i\in F_j}(\widetilde E_i^*)^{m_i}\Bigr)
  \widetilde\Upsilon_\vartheta(x),\\
\widetilde\Upsilon_\vartheta(\widetilde F_j^*x)
&=\Bigl(\prod_{i\in F_j}(\widetilde F_i^*)^{m_i}\Bigr)
  \widetilde\Upsilon_\vartheta(x),
\end{aligned}
\end{equation}
\begin{equation}\label{eq:localized-crystal-functions-star}
\ep_i^{\mathrm{loc},*}(\widetilde\Upsilon_\vartheta(x))
=m_i\ep_j^{\mathrm{loc},*}(x),\qquad
\varphi_i^{\mathrm{loc},*}(\widetilde\Upsilon_\vartheta(x))
=m_i\varphi_j^{\mathrm{loc},*}(x).
\end{equation}
\item\label{item:localized-similarity-weight}
\begin{equation}\label{eq:localized-crystal-weight}
\wt(\widetilde\Upsilon_\vartheta(x))=\iota(\wt(x)).
\end{equation}
\end{enumerate}
The products in \eqref{eq:localized-crystal-E-F} and
\eqref{eq:localized-crystal-E-F-star} are independent of the order
of their factors.
Among maps from $\LB_\vartheta(w)$ to $\LB(\Theta(w))$,
$\widetilde\Upsilon_\vartheta$ is uniquely determined by
$\widetilde\Upsilon_\vartheta(\one)=\one$ and
\eqref{eq:localized-crystal-E-F}.
\end{theorem}

\begin{proof}
The injectivity and the equality
$\widetilde\Upsilon_\vartheta(\one)=\one$ follow from
Proposition~\ref{prop:S-extends-to-localized}.

\textup{(\ref{item:localized-similarity-ordinary})}
Since $\mathrm K_\bfj$ and $\mathrm K_\bfi$ are crystal isomorphisms,
\eqref{eq:cellular-K-similarity} and
Lemma~\ref{lem:cellular-tensor-similarity} give
\[
\begin{aligned}
\widetilde\Upsilon_\vartheta(\widetilde E_jx)
&=\mathrm K_\bfi^{-1}\bigl(T(\te_j\mathrm K_\bfj(x))\bigr)
=\Bigl(\prod_{i\in F_j}\widetilde E_i^{m_i}\Bigr)
  \widetilde\Upsilon_\vartheta(x),\\
\widetilde\Upsilon_\vartheta(\widetilde F_jx)
&=\mathrm K_\bfi^{-1}\bigl(T(\tf_j\mathrm K_\bfj(x))\bigr)
=\Bigl(\prod_{i\in F_j}\widetilde F_i^{m_i}\Bigr)
  \widetilde\Upsilon_\vartheta(x).
\end{aligned}
\]
Equations~\eqref{eq:KN-K-functions},
\eqref{eq:cellular-K-similarity}, and
\eqref{eq:cellular-tensor-functions-similarity} give
\[
\begin{aligned}
\ep_i^{\mathrm{loc}}(\widetilde\Upsilon_\vartheta(x))
&=\ep_i(T(\mathrm K_\bfj(x)))
=m_i\ep_j(\mathrm K_\bfj(x))
=m_i\ep_j^{\mathrm{loc}}(x),\\
\varphi_i^{\mathrm{loc}}(\widetilde\Upsilon_\vartheta(x))
&=\varphi_i(T(\mathrm K_\bfj(x)))
=m_i\varphi_j(\mathrm K_\bfj(x))
=m_i\varphi_j^{\mathrm{loc}}(x).
\end{aligned}
\]

\textup{(\ref{item:localized-similarity-starred})}
Let $T^{\mathrm{rev}}$ be the map defined by
\eqref{eq:cellular-factor-maps} for the sequence
$(j_\ell,\ldots,j_1)$, with the order on each fiber reversed.
Thus $T^{\mathrm{rev}}$ maps
$\mathcal B_{j_\ell}\otimes\cdots\otimes\mathcal B_{j_1}$ to
$\mathcal B_{i_L}\otimes\cdots\otimes\mathcal B_{i_1}$.
For $x\in\LB_\vartheta(w)$, write
$\STR_\bfj(x)=(a_1,\ldots,a_\ell)$.
Theorem~\ref{thm:four-param-localized}\textup{(\ref{item:four-param-STR})} gives
\[
\STR_\bfi(\widetilde\Upsilon_\vartheta(x))
=\Gamma_\bfj(a_1,\ldots,a_\ell).
\]
By \eqref{eq:cellular-starred-string-map} and the definition of
$\Gamma_\bfj$,
\begin{equation}\label{eq:cellular-K-star-similarity}
\mathrm K_\bfi^*(\widetilde\Upsilon_\vartheta(x))
=T^{\mathrm{rev}}(\mathrm K_\bfj^*(x)).
\end{equation}

Lemma~\ref{lem:cellular-tensor-similarity} gives
\begin{equation}\label{eq:cellular-reverse-tensor-operators}
T^{\mathrm{rev}}\circ\te_j
=\Bigl(\prod_{i\in F_j}\te_i^{m_i}\Bigr)\circ T^{\mathrm{rev}},\qquad
T^{\mathrm{rev}}\circ\tf_j
=\Bigl(\prod_{i\in F_j}\tf_i^{m_i}\Bigr)\circ T^{\mathrm{rev}}.
\end{equation}
Equations~\eqref{eq:cellular-starred-string-operators},
\eqref{eq:cellular-K-star-similarity}, and
\eqref{eq:cellular-reverse-tensor-operators} give
\[
\begin{aligned}
\widetilde\Upsilon_\vartheta(\widetilde E_j^*x)
&=(\mathrm K_\bfi^*)^{-1}
  \bigl(T^{\mathrm{rev}}(\te_j\mathrm K_\bfj^*(x))\bigr)\\
&=\Bigl(\prod_{i\in F_j}(\widetilde E_i^*)^{m_i}\Bigr)
  \widetilde\Upsilon_\vartheta(x),\\
\widetilde\Upsilon_\vartheta(\widetilde F_j^*x)
&=(\mathrm K_\bfi^*)^{-1}
  \bigl(T^{\mathrm{rev}}(\tf_j\mathrm K_\bfj^*(x))\bigr)\\
&=\Bigl(\prod_{i\in F_j}(\widetilde F_i^*)^{m_i}\Bigr)
  \widetilde\Upsilon_\vartheta(x).
\end{aligned}
\]
Equations~\eqref{eq:cellular-starred-string-functions},
\eqref{eq:cellular-K-star-similarity}, and
\eqref{eq:cellular-tensor-functions-similarity}, applied to
$T^{\mathrm{rev}}$, give
\[
\begin{aligned}
\ep_i^{\mathrm{loc},*}(\widetilde\Upsilon_\vartheta(x))
&=\ep_i(T^{\mathrm{rev}}(\mathrm K_\bfj^*(x)))\\
&=m_i\ep_j(\mathrm K_\bfj^*(x))
=m_i\ep_j^{\mathrm{loc},*}(x),\\
\varphi_i^{\mathrm{loc},*}(\widetilde\Upsilon_\vartheta(x))
&=\varphi_i(T^{\mathrm{rev}}(\mathrm K_\bfj^*(x)))\\
&=m_i\varphi_j(\mathrm K_\bfj^*(x))
=m_i\varphi_j^{\mathrm{loc},*}(x).
\end{aligned}
\]

\textup{(\ref{item:localized-similarity-weight})}
Equations~\eqref{eq:KN-K-functions},
\eqref{eq:cellular-K-similarity}, and
\eqref{eq:cellular-tensor-functions-similarity} give
\[
\begin{aligned}
\wt(\widetilde\Upsilon_\vartheta(x))
&=\wt(T(\mathrm K_\bfj(x)))\\
&=\iota(\wt(\mathrm K_\bfj(x)))
=\iota(\wt(x)).
\end{aligned}
\]

The crystal isomorphism $\mathrm K_\bfi$,
\eqref{eq:cellular-starred-string-operators}, and
Lemma~\ref{lem:elementary-tensor-orthogonal}\textup{(\ref{item:elementary-tensor-commute})}
give the independence of the order of the products in
\eqref{eq:localized-crystal-E-F} and
\eqref{eq:localized-crystal-E-F-star}.

By \cite[Theorem~9.4]{KN25}, the crystal $\LB_\vartheta(w)$ is
connected. Thus every $x\in\LB_\vartheta(w)$ is obtained from
$\one$ by finitely many applications of
$\widetilde E_j$ and $\widetilde F_j$, $j\in J$.
Starting with $\widetilde\Upsilon_\vartheta(\one)=\one$,
repeated use of \eqref{eq:localized-crystal-E-F} therefore determines
$\widetilde\Upsilon_\vartheta(x)$ uniquely.
\end{proof}

\section{Similarity maps and quantum twists}\label{sec:similarity-twist}

\subsection{Compatibility with quantum twists}

For $w\in\weyl$, Kimura and Oya constructed the quantum twist automorphism
\begin{equation}\label{eq:eta-w}
\eta_w\colon \widetilde{\qA}_q(\n(w))
\xrightarrow{\sim}
\widetilde{\qA}_q(\n(w))
\end{equation}
and proved that $\eta_w$ preserves the upper global basis
$\tGup(w)$ of $\widetilde{\qA}_q(\n(w))$; see
\cite[Theorems~4.13 and~6.1]{KiOy21}.
For $\bfi=(i_1,\ldots,i_\ell)\in\rR(w)$, the quantum twist
$\eta_w$ satisfies
\begin{equation}\label{eq:twist-GLS-minor}
\eta_w(\qmD_p^{\bfi})
\equiv_q (\qmC_{i_p}^w)^{-1}
\qmD(w\Lambda_{i_p},w_{\le p}\Lambda_{i_p})
\qquad(1\le p\le\ell);
\end{equation}
see \cite[Theorems~4.13 and~6.1]{KiOy21}.
By \cite[Theorem~6.1]{KiOy21},
\[
\eta_w\bigl((\qmC_i^w)^{-1}\bigr)
\equiv_q\qmC_i^w
\qquad(i\in I).
\]
Via the bijection $\tGup\colon\LB(w)\to\tGup(w)$, the inverse automorphism
$\eta_w^{-1}$ induces the permutation
\begin{equation}\label{eq:crystal-twist}
\cdual_w\colon \LB(w)\xrightarrow{\sim}\LB(w),
\qquad
\tGup(\cdual_w(x))=\eta_w^{-1}(\tGup(x)).
\end{equation}
For $u\in\weyl_\vartheta$, we write
$\eta_u$ for the quantum twist automorphism for $\g_\vartheta$
and $u$, and $\cdual_u$ for the permutation of
$\LB_\vartheta(u)$ characterized by \eqref{eq:crystal-twist}
with $\g_\vartheta$ in place of $\g$ and $u$ in place of $w$.
Theorem~\ref{thm:old-main} expresses the action of $\cdual_w$
in the PBW parametrization.

\begin{theorem}[{\cite[Theorem~3.3(iv)]{JP25}}]\label{thm:old-main}
Let $w\in\weyl$, let
$\bfi\in\rR(w)$, and let $x\in\LB(w)$. Then
$\displaystyle
\PBW_\bfi(\cdual_w(x))
=\psi_\bfi^{-1}\circ\gM_\bfi\circ\gN_\bfi
  \bigl(-\PBW_\bfi(x)\bigr)$.
The formula for $\PBW_\bfi(\cdual_w(x))$ also holds for $\g_\vartheta$,
with $J$ in place of $I$.
\end{theorem}

Fix $w\in\weyl_\vartheta$, a reduced expression
$\bfj\in\rR(w)$, and the reduced expression
$\bfi\in\rR(\Theta(w))$ defined by
\eqref{eq:expanded-reduced-expression}.

\begin{theorem}\label{thm:main-goal}
For every $x\in\LB_\vartheta(w)$, one has
\begin{equation}\label{eq:twist-and-similarity-main}
\widetilde{\Upsilon}_\vartheta(\cdual_w(x))
=
\cdual_{\Theta(w)}(\widetilde{\Upsilon}_\vartheta(x)).
\end{equation}
\end{theorem}

That is, the following diagram commutes:
\[
\xymatrix@C=5em@R=3em{
\LB_\vartheta(w)
  \ar[r]^{\widetilde{\Upsilon}_\vartheta}
  \ar[d]_{\cdual_w}
&
\LB(\Theta(w))
  \ar[d]^{\cdual_{\Theta(w)}}
\\
\LB_\vartheta(w)
  \ar[r]_{\widetilde{\Upsilon}_\vartheta}
&
\LB(\Theta(w)).
}
\]

\begin{proof}
By Theorem~\ref{thm:old-main} and
Theorem~\ref{thm:four-param-localized}\textup{(\ref{item:four-param-PBW})}, we have
\[
\PBW_{\bfi}\bigl(\widetilde{\Upsilon}_\vartheta(\cdual_w(x))\bigr)
=
\Gamma_{\bfj}\circ\psi_\bfj^{-1}\circ\gM_\bfj\circ\gN_\bfj
\bigl(-\PBW_\bfj(x)\bigr).
\]
Since $\psi_{\bfj}$ and $\psi_{\bfi}$ are bijections,
Corollary~\ref{cor:intertwining}\textup{(\ref{item:intertwining-psi})} also gives, on
$\tcS_{\bfj}(w)$,
\begin{equation}\label{eq:intertwine-psi-inverse}
\psi_{\bfi}^{-1}\circ\Gamma_{\bfj}=\Gamma_{\bfj}\circ\psi_{\bfj}^{-1}.
\end{equation}

Corollary~\ref{cor:intertwining}\textup{(\ref{item:intertwining-N})--(\ref{item:intertwining-M})},
\eqref{eq:intertwine-psi-inverse}, and the linearity of $\Gamma_{\bfj}$ give
\begin{align*}
\PBW_{\bfi}\bigl(\widetilde{\Upsilon}_\vartheta(\cdual_w(x))\bigr)
&=
\psi_{\bfi}^{-1}\circ\gM_{\bfi}\circ\gN_{\bfi}
\bigl(-\Gamma_{\bfj}(\PBW_\bfj(x))\bigr)\\
&=
\psi_{\bfi}^{-1}\circ\gM_{\bfi}\circ\gN_{\bfi}
\bigl(-\PBW_{\bfi}(\widetilde{\Upsilon}_\vartheta(x))\bigr).
\end{align*}
Applying Theorem~\ref{thm:old-main} to $\Theta(w)$ gives
\[
\psi_{\bfi}^{-1}\circ\gM_{\bfi}\circ\gN_{\bfi}
\bigl(-\PBW_{\bfi}(\widetilde{\Upsilon}_\vartheta(x))\bigr)
=\PBW_{\bfi}(\cdual_{\Theta(w)}(\widetilde{\Upsilon}_\vartheta(x))).
\]
The injectivity of $\pbwP_{\bfi}$ proves
\eqref{eq:twist-and-similarity-main}.
\end{proof}

\subsection{Equivalence and periods}\label{subsec:equivalence-periods}\leavevmode

Fix $w\in\weyl$.
Following \cite[Section~2]{JP25}, for $x,x'\in\LB(w)$, define
\begin{equation}\label{eq:minor-equivalence}
x\eqfr x'
\quad\Longleftrightarrow\quad
x=x'\cdot\fz(\bfa)
\quad\text{for some }\bfa\in\Z^{\oplus I}.
\end{equation}
For $\bfi\in\rR(w)$, we also write
$\STR_\bfi(x)\eqfr\STR_\bfi(x')$ when $x\eqfr x'$.

\begin{lemma}[{\cite[Theorem~6.1]{KiOy21}}]\label{lem:cdual-preserves-frozen-equivalence}
For $x\in\LB(w)$ and $\bfa\in\Z^{\oplus I}$, one has
\begin{equation}\label{eq:twist-minor-powers}
\cdual_w(x\cdot\fz(\bfa))
=\cdual_w(x)\cdot\fz(-\bfa).
\end{equation}
Consequently, $\cdual_w$ induces a permutation of
$\LB(w)/\!\eqfr$.
\end{lemma}

Following \cite[Section~5]{JP25}, define
\begin{equation}\label{eq:period-definition}
\PD(w):=\min\{m\in\Z_{>0}\mid
\cdual_w^m(x)\eqfr x\text{ for every }x\in\LB(w)\},
\end{equation}
with $\min\varnothing=\infty$.
By Lemma~\ref{lem:cdual-preserves-frozen-equivalence}, $\PD(w)$ is the
order of the permutation induced by $\cdual_w$ on $\LB(w)/\!\eqfr$.

For $u\in\weyl_\vartheta$, define the relation
$\eqfr$ on $\LB_\vartheta(u)$ and the period $\PD(u)$ by
\eqref{eq:minor-equivalence} and \eqref{eq:period-definition}, using
$\g_\vartheta,J,u$ in place of $\g,I,w$.
Lemma~\ref{lem:cdual-preserves-frozen-equivalence} applies to
the Cartan datum of $\g_\vartheta$.

\begin{proposition}\label{prop:similarity-frozen-equivalence}
Let $w\in\weyl_\vartheta$.
For $x,x'\in\LB_\vartheta(w)$, one has
\[
x\eqfr x'
\quad\Longleftrightarrow\quad
\widetilde{\Upsilon}_\vartheta(x)
\eqfr\widetilde{\Upsilon}_\vartheta(x').
\]
Thus $\widetilde{\Upsilon}_\vartheta$ induces an injective map
from $\LB_\vartheta(w)/\!\eqfr$ to
$\LB(\Theta(w))/\!\eqfr$, commuting with the permutations
induced by $\cdual_w$ and $\cdual_{\Theta(w)}$, respectively.
\end{proposition}

\begin{proof}

Write $x=b\cdot\fz(\mathbf d)$ with $b\in B_\vartheta(w)$ and
$\mathbf d\in\Z^{\oplus J}$.
For $\bfa\in\Z^{\oplus J}$, \eqref{eq:localized-S-formula} and
\eqref{eq:y-a} give
\begin{equation}\label{eq:localized-shift-compatibility}
\widetilde\Upsilon_\vartheta(x\cdot\fz(\bfa))
=\Upsilon_\vartheta(b)\cdot y(\mathbf d+\bfa)
=\widetilde\Upsilon_\vartheta(x)\cdot y(\bfa).
\end{equation}
If $x'=x\cdot\fz(\bfa)$, then \eqref{eq:localized-shift-compatibility} gives
$\widetilde\Upsilon_\vartheta(x')
\eqfr\widetilde\Upsilon_\vartheta(x)$.

Conversely, suppose that
\[
\widetilde\Upsilon_\vartheta(x')=
\widetilde\Upsilon_\vartheta(x)\cdot
\prod_{i\in I}(\fz_i^{\Theta(w)})^{c_i}
\qquad(c_i\in\Z).
\]
We will show that, for every $j\in J$, there is an integer $a_j$
such that $c_i=m_i a_j$ for all $i\in F_j$.
Choose $\bfj=(j_1,\dots,j_\ell)\in\rR(w)$ and let
$\bfi=(i_1,\dots,i_L)$ be defined by
\eqref{eq:expanded-reduced-expression}, with the blocks $D_k$ of
\eqref{eq:block-interval}. Put
\[
\bfv=(v_1,\dots,v_\ell):=\PBW_\bfj(x')-\PBW_\bfj(x).
\]
Theorem~\ref{thm:four-param-localized}\textup{(\ref{item:four-param-PBW})} and
\eqref{eq:localized-affine-pbw} give
\begin{equation}\label{eq:frozen-equivalence-pbw-difference}
\Gamma_{\bfj}(\bfv)
=\PBW_{\bfi}(\widetilde\Upsilon_\vartheta(x'))
 -\PBW_{\bfi}(\widetilde\Upsilon_\vartheta(x))
=\sum_{i\in I}c_i\PBW_{\bfi}(\fz_i^{\Theta(w)}).
\end{equation}
By \cite[Corollary~6.4]{Kimura12},
\[
\bigl(\PBW_{\bfi}(\fz_h^{\Theta(w)})\bigr)_r=\delta_{i_r,h}
\qquad(h\in I,\ 1\le r\le L).
\]
Taking the $r$-th component of
\eqref{eq:frozen-equivalence-pbw-difference} and using
\eqref{eq:upsilon-formula}, we obtain
\[
m_{i_r}v_k=\sum_{h\in I}c_h\delta_{i_r,h}=c_{i_r}
\qquad(1\le k\le\ell,\ r\in D_k).
\]
Since $(i_r)_{r\in D_k}$ lists the elements of $F_{j_k}$ once each,
\[
c_i=m_i v_k\qquad(1\le k\le\ell,\ i\in F_{j_k}).
\]
For each $j\in J$, the assumption $\supp(w)=J$ allows us to choose
$k$ with $j_k=j$ and set $a_j=v_k\in\Z$.
Then $\bfa=(a_j)_{j\in J}\in\Z^{\oplus J}$ satisfies
$c_i=m_i a_{\vartheta(i)}$ for all $i\in I$.
Hence, by \eqref{eq:y-a} and \eqref{eq:localized-shift-compatibility},
\[
\prod_{i\in I}(\fz_i^{\Theta(w)})^{c_i}=y(\bfa),\qquad
\widetilde\Upsilon_\vartheta(x')
=\widetilde\Upsilon_\vartheta(x)\cdot y(\bfa)
=\widetilde\Upsilon_\vartheta(x\cdot\fz(\bfa)).
\]
The injectivity in Proposition~\ref{prop:S-extends-to-localized}
therefore gives $x'=x\cdot\fz(\bfa)$.

The equivalence in Proposition~\ref{prop:similarity-frozen-equivalence}
shows that $\widetilde{\Upsilon}_\vartheta$ induces a well-defined
injective map from $\LB_\vartheta(w)/\!\eqfr$ to
$\LB(\Theta(w))/\!\eqfr$.
Theorem~\ref{thm:main-goal} and
Lemma~\ref{lem:cdual-preserves-frozen-equivalence} show that the map
induced by $\widetilde{\Upsilon}_\vartheta$ on
$\LB_\vartheta(w)/\!\eqfr$ commutes with the permutations induced
by $\cdual_w$ and $\cdual_{\Theta(w)}$.
\end{proof}

\begin{corollary}\label{prop:periodicity-descends-under-similarity}
Let $w\in\weyl_\vartheta$.
If $\PD(\Theta(w))$ is finite, then $\PD(w)$ is finite and
\[
\PD(w)\mid\PD(\Theta(w)).
\]
\end{corollary}

\begin{proof}
Set $m=\PD(\Theta(w))$.
For every $x\in\LB_\vartheta(w)$, Theorem~\ref{thm:main-goal} gives
\[
\widetilde{\Upsilon}_\vartheta(\cdual_w^m(x))
=
\cdual_{\Theta(w)}^m(\widetilde{\Upsilon}_\vartheta(x))
\eqfr\widetilde{\Upsilon}_\vartheta(x).
\]
Proposition~\ref{prop:similarity-frozen-equivalence} gives
$\cdual_w^m(x)\eqfr x$.
Hence the permutation induced by $\cdual_w$ on
$\LB_\vartheta(w)/\!\eqfr$ has $m$-th power equal to the
identity, so $\PD(w)$ is finite and divides $m$.
\end{proof}
\section{Periods in types \texorpdfstring{$B_n$}{Bn} and \texorpdfstring{$C_n$}{Cn}}\label{sec:application-BC-periods}

We prove Conjecture~1\textup{(i)(b)} of \cite[Section~5.2.1]{JP25}.
Fix $n\ge2$. Let $x_n^{B_n}$ and $x_n^{C_n}$ be the longest
elements among the minimal length representatives for
\[
\langle s_1,\ldots,s_{n-1}\rangle\backslash\weyl(B_n)
\qquad\text{and}\qquad
\langle s_1,\ldots,s_{n-1}\rangle\backslash\weyl(C_n),
\]
respectively. The subgroups $\langle s_1,\ldots,s_{n-1}\rangle$ are
the maximal standard parabolic subgroups obtained by deleting vertex $n$.

\begin{theorem}\label{prop:type-B-xn-period}\label{thm:type-C-xn-period}
For $n\ge2$, one has
\[
\PD(x_n^{B_n})=\PD(x_n^{C_n})=
\begin{cases}
2 & \text{if }n=2,\\[4pt]
4 & \text{if }n\ge3.
\end{cases}
\]
\end{theorem}

\subsection{Similarity maps in types \texorpdfstring{$B_n$}{Bn} and \texorpdfstring{$C_n$}{Cn}}

We first obtain upper bounds for $\PD(x_n^{B_n})$ and
$\PD(x_n^{C_n})$ from the type $A_{2n-1}$ period formula in
\cite[Theorem~5.18]{JP25}.

\label{ex:type-A-similarity-data}
Let $\g$ have type $A_{2n-1}$, and set
\[
I=\{1,\ldots,2n-1\},\qquad J=\{1,\ldots,n\},\qquad
\vartheta(a)=\min(a,2n-a).
\]
Then
\[
F_j=\{j,2n-j\}\quad(1\le j<n),\qquad F_n=\{n\}.
\]
Taking $m_i=1$ for every $i\in I$ gives $\g_\vartheta$ of type $B_n$
(Figure~\ref{fig:similarity-data}\textup{(a)} for $n=3$).
Taking $m_n=2$ and $m_i=1$ for $i\ne n$ gives $\g_\vartheta$ of type $C_n$
(Figure~\ref{fig:similarity-data}\textup{(b)} for $n=3$).
We write $x_n$ for $x_n^{B_n}$ or $x_n^{C_n}$ when the type is fixed.
For either choice of $(m_i)$, the embedding $\Theta$ of
Section~\ref{subsec:similarity-data} is given by
\[
\Theta(s_j^\vartheta)=s_js_{2n-j}\quad(1\le j<n),
\qquad
\Theta(s_n^\vartheta)=s_n.
\]

Let $x_n^{A_{2n-1}}$ be the longest element among the minimal length
representatives for
\[
\langle s_1,\ldots,s_{n-1},s_{n+1},\ldots,s_{2n-1}\rangle
\backslash\weyl(A_{2n-1}).
\]

By Lemma~\ref{lem:expanded-reduced-expression} and the commutation
of the reflections $s_i$, $i\in F_j$, the
embedding $\Theta$ sends the longest elements of
$\weyl_\vartheta$ and $\langle s_j^\vartheta\mid j\in J\setminus\{n\}\rangle$
to the longest elements of $\weyl(A_{2n-1})$ and
$\langle s_i\mid i\in I\setminus\{n\}\rangle$, respectively.
The longest minimal length representative is the product of the longest
element of the standard parabolic subgroup and the longest element of
the Weyl group. Hence
\begin{equation}\label{eq:vartheta-xn-type-A-identification}
\Theta(x_n)=x_n^{A_{2n-1}}.
\end{equation}

By \cite[Theorem~5.18]{JP25},
\begin{equation}\label{eq:type-A-xn-period}
\PD(x_n^{A_{2n-1}})=
\begin{cases}
2 & \text{if }n=2,\\[4pt]
4 & \text{if }n\ge3.
\end{cases}
\end{equation}
By \eqref{eq:vartheta-xn-type-A-identification} and
Corollary~\ref{prop:periodicity-descends-under-similarity},
\begin{equation}\label{eq:BC-period-divisibility}
\PD(x_n^{B_n})\mid\PD(x_n^{A_{2n-1}}),
\qquad
\PD(x_n^{C_n})\mid\PD(x_n^{A_{2n-1}}).
\end{equation}

\subsection{Exact periods}\leavevmode

We exclude $\PD(x_n)=1$ when $n=2$ and
$\PD(x_n)\in\{1,2\}$ when $n\ge3$ by considering $\tf_n\one$.
The elements $x_n^{B_n}$ and $x_n^{C_n}$ have the common reduced expression
\begin{equation}\label{eq:BC-standard-word}
\bfj_n=(n,n-1,\ldots,1,\;n,n-1,\ldots,2,\;\ldots,\;n);
\end{equation}
see \cite[Section~4.2]{JP25} for type $B_n$
(the Coxeter matrices of types $B_n$ and $C_n$ coincide).
Since $\bfj_n$ begins with $n$, one has
$\tf_n\one\in B_\vartheta(x_n)$.
Since $F_n=\{n\}$, \eqref{eq:similarity-e-f} and
\eqref{eq:localized-S-formula} give
\begin{equation}\label{eq:BC-test-element-image}
\widetilde\Upsilon_\vartheta(\tf_n\one)
=\tf_n^{m_n}\one,
\qquad
m_n=\begin{cases}
1 & \text{in type }B_n,\\
2 & \text{in type }C_n.
\end{cases}
\end{equation}
By Theorem~\ref{thm:main-goal} and
Proposition~\ref{prop:similarity-frozen-equivalence}, it suffices to show that
$\tf_n^a\one$ is not fixed modulo $\eqfr$ by
$\cdual_{x_n^{A_{2n-1}}}$ when $n=2$, and by
$\cdual_{x_n^{A_{2n-1}}}^{\,2}$ when $n\ge3$, for $a=1,2$.

We use the reduced expression
\begin{equation}\label{eq:JP25-type-A-word}
\bfi_n=
(n,n+1,\ldots,2n-1;\;
n-1,n,\ldots,2n-2;\;
\ldots;\;
1,2,\ldots,n)
\in\rR(x_n^{A_{2n-1}})
\end{equation}
from \cite[Equation~(4.5)]{JP25}.
We arrange the coordinates for $\bfi_n$ in an $n\times n$ array:
the entry at $(p,q)$, $1\le p,q\le n$, corresponds to the $((p-1)n+q)$-th letter
of $\bfi_n$, namely $n-p+q$.

For $1\le p,q\le n$, let $e_{p,q}$ denote the standard basis vector
corresponding to position $(p,q)$.
For $0\le p,q\le n$, set
\[
b_{p,q}:=\sum_{1\le i\le p,\ 1\le j\le q}e_{i,j},
\]
with empty sums equal to zero, as in \cite[Equation~(5.5)]{JP25}.

Since $\bfi_n$ begins with $n$, one has
$\tf_n^a\one\in B(x_n^{A_{2n-1}})$ for $a\in\{1,2\}$.

\begin{samepage}
\begin{lemma}\label{lem:type-A-first-nonreturn}

Let $a\in\{1,2\}$, and set
$y_a:=\tf_n^a\one\in\LB(x_n^{A_{2n-1}})$.\ 
If $n=2$, then

\[
\cdual_{x_2^{A_3}}(y_a)\not\equiv_{\mathrm{fr}} y_a.
\]
If $n\ge3$, then

\[
\cdual_{x_n^{A_{2n-1}}}^{\,2}(y_a)\not\equiv_{\mathrm{fr}} y_a.
\]
\end{lemma}
\end{samepage}

\begin{proof}

First take $a=1$. By \eqref{eq:STR-on-Bw},
$\STR_{\bfi_n}(y_1)=e_{1,1}$.
Applying \cite[Proposition~5.10\textup{(ii)}]{JP25} to $e_{1,1}$ gives
\[
\STR_{\bfi_n}\bigl(\cdual_{x_n^{A_{2n-1}}}(y_1)\bigr)
\eqfr b_{1,n}-e_{1,1}.
\]
The vector $b_{1,n}-e_{1,1}$ has entry $1$ at
$(1,2),\ldots,(1,n)$ and entry $0$ elsewhere.
Applying \cite[Proposition~5.10\textup{(i)}]{JP25} to $b_{1,n}-e_{1,1}$
and using Lemma~\ref{lem:cdual-preserves-frozen-equivalence}, we obtain
\[
\begin{aligned}
\STR_{\bfi_n}\bigl(\cdual_{x_n^{A_{2n-1}}}^{\,2}(y_1)\bigr)
&\eqfr\sum_{1\le p,q\le n-1}e_{p,q}\\
&=b_{n-1,n-1}.
\end{aligned}
\]

Define the linear map $\Delta\colon\Z^{n\times n}\to\Z$ by
\begin{equation}\label{eq:type-A-string-difference}
\Delta\bigl((v_{p,q})_{1\le p,q\le n}\bigr)
:=v_{1,1}-v_{1,2}-v_{2,1}+v_{2,2}.
\end{equation}
By \cite[Lemmas~5.8 and~5.9]{JP25}, the frozen elements
$\fz_i^{x_n^{A_{2n-1}}}$, $1\le i\le2n-1$, have string data
with respect to $\bfi_n$
\[
b_{n,q}\quad(1\le q\le n),\qquad
b_{p,n}\quad(1\le p<n).
\]
The definition of $\Delta$ gives $\Delta(b_{n,q})=0$ for
$1\le q\le n$ and $\Delta(b_{p,n})=0$ for $1\le p<n$.
For $x,x'\in\LB(x_n^{A_{2n-1}})$,
\eqref{eq:localized-affine-string} therefore gives
$\Delta(\STR_{\bfi_n}(x)-\STR_{\bfi_n}(x'))=0$ if $x\eqfr x'$.
Consequently,
\begin{equation}\label{eq:type-A-string-test}
\Delta\bigl(\STR_{\bfi_n}(x)-\STR_{\bfi_n}(x')\bigr)\ne0
\quad\Longrightarrow\quad x\not\equiv_{\mathrm{fr}}x'.
\end{equation}

When $n=2$, the array $b_{1,2}-2e_{1,1}$ is
\[
\begin{pmatrix}
-1&1\\
0&0
\end{pmatrix}.
\]
Hence $\Delta(b_{1,2}-2e_{1,1})=-1-1-0+0=-2$.
When $n\ge3$, the entries of $b_{n-1,n-1}-e_{1,1}$ in rows
$1,2$ and columns $1,2$ form the matrix
\[
\begin{pmatrix}
0&1\\
1&1
\end{pmatrix}.
\]
Hence $\Delta(b_{n-1,n-1}-e_{1,1})=0-1-1+1=-1$.

For $a=2$, take $\vartheta=\mathrm{id}_I$ and $m_i=2$ for every $i\in I$.
Then $A_\vartheta=A$, and we identify $B_\vartheta(\infty)$ with
$B(\infty)$ by matching $\one$ and the Kashiwara operators.
Proposition~\ref{prop:S-extends-to-localized} therefore gives a self-map
of $\LB(x_n^{A_{2n-1}})$ sending $y_1$ to $y_2$
by \eqref{eq:similarity-e-f} and \eqref{eq:localized-S-formula}.
The fibers of $\mathrm{id}_I$ are singletons, so the linear map defined
by \eqref{eq:upsilon-formula} is multiplication by $2$.
Theorems~\ref{thm:four-param-localized}\textup{(\ref{item:four-param-STR})} and~\ref{thm:main-goal} therefore give, for
$r=0,1,2$,
\[
\STR_{\bfi_n}\bigl(\cdual_{x_n^{A_{2n-1}}}^{\,r}(y_2)\bigr)
=2\,\STR_{\bfi_n}\bigl(\cdual_{x_n^{A_{2n-1}}}^{\,r}(y_1)\bigr).
\]
The linearity of $\Delta$ and the constancy of
$\Delta\circ\STR_{\bfi_n}$ on the equivalence classes for $\eqfr$
give, for $a\in\{1,2\}$,
\[
\begin{aligned}
\Delta\bigl(\STR_{\bfi_2}(\cdual_{x_2^{A_3}}(y_a))-
\STR_{\bfi_2}(y_a)\bigr)&=-2a\ne0 &&(n=2),\\
\Delta\bigl(\STR_{\bfi_n}(\cdual_{x_n^{A_{2n-1}}}^{\,2}(y_a))-
\STR_{\bfi_n}(y_a)\bigr)&=-a\ne0 &&(n\ge3).
\end{aligned}
\]
Equation~\eqref{eq:type-A-string-test} gives
$\cdual_{x_2^{A_3}}(y_a)\not\equiv_{\mathrm{fr}}y_a$ when $n=2$,
and $\cdual_{x_n^{A_{2n-1}}}^{\,2}(y_a)
\not\equiv_{\mathrm{fr}}y_a$ when $n\ge3$.
\end{proof}

\begin{proof}[Proof of Theorem~\ref{thm:type-C-xn-period}]
For $r=1,2$, \eqref{eq:BC-test-element-image},
Theorem~\ref{thm:main-goal}, and
Proposition~\ref{prop:similarity-frozen-equivalence} give
\begin{equation}\label{eq:BC-twist-test}
\cdual_{x_n}^{\,r}(\tf_n\one)\eqfr\tf_n\one
\quad\Longleftrightarrow\quad
\cdual_{x_n^{A_{2n-1}}}^{\,r}(\tf_n^{m_n}\one)
\eqfr\tf_n^{m_n}\one.
\end{equation}
By \eqref{eq:BC-twist-test}, Lemma~\ref{lem:type-A-first-nonreturn}
with $a=m_n$ excludes $\PD(x_n)=1$ when $n=2$ and
$\PD(x_n)\in\{1,2\}$ when $n\ge3$.
Equations~\eqref{eq:type-A-xn-period} and \eqref{eq:BC-period-divisibility}
therefore give
$\PD(x_n^{B_n})=\PD(x_n^{C_n})=2$ for $n=2$, and
$\PD(x_n^{B_n})=\PD(x_n^{C_n})=4$ for $n\ge3$.
\end{proof}

\begin{example}\label{ex:type-A-twist-iterates}
We illustrate the applications of \cite[Proposition~5.10]{JP25} in the
proof of Lemma~\ref{lem:type-A-first-nonreturn} for $n=2,3$ and $a=1$.
For each $n$, the grids labelled $k=0,1,2$ show the string data
$e_{1,1}$, $b_{1,n}-e_{1,1}$, and $b_{n-1,n-1}$
of representatives, modulo $\eqfr$, of
\[
\tf_n\one,\qquad
\cdual_{x_n^{A_{2n-1}}}(\tf_n\one),\qquad
\cdual_{x_n^{A_{2n-1}}}^{\,2}(\tf_n\one),
\]
respectively.
Each grid uses the coordinates for $\bfi_n$, with rows counted
from top to bottom and columns from left to right.
Shaded cells have entry $1$, and all other cells have entry $0$.

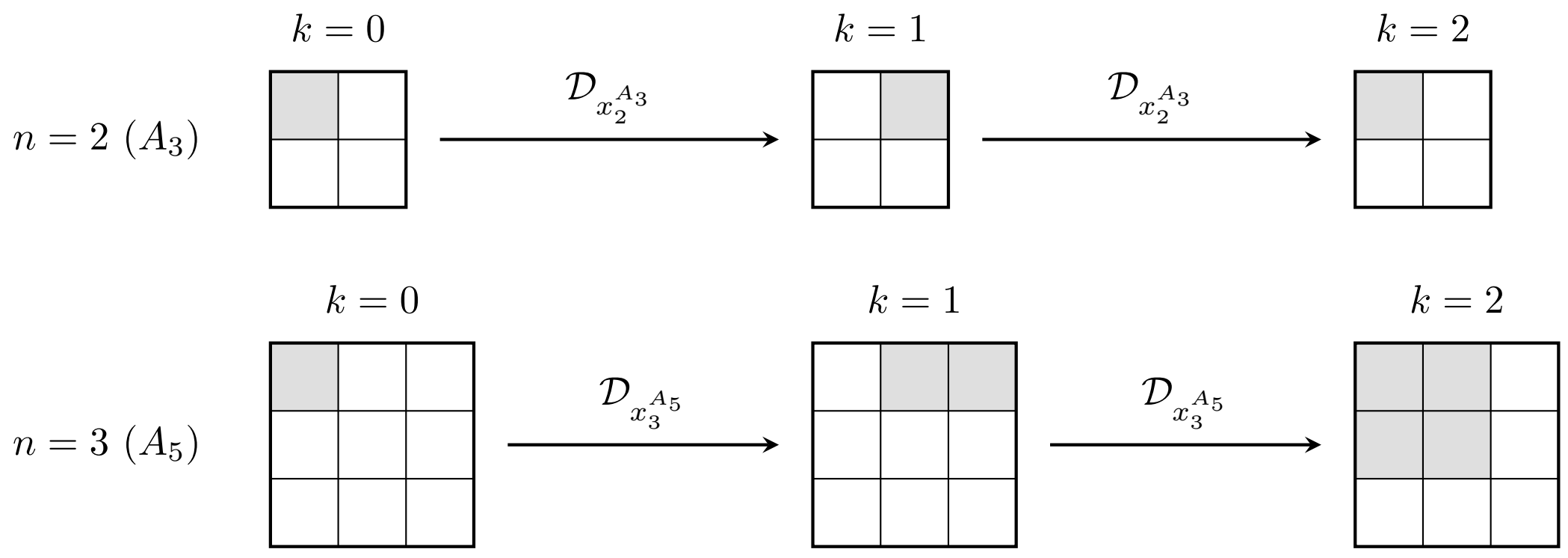
\begin{center}
\colorlet{typeAtwistexamplecolor}{.}
\begin{tikzpicture}[x=0.6cm,y=0.6cm,>=stealth,font=\small,
  draw=typeAtwistexamplecolor,text=typeAtwistexamplecolor]
  \node[anchor=east] at (-0.8,6) {$n=2$ ($A_3$)};
  \begin{scope}[shift={(0,5)}]
    \fill[gray!25] (0,1) rectangle (1,2);
    \draw[thick] (0,0) rectangle (2,2);
    \draw[thin] (1,0)--(1,2);
    \draw[thin] (0,1)--(2,1);
    \node at (1,2.65) {$k=0$};
  \end{scope}
  \draw[->,thick] (2.5,6)--(7.5,6)
    node[midway,above] {$\cdual_{x_2^{A_3}}$};
  \begin{scope}[shift={(8,5)}]
    \fill[gray!25] (1,1) rectangle (2,2);
    \draw[thick] (0,0) rectangle (2,2);
    \draw[thin] (1,0)--(1,2);
    \draw[thin] (0,1)--(2,1);
    \node at (1,2.65) {$k=1$};
  \end{scope}
  \draw[->,thick] (10.5,6)--(15.5,6)
    node[midway,above] {$\cdual_{x_2^{A_3}}$};
  \begin{scope}[shift={(16,5)}]
    \fill[gray!25] (0,1) rectangle (1,2);
    \draw[thick] (0,0) rectangle (2,2);
    \draw[thin] (1,0)--(1,2);
    \draw[thin] (0,1)--(2,1);
    \node at (1,2.65) {$k=2$};
  \end{scope}
  \node[anchor=east] at (-0.8,1.5) {$n=3$ ($A_5$)};
  \begin{scope}[shift={(0,0)}]
    \fill[gray!25] (0,2) rectangle (1,3);
    \draw[thick] (0,0) rectangle (3,3);
    \draw[thin] (1,0)--(1,3);
    \draw[thin] (0,1)--(3,1);
    \draw[thin] (2,0)--(2,3);
    \draw[thin] (0,2)--(3,2);
    \node at (1.5,3.65) {$k=0$};
  \end{scope}
  \draw[->,thick] (3.5,1.5)--(7.5,1.5)
    node[midway,above] {$\cdual_{x_3^{A_5}}$};
  \begin{scope}[shift={(8,0)}]
    \fill[gray!25] (1,2) rectangle (2,3);
    \fill[gray!25] (2,2) rectangle (3,3);
    \draw[thick] (0,0) rectangle (3,3);
    \draw[thin] (1,0)--(1,3);
    \draw[thin] (0,1)--(3,1);
    \draw[thin] (2,0)--(2,3);
    \draw[thin] (0,2)--(3,2);
    \node at (1.5,3.65) {$k=1$};
  \end{scope}
  \draw[->,thick] (11.5,1.5)--(15.5,1.5)
    node[midway,above] {$\cdual_{x_3^{A_5}}$};
  \begin{scope}[shift={(16,0)}]
    \fill[gray!25] (0,2) rectangle (1,3);
    \fill[gray!25] (1,2) rectangle (2,3);
    \fill[gray!25] (0,1) rectangle (1,2);
    \fill[gray!25] (1,1) rectangle (2,2);
    \draw[thick] (0,0) rectangle (3,3);
    \draw[thin] (1,0)--(1,3);
    \draw[thin] (0,1)--(3,1);
    \draw[thin] (2,0)--(2,3);
    \draw[thin] (0,2)--(3,2);
    \node at (1.5,3.65) {$k=2$};
  \end{scope}
\end{tikzpicture}
\end{center}

By \eqref{eq:type-A-string-difference},
\[
\bigl(\Delta(e_{1,1}),\Delta(b_{1,n}-e_{1,1}),
\Delta(b_{n-1,n-1})\bigr)
=\begin{cases}
(1,-1,1)&\text{if }n=2,\\
(1,-1,0)&\text{if }n=3.
\end{cases}
\]
Equation~\eqref{eq:type-A-string-test} therefore gives

\[
\cdual_{x_2^{A_3}}(\tf_2\one)\not\equiv_{\mathrm{fr}}\tf_2\one,
\qquad
\cdual_{x_3^{A_5}}^{\,2}(\tf_3\one)\not\equiv_{\mathrm{fr}}\tf_3\one.
\]
\end{example}

\bibliographystyle{amsplain}
\bibliography{ref}

@article{BFZ96,
  author  = {Berenstein, A. and Fomin, S. and Zelevinsky, A.},
  title   = {Parametrizations of canonical bases and totally positive matrices},
  journal = {Adv. Math.},
  year    = {1996},
  volume  = {122},
  number  = {1},
  pages   = {49--149},
}

@article{BZ97,
  author  = {Berenstein, A. and Zelevinsky, A.},
  title   = {Total positivity in {Schubert} varieties},
  journal = {Comment. Math. Helv.},
  year    = {1997},
  volume  = {72},
  number  = {1},
  pages   = {128--166},
}

@article{GLS13,
  author  = {Gei\ss, C. and Leclerc, B. and Schr\"oer, J.},
  title   = {Cluster structures on quantum coordinate rings},
  journal = {Selecta Math. (N.S.)},
  year    = {2013},
  volume  = {19},
  number  = {2},
  pages   = {337--397},
}

@book{HK02,
  author    = {Hong, J. and Kang, S.-J.},
  title     = {Introduction to quantum groups and crystal bases},
  year      = {2002},
  series    = {Graduate Studies in Mathematics},
  volume    = {42},
  publisher = {American Mathematical Society},
  address   = {Providence, RI},
}

@misc{IO21_arXiv2,
  author = {Ishibashi, T. and Oya, H.},
  title  = {Wilson lines and their {Laurent} positivity},
  year   = {2021},
  note   = {arXiv:2011.14260v2},
}

@article{IO21,
  author  = {Ishibashi, T. and Oya, H.},
  title   = {Wilson lines and their {Laurent} positivity},
  journal = {Math. Z.},
  year    = {2023},
  volume  = {305},
  number  = {2},
  pages   = {Paper No.~34, 60~pp.},
}

@article{KKKO18,
  author  = {Kang, S.-J. and Kashiwara, M. and Kim, M. and Oh, S.-j.},
  title   = {Monoidal categorification of cluster algebras},
  journal = {J. Amer. Math. Soc.},
  year    = {2018},
  volume  = {31},
  number  = {2},
  pages   = {349--426},
}

@article{K91,
  author  = {Kashiwara, M.},
  title   = {On crystal bases of the {$Q$}-analogue of universal enveloping algebras},
  journal = {Duke Math. J.},
  year    = {1991},
  volume  = {63},
  number  = {2},
  pages   = {465--516},
}

@article{K93,
  author  = {Kashiwara, M.},
  title   = {Global crystal bases of quantum groups},
  journal = {Duke Math. J.},
  year    = {1993},
  volume  = {69},
  number  = {2},
  pages   = {455--485},
}

@article{KKOP21C,
  author  = {Kashiwara, M. and Kim, M. and Oh, S.-j. and Park, E.},
  title   = {Localizations for quiver {Hecke} algebras},
  journal = {Pure Appl. Math. Q.},
  year    = {2021},
  volume  = {17},
  number  = {4},
  pages   = {1465--1548},
}

@article{KKOP23,
  author  = {Kashiwara, M. and Kim, M. and Oh, S.-j. and Park, E.},
  title   = {Localizations for quiver {Hecke} algebras {II}},
  journal = {Proc. Lond. Math. Soc. (3)},
  year    = {2023},
  volume  = {127},
  number  = {4},
  pages   = {1134--1184},
}

@article{Kimura12,
  author  = {Kimura, Y.},
  title   = {Quantum unipotent subgroup and dual canonical basis},
  journal = {Kyoto J. Math.},
  year    = {2012},
  volume  = {52},
  number  = {2},
  pages   = {277--331},
}

@article{KiOy21,
  author  = {Kimura, Y. and Oya, H.},
  title   = {Twist automorphisms on quantum unipotent cells and dual canonical bases},
  journal = {Int. Math. Res. Not. IMRN},
  year    = {2021},
  number  = {9},
  pages   = {6772--6847},
}

@book{LusztigBook,
  author    = {Lusztig, G.},
  title     = {Introduction to quantum groups},
  year      = {1993},
  series    = {Progress in Mathematics},
  volume    = {110},
  publisher = {Birkh\"auser Boston},
  address   = {Boston, MA},
}

@misc{Nakashima22,
  author = {Nakashima, T.},
  title  = {Categorified crystal structure on localized quantum coordinate rings},
  year   = {2022},
  note   = {arXiv:2208.08396},
}

@misc{Qin20,
  author = {Qin, F.},
  title  = {Dual canonical bases and quantum cluster algebras},
  year   = {2020},
  note   = {arXiv:2003.13674},
}

@misc{sage,
  author       = {{The Sage Developers}},
  title        = {{SageMath}, the {Sage Mathematics Software System}},
  year         = {2026},
  howpublished = {\url{https://www.sagemath.org}},
}

@article{GY21,
  author  = {Goodearl, K. R. and Yakimov, M. T.},
  title   = {Integral quantum cluster structures},
  journal = {Duke Math. J.},
  year    = {2021},
  volume  = {170},
  number  = {6},
  pages   = {1137--1200},
}

@misc{JP25,
  author = {Jung, W.-S. and Park, E.},
  title  = {Crystals and quantum twist automorphisms},
  year   = {2025},
  note   = {arXiv:2507.01306},
}

@article{K93D,
  author  = {Kashiwara, M.},
  title   = {The crystal base and {Littelmann's} refined {Demazure} character formula},
  journal = {Duke Math. J.},
  year    = {1993},
  volume  = {71},
  number  = {3},
  pages   = {839--858},
}

@incollection{K96,
  author    = {Kashiwara, M.},
  title     = {Similarity of crystal bases},
  booktitle = {{Lie} algebras and their representations (Seoul, 1995)},
  series    = {Contemp. Math.},
  volume    = {194},
  publisher = {Amer. Math. Soc.},
  address   = {Providence, RI},
  year      = {1996},
  pages     = {177--186},
}

@article{KN25,
  author  = {Kashiwara, M. and Nakashima, T.},
  title   = {Crystal structure of localized quantum unipotent coordinate category},
  journal = {Proc. Lond. Math. Soc. (3)},
  year    = {2026},
  volume  = {133},
  number  = {3},
  pages   = {Paper No.~e70218},
}

@article{K17,
  author  = {Kimura, Y.},
  title   = {Remarks on quantum unipotent subgroups and the dual canonical basis},
  journal = {Pacific J. Math.},
  year    = {2017},
  volume  = {286},
  number  = {1},
  pages   = {125--151},
}

@misc{LP26,
  author = {Lee, K.-H. and Park, E.},
  title  = {Imaginary elements and crystals for quantum groups},
  year   = {2026},
  note   = {arXiv:2610.32040},
}

@incollection{NS04,
  author    = {Naito, S. and Sagaki, D.},
  title     = {Crystal bases and diagram automorphisms},
  booktitle = {Representation theory of algebraic groups and quantum groups},
  series    = {Adv. Stud. Pure Math.},
  volume    = {40},
  publisher = {Math. Soc. Japan},
  address   = {Tokyo},
  year      = {2004},
  pages     = {321--341},
}

@article{Park23,
  author  = {Park, E.},
  title   = {Braid group action on extended crystals},
  journal = {Adv. Math.},
  year    = {2023},
  volume  = {429},
  pages   = {Paper No.~109193},
}

@article{Sai94,
  author  = {Saito, Y.},
  title   = {{PBW} basis of quantized universal enveloping algebras},
  journal = {Publ. Res. Inst. Math. Sci.},
  year    = {1994},
  volume  = {30},
  number  = {2},
  pages   = {209--232},
  note    = {doi:10.2977/PRIMS/1195166130},
}

@article{TW16,
  author  = {Tingley, P. and Webster, B.},
  title   = {{Mirkovi\'{c}--Vilonen} polytopes and {Khovanov--Lauda--Rouquier} algebras},
  journal = {Compos. Math.},
  year    = {2016},
  volume  = {152},
  number  = {8},
  pages   = {1648--1696},
}
\end{document}